\documentclass[a4paper,11pt]{amsart}
\usepackage{mathrsfs}
\usepackage{cases}
\usepackage{epic}
\usepackage{amsfonts}
\usepackage{graphicx}
\usepackage{amsmath}
\usepackage{amssymb, upgreek}
\usepackage{bm}
\usepackage{latexsym,todonotes}
\usepackage{pdflscape}
\usepackage[all]{xypic}
\usepackage{ytableau}
\usepackage[enableskew,vcentermath]{youngtab}
\usepackage[all]{xy}
\usepackage{color}
\usepackage{colordvi}
\usepackage{multicol}
\usepackage[linktocpage=true]{hyperref}
\hypersetup{colorlinks,linkcolor=blue,urlcolor=cyan,citecolor=blue}

\renewcommand{\ker}{\operatorname{Ker}\nolimits }
\allowdisplaybreaks

\begin{document}
\input xy
\xyoption{all}

\newtheorem{innercustomthm}{{\bf Theorem}}
\newenvironment{customthm}[1]
  {\renewcommand\theinnercustomthm{#1}\innercustomthm}
  {\endinnercustomthm}

  \newtheorem{innercustomcor}{{\bf Corollary}}
\newenvironment{customcor}[1]
  {\renewcommand\theinnercustomcor{#1}\innercustomcor}
  {\endinnercustomthm}

  \newtheorem{innercustomprop}{{\bf Proposition}}
\newenvironment{customprop}[1]
  {\renewcommand\theinnercustomprop{#1}\innercustomprop}
  {\endinnercustomthm}

\def \la{\lambda}
\newcommand{\LaK}{\la_{\texttt{Kr}}}
\newcommand{\Q}{\mathbb Q}
\newcommand{\QK}{Q_{\texttt{Kr}}}

\newcommand{\LaJ}{\La_{\texttt J}}
\newcommand{\tMHL}{\cs\cd\widetilde{\ch}(\LaK^\imath)}
\newcommand{\tMHg}{{}^\imath\widetilde{\ch}(\QJ)}
\newcommand{\PL}{\bbP^1_{\bfk}}
\def\bbP{\mathbb P}
\def\bfk{\Bbbk}
\renewcommand{\mod}{\operatorname{mod^{\rm nil}}\nolimits}
\newcommand{\Aut}{\operatorname{Aut}\nolimits}

\newcommand{\End}{\operatorname{End}\nolimits}
\newcommand{\Iso}{\operatorname{Iso}\nolimits}

\newcommand{\sqq}{{\bf v}}
\newcommand{\bq}{{q}}

\newcommand{\tMHLJ}{{}^\imath\widetilde{\ch}(\bfk \QJ)}

\newcommand{\bB}{{\bf B}}
\newcommand{\iH}{{H}^{\imath}}
\newcommand{\iP}{P^{\imath}}
\newcommand{\iB}{{B}^{\imath}}
\newcommand{\iQ}{{Q}^{\imath}}
\newcommand{\iV}{{V}^{\imath}}
\newcommand{\iVh}{\widehat{V}^{\imath}}
\newcommand{\is}{{s}^{\imath}}
\newcommand{\haE}{\widehat{E}}
\newcommand{\haH}{\widehat{H}}
\newcommand{\haP}{\widehat{P}}
\newcommand{\haQ}{\widehat{Q}}
\newcommand{\haV}{\widehat{V}}
\newcommand{\haT}{\widehat{\Theta}}
\newcommand{\ka}{\kappa}
\newcommand{\vv}{v}
\newcommand{\x}{{\bf x}}

\newcommand{\ov}{\overline}
\newcommand{\und}{\underline}
\newcommand{\tk}{\widetilde{k}}
\newcommand{\tK}{\widetilde{K}}
\newcommand{\tTT}{\operatorname{\widetilde{\texttt{\rm T}}}\nolimits}
\newcommand{\iRH}{\operatorname{{}^\imath \widetilde{\ch}}\nolimits}

\newcommand{\aut}{\operatorname{Aut}\nolimits}
\newcommand{\res}{\operatorname{res}\nolimits}

\newcommand{\QJ}{Q_{\texttt J}}
\def\bfk{\Bbbk}
\def\calc {\mathcal C}

\def \cI{\mathcal I}
\def \cJ{\mathcal J}

\newcommand{\rep}{\operatorname{rep}\nolimits}
\newcommand{\Ext}{\operatorname{Ext}\nolimits}
\newcommand{\Hom}{\operatorname{Hom}\nolimits}
\renewcommand{\Im}{\operatorname{Im}\nolimits}
\newcommand{\Ker}{\operatorname{Ker}\nolimits}
\newcommand{\Coh}{\operatorname{Coh}\nolimits}
\newcommand{\Id}{\operatorname{Id}\nolimits}

\newcommand{\coker}{\operatorname{Coker}\nolimits}
\newcommand{\qbinom}[2]{\begin{bmatrix} #1\\#2 \end{bmatrix} }

\newcommand{\gr}{\operatorname{gr}\nolimits}
\newcommand{\id}{\operatorname{Id}\nolimits}
\newcommand{\Res}{\operatorname{Res}\nolimits}
\def \tT{\widetilde{\mathcal T}}
\def \tTL{\tT(\la^\imath)}

\newcommand{\mbf}{\mathbf}
\newcommand{\mbb}{\mathbb}
\newcommand{\mrm}{\mathrm}

\newcommand{\LR}[2]{\left\llbracket \begin{matrix} #1\\#2 \end{matrix} \right\rrbracket}
\newcommand{\ext}{{ \mathfrak{Ext}}}
\def\scrP{\mathscr{P}}
\newcommand{\bk}{{\mathbb K}}
\newcommand{\cc}{{\mathcal C}}
\newcommand{\gc}{{\mathcal GC}}
\newcommand{\dg}{{\rm dg}}
\newcommand{\ce}{{\mathcal E}}
\newcommand{\cs}{{\mathcal S}}
\newcommand{\cP}{{\mathcal P}}
\newcommand{\cl}{{\mathcal L}}
\newcommand{\cf}{{\mathcal F}}
\newcommand{\cx}{{\mathcal X}}
\newcommand{\cy}{{\mathcal Y}}
\newcommand{\ct}{{\mathcal T}}
\newcommand{\cu}{{\mathcal U}}
\newcommand{\cv}{{\mathcal V}}
\newcommand{\cn}{{\mathcal N}}
\newcommand{\mcr}{{\mathcal R}}
\newcommand{\ch}{{\mathcal H}}
\newcommand{\ca}{{\mathcal A}}
\newcommand{\cb}{{\mathcal B}}
\newcommand{\ci}{{\I}_{\btau}}
\newcommand{\cj}{{\mathcal J}}
\newcommand{\cm}{{\mathcal M}}
\newcommand{\cp}{{\mathcal P}}
\newcommand{\cg}{{\mathcal G}}
\newcommand{\cw}{{\mathcal W}}
\newcommand{\co}{{\mathcal O}}
\newcommand{\cq}{{Q^{\rm dbl}}}
\newcommand{\cd}{{\mathcal D}}
\newcommand{\ck}{\widetilde{\mathcal K}}
\newcommand{\calr}{{\mathcal R}}
\newcommand{\dmno}{\ytableausetup{boxsize=3pt}\ydiagram{1,1}}
\newcommand{\dmnoB}{\ytableausetup{boxsize=3pt}\ydiagram{2}}
\newcommand{\mno}{\ytableausetup{boxsize=3pt}\ydiagram{1}}
\newcommand{\iLa}{\la^{\imath}}
\newcommand{\La}{\Lambda}
\newcommand{\ol}{\overline}
\newcommand{\ul}{\underline}
\newcommand{\ow}{\widetilde}

\newtheorem{theorem}{Theorem}[section]
\newtheorem{acknowledgement}[theorem]{Acknowledgement}
\newtheorem{algorithm}[theorem]{Algorithm}
\newtheorem{assumption}[theorem]{Assumption}
\newtheorem{axiom}[theorem]{Axiom}
\newtheorem{case}[theorem]{Case}
\newtheorem{claim}[theorem]{Claim}
\newtheorem{conclusion}[theorem]{Conclusion}
\newtheorem{condition}[theorem]{Condition}
\newtheorem{conjecture}[theorem]{Conjecture}
\newtheorem{construction}[theorem]{Construction}
\newtheorem{corollary}[theorem]{Corollary}
\newtheorem{criterion}[theorem]{Criterion}
\newtheorem{definition}[theorem]{Definition}
\newtheorem{example}[theorem]{Example}
\newtheorem{exercise}[theorem]{Exercise}
\newtheorem{lemma}[theorem]{Lemma}
\newtheorem{notation}[theorem]{Notation}
\newtheorem{problem}[theorem]{Problem}
\newtheorem{proposition}[theorem]{Proposition}
\newtheorem{solution}[theorem]{Solution}
\newtheorem{summary}[theorem]{Summary}
\numberwithin{equation}{section}

\theoremstyle{remark}
\newtheorem{remark}[theorem]{Remark}

\def \cz{\mathcal Z}

\def \bmu{\boldsymbol{\mu}}
\def \bnu{\boldsymbol{\nu}}
\def \bla{\boldsymbol{\la}}

\def \bfK{{\mathbf K}}

\def \bA{{\mathbf A}}
\def \ba{{\mathbf a}}
\def \bL{{\mathbf L}}
\def \bF{{\mathbf F}}
\def \bS{{\mathbf S}}
\def \bC{{\mathbf C}}
\def \bU{{\mathbf U}}
\def \bc{{\mathbf c}}
\def \fpi{\mathfrak{P}^\imath}
\def \Ni{N^\imath}
\def \fp{\mathfrak{P}}
\def \fg{\mathfrak{g}}
\def \fk{\fg^\theta}  

\def \fn{\mathfrak{n}}
\def \fh{\mathfrak{h}}
\def \fu{\mathfrak{u}}
\def \fv{\mathfrak{v}}
\def \fa{\mathfrak{a}}
\def \fq{\mathfrak{q}}
\def \Z{{\mathbb Z}}
\def \F{{\mathbb F}}
\def \D{{\mathbb D}}
\def \bB{{\mathbb B}}
\def \C{{\mathbb C}}
\def \N{{\mathbb N}}
\def \Q{{\mathbb Q}}
\def \G{{\Bbb G}}
\def \P{{\Bbb P}}
\def \K{{\mathbf k}}
\def \bK{{\mathbb K}}

\def \E{{\Bbb E}}
\def \A{{\Bbb A}}
\def \L{{\Bbb L}}
\def \R{{\Bbb R}}
\def \I{{\Bbb I}}
\def \BH{{\Bbb H}}
\def \T{{\Bbb T}}
\def \de{{\delta}}
\def \vth{{\theta}}

\def \cN{{\mathcal N}}

\newcommand{\arxiv}[1]{\href{http://arxiv.org/abs/#1}{\tt arXiv:\nolinkurl{#1}}}

\newcommand{\nc}{\newcommand}
\newcommand{\browntext}[1]{\textcolor{brown}{#1}}
\newcommand{\greentext}[1]{\textcolor{green}{#1}}
\newcommand{\redtext}[1]{\textcolor{red}{#1}}
\newcommand{\bluetext}[1]{\textcolor{blue}{#1}}
\newcommand{\brown}[1]{\browntext{ #1}}
\newcommand{\green}[1]{\greentext{ #1}}
\newcommand{\red}[1]{\redtext{ #1}}
\newcommand{\blue}[1]{\bluetext{ #1}}

\newcommand{\wtodo}{\todo[inline,color=orange!20, caption={}]}
\newcommand{\lutodo}{\todo[inline,color=green!20, caption={}]}

\title[$\imath$Hall algebra of Jordan quiver]{$\imath$Hall algebra of Jordan quiver and $\imath$Hall-Littlewood functions}

\author[Ming Lu]{Ming Lu}
\address{Department of Mathematics, Sichuan University, Chengdu 610064, P.R. China}
\email{luming@scu.edu.cn}

\author[Shiquan Ruan]{Shiquan Ruan}
\address{ School of Mathematical Sciences, Xiamen University, Xiamen 361005, P.R. China}
\email{sqruan@xmu.edu.cn}

\author[Weiqiang Wang]{Weiqiang Wang}
\address{Department of Mathematics, University of Virginia, Charlottesville, VA 22904, USA}
\email{ww9c@virginia.edu}


\begin{abstract}
We show that the $\imath$Hall algebra of the Jordan quiver is a polynomial ring in infinitely many generators and obtain transition relations among several generating sets. We establish a ring isomorphism from this $\imath$Hall algebra to the ring of symmetric functions in two parameters $t, \theta$, which maps the $\imath$Hall basis to a class of (modified) inhomogeneous Hall-Littlewood ($\imath$HL) functions. The (modified) $\imath$HL functions admit a formulation via raising and lowering operators. We formulate and prove Pieri rules for (modified) $\imath$HL functions. The modified $\imath$HL functions specialize at $\theta=0$ to the modified HL functions; they specialize at $\theta=1$ to the deformed universal characters of type C, which further specialize at $(t=0, \theta =1)$  to the universal characters of type C.
\end{abstract}

\maketitle
 \setcounter{tocdepth}{1}
 \tableofcontents

\section{Introduction}

\subsection{Hall to $\imath$Hall}

Hall algebra is a fundamental concept in representation theory with intimate connections and applications to combinatorics of symmetric functions,  quantum groups and algebraic geometry.
Hall algebra of the Jordan quiver was historically the first example (studied by E.~Steinitz and then by P.~Hall), and its Hall basis leads to Hall-Littlewood (HL) polynomials, a distinguished basis for the ring of symmetric functions; cf. \cite{Mac95}. Ringel's Hall algebras \cite{Rin90} (see also Green \cite{Gr95}) for general quivers have provided a realization of halves of quantum groups and led to Lusztig's construction of canonical basis \cite{Lus93}.

Jing \cite{J91} realized Hall-Littewood polynomials via vertex operators which are closely related to the boson-fermion correspondence. 
It was shown by Tsilevich \cite{Ts06} that Hall-Littlewood polynomials are  the off-shell wave functions of the $t$-boson model. The $t$-boson model arose in \cite{BB92} as a discretization of the Bose gas model. The connection between Hall-Littlewood polynomials and the $t$-boson model has been further explored by Wheeler and Zinn-Justin \cite{WZ16}. 

Motivated by Bridgeland's realization \cite{Br13} of quantum groups via Hall algebras of complexes associated to general quivers, two of the authors \cite{LW19, LW20, LW21} have recently realized   the $\imath$quantum groups arising from quantum symmetric pairs via $\imath$Hall algebras of quivers. (See \cite{LuoW22} for an exposition on quantum symmetric pairs and applications.) The $\imath$Hall algebras are set up in the framework of semi-derived Hall algebras \cite{LP21, Lu19}, which were in turn built on the earlier foundational works of Bridgeland and M.~Gorsky \cite{Gor2}. The $\imath$Hall algebra of the projective line has recently been developed \cite{LRW20} to realize a Drinfeld type presentation of $q$-Onsager algebra (the simplest affine $\imath$quantum group); compare \cite{Ka97, BKa01, Sch04, Sch06}.

\subsection{$\imath$Hall algebra of the Jordan quiver}

We aim to develop in this paper the $\imath$Hall algebra of the Jordan quiver and connections to symmetric functions, adding a new chapter to the classic constructions of Steinitz, Hall, and Macdonald.

By definition, this $\imath$Hall algebra is based on the Jordan quiver enhanced by an additional arrow $\varepsilon$ with $\varepsilon^2=0$ and $\alpha\varepsilon=\varepsilon \alpha$:
\vspace{-2mm}
\begin{figure}[h]
\centering
\begin{tikzpicture}[scale=1.1]
\draw[color=purple] (7,0) .. controls (7.2,0.35) and (7.6,0.35) .. (7.6,-0.05);
\draw[-latex,color=purple] (7.6,-0.05) .. controls (7.6,-0.45) and (7.2,-0.45) .. (7,-0.1);
\draw[-] (6.7,0) .. controls (6.5,0.35) and (6.1,0.35) .. (6.1,-0.05);
\draw[-latex] (6.1,-0.05) .. controls (6.1,-0.45) and (6.5,-0.45) .. (6.7,-0.1);
\node at (6.85,0 ){ \small$1$};
\node at (5.9,0 ){ \small$\alpha$};
\node at (7.8,0 ){ \small\color{purple}{$\varepsilon$}};
\end{tikzpicture}
\end{figure}
\vspace{-3mm}

\noindent Let $\bfk =\mathbb F_q$ be a finite field of $q$ elements. We show that the $\imath$Hall algebra $\iRH(\bfk \QJ)$ of the Jordan quiver $\QJ$ is isomorphic to a polynomial algebra in infinitely many generators: $[K_S]$ and $[S^{(r)}]$, for $r\ge 1$, where $S^{(r)}$ denotes the indecomposable $\bfk \QJ$-module of length $r$. The new generator $[K_S]$ arises from an acyclic complex $K_S$ associated to the simple module $S=S^{(1)}$; central generators like $[K_S]$ arise naturally in other $\imath$Hall algebras and can be treated as parameters \cite{LW19, LRW20}. The semi-derived Hall algebras and $\imath$Hall algebras used in this paper are slight modifications from earlier versions so that $[K_S]$ is not required to be invertible in $\iRH(\bfk \QJ)$.

The algebra $\iRH(\bfk \QJ)$ is naturally a $\N$-graded algebra with degrees of generators given by $\big |[S^{(r)}] \big | =r$ and $\big | [K_S] \big | =2$.
The $\imath$Hall algebra admits a natural $\imath$Hall basis $\{[S^\la]\}$ over $\Q(t)\big[[K_S] \big]$, where $\la$ runs over all partitions. We obtain several explicit formulas for the multiplication with $[S^{(r)}]$ (or with $[S^{1^r)}]$) in the $\imath$Hall basis (that is, horizontal and vertical Pieri rules). Remarkably, different Pieri rules are related by Green's formula \cite{Gr95}.

We can formulate a generic $\imath$Hall algebra $\iRH(\QJ)$, when $q^{-1}$ in $\iRH(\bfk \QJ)$ is replaced by a generic (HL) parameter $t$ throughout. We regard $\vth :=t [K_S]$ as a second deformation parameter. Setting $\vth=0$ in the $\imath$Hall algebra recovers (Steinitz's) Hall algebra of the Jordan quiver. Various formulas in the Hall algebra as presented in \cite[III]{Mac95} are the leading terms (i.e., at $\vth=0$) of formulas in the $\imath$Hall algebra established in this paper.

Our initial computations with the $\imath$Hall basis suggested a generating function formulation generalizing that for Hall-Littlewood polynomials; see \eqref{HL functions}; compare \cite[III, (2.15)]{Mac95}. This naturally leads to a formulation via raising and lowering operators $R_{ij}$ and $L_{ij}$, which turns out to be one of the key points in our approach. We have thus adopted this as our definition of polynomials $\iV_\la$ in commuting variables $v_r$, for $r\ge 1$ (see Definition~ \ref{def:V}), which shall give rise to (modified) $\imath$HL functions:
\begin{align*}
\iV_\la = \prod_{1\le i<j} \frac{1 - \vth L_{ij}}{1 - \vth t L_{ij}} \frac{1 -R_{ij}}{1 -tR_{ij}} \vv_\la.
\end{align*}
This formula at $\vth=0$ recovers \cite[III, (2.15$'$)]{Mac95} when interpreting $v_r$ as the symmetric functions $q_r$. Following \cite[\S2.1,(6)]{Tam11}, we shall refer to $\iV_\la$ as Giambelli type polynomials.

In this way, we have set up a {\em linear} isomorphism from the generic $\imath$Hall algebra to the ring of symmetric functions over $\Q(t)[\vth]$,
\[
\Phi:  \iRH(\QJ) \longrightarrow \La_{t,\vth},
\]
which identify the $\imath$Hall basis with the Giambelli type polynomials associated to partitions.
We shall show that $\Phi$ is an algebra isomorphism.

\subsection{Combinatorics}

To that end, we shall establish Pieri rules for the Giambelli type polynomials $\iV_\la$ (or $\imath$HL functions), by extending an approach of Tamvakis.

The raising operator formalism played a fundamental role in the Schubert calculus of classical type developed by Buch, Kresch and Tamvakis \cite{BKT08}. Tamvakis developed the raising operator approach further \cite{Tam11}; in particular, starting from a raising operator definition of Giambelli polynomials (and then HL functions), he  found a direct proof of the horizontal Pieri rule \cite[III, (5.7$'$)]{Mac95}. Such a Pieri rule was earlier obtained by Morris \cite{Mor64}.

In spite of the new presence of lowering operators, we are able to generalize Tamvakis' approach to establish a first horizontal Pieri rule for the $\iV_\la$. Two mirror identities among HL functions, which Tamvakis established by the raising operator calculus, are instrumental in this approach. It is remarkable that these mirror identities remain unchanged in our new setup involving lowering operators. While the first identity is natural in terms of raising operators, the second one is most natural in terms of lowering operators.

Actually, we have some freedom in passing from Giambelli type polynomials to symmetric functions. If we identify $\vv_r =q_r$ (the HL function associated to one-row partitions), then we obtain a version of $\imath$HL functions $\iV_\la =\iQ_\la$ which specialize at $\vth=0$ to HL functions $Q_\la$. As we learned from S.~Okada, a specialization of the iHL functions $Q^\imath_\lambda$ (for strict partitions $\la$) coincides with the universal symplectic Q-functions $ \mathcal Q^C_\lambda$ introduced in his paper \cite[\S3.1]{Ok21}, i.e.,
$Q^\imath_\lambda|_{t=-1,\theta=1} = \mathcal Q^C_\lambda;$ see Remark~\ref{rem:Okada}.
On the other hand, by identifying $\vv_r =h_r$, then we obtain the modified $\imath$HL functions $\iV_\la =\iH_\la$ which specialize at $\vth=0$ to the (Garsia's) modified HL functions $H_\la$.

In the process of developing the $\imath$HL functions, we have come to realize that  our $\imath$HL functions admit another incarnation!
In a remarkable work \cite{SZ06}, Shimozono and Zabrocki used creation operators (i.e., modified versions of Jing's HL vertex operators \cite{J91, SZ01}) to introduce $t$-deformed universal characters of classical types, denoted by $\mathbb B_\la^\diamond$ (where $\diamond = \emptyset, \dmnoB,$ $\dmno,$ $\mno$). At $t=0$,  $\mathbb B_\la^\diamond$ specializes to the universal characters $s_\la^\diamond$ introduced in \cite{KT87} as suitable stable limits of classical irreducible characters. 

In the notation of \cite{SZ06}, the $s_\la^{\dmno}$ are universal characters of type C.
It turns out that we have the identification ${\iH_\la} \big |_{\vth=1} = {\mathbb B}_\la^{\dmno}$, as one can convert the definitions via creation operators and via raising/lowering operators. (We thank Mark Shimozono for his help in confirming this expectation.)
As ${\bB}_\la^{\dmno}$ is inhomogeneous, inserting powers of $\vth$ (with degree $|\vth |= 2$) to ${\bB}_\la^{\dmno}$ in a unique homogeneous way we reconstruct $\iH_\la$. 
 The Pieri rules were not addressed in \cite{SZ06}.

The above unexpected type C interpretation fits well with the philosophy of the ``$\imath$fication"-program and has a precedent. In \cite{BW18} (and also \cite{BKLW18}), Kazhdan-Lusztig basis of type C was formulated in terms of $\imath$canonical basis arising from $\imath$quantum groups of type AIII (with a type A diagram together with the diagram involution as the input), and $\imath$quantum groups admit an $\imath$Hall algebra realization \cite{LW19}. Our $\imath$Hall algebra $\iRH(\QJ)$ arises from Jordan $\imath$quiver (i.e., the Jordan quiver coupled with a trivial involution). An $\imath$fication process produces constructions of ``type A with involution", which may happen to have a type BCD reincarnation.

The raising operator formulation also served as a starting point and played a crucial role in the recent breakthrough \cite{BMPS19} in the theory of Macdonald polynomials, Catalan functions and $k$-Schur functions. 
We expect a natural {\em classical} synthesis of the BMPS framework with this work (and \cite{SZ06}); this will be explored elsewhere.

Below let us formulate our main results more precisely, starting from the combinatorics side.

\subsection{$\imath$Pieri rules}

For partitions $\mu, \nu$, we denote $\nu \leq \mu$ if $\nu_i \le \mu_i$ for all $i$, and denote $\nu \stackrel{a}{\rightarrow} \mu$ if $\nu \leq \mu$ and $\mu-\nu$ is a horizontal $a$-strip. The definitions for the polynomials $\varphi_{\mu/\nu}(t), \psi_{\la/\nu}(t)$ can be found in \eqref{eq:phipsi}, following \cite[III, (5.8)-(5.8$'$)]{Mac95}.

\begin{customthm}{{\bf A(I)}}   [Down-up horizontal $\imath$Pieri rule, Theorem~\ref{thm:Pieri-Ho}]
For $r\ge 1$ and any partition $\mu$, we have
\begin{align*}
\iV_\mu  \cdot \vv_r
&=\sum_{a+b=r} \sum_{\nu \stackrel{a}{\rightarrow} \mu ,\nu\stackrel{b}{\rightarrow} \la} \vth^{a}\, \varphi_{\mu/\nu}(t) \psi_{\la/\nu}(t) \;  \iV_\la.
\end{align*}
\end{customthm}
We refer to this horizontal $\imath$Pieri rule as a ``down-up" version as we first go down to intermediate Young diagram $\nu$ from the initial diagram $\mu$ and then up to an end diagram $\la$. The leading term (= the $\vth^0$-term) in the above formula recovers \cite[III, (5.7$'$)]{Mac95}, while the lowest term  (= the $\vth^r$-term) encodes a dual Pieri rule for HL functions. The down-up horizontal $\imath$Pieri rule in Theorem~{\bf A(I)}  is proved by generalizing the approach of Tamvakis \cite{Tam11}.

Recall the specialization of $\iV_\mu$ at $t=0$ is essentially the universal character $s_\mu^{\dmno}$. The down-up horizontal and up-down vertical Pieri rules for (universal) characters of type C were well known \cite{Su90} (also cf. \cite{Ok16}).  Theorem~{\bf A(I)} specializes at $t=0$ (where all the polynomial coefficients in t become 1) to such a down-up horizontal Pieri rule.

We shall formulate several additional $\imath$Pieri rules. First, let us present an ``up-down" version of the horizontal $\imath$Pieri rule. This up-down horizontal version at $t=0$ can be easily derived from the up-down {\em vertical} Pieri rule by applying symmetries of universal characters (cf. \cite{KT87, SZ06}).

\begin{customthm}{{\bf A(II)}} [Up-down horizontal $\imath$Pieri rule]
For $r\ge 1$ and any partition $\mu$, we have
\begin{align*}
\iV_\mu*  v_r =&\sum_{a+b =r}
\sum_{\la \stackrel{a}{\rightarrow} \xi} \sum_{\mu\stackrel{b}{\rightarrow} \xi}
\vth^a  \varphi_{\xi/\la}(t) \psi_{\xi/\mu}(t) \;  \iV_\la
\\\notag
&+\sum_{a+b=r}\sum\limits_{i=1}^{\min(a,b)}\sum_{\la \stackrel{a-i}{\rightarrow} \xi}\sum_{\mu \stackrel{b-i}{\rightarrow} \xi}
\vth^a  (t^i -t^{i-1})\varphi_{\xi/\la}(t) \psi_{\xi/\mu}(t) \; \iV_\la.
\end{align*}
\end{customthm}

We shall now formulate 2 versions, down-up and up-down, of vertical $\imath$Pieri rules.
For partitions $\mu, \nu$, we use $\mu \stackrel{a}{\downarrow} \nu$ to denote that $\nu \leq \mu$ and $\mu-\nu$ is a vertical $a$-strip. The definitions for $\varphi_r(t)$, $b_{\nu}(t)$, and $f_{\nu,(1^{m})}^\mu (t)$ can be found in \eqref{eq:varphi:r}--\eqref{def:bla} and \eqref{eq:fmu}, following \cite[III, (2.12), (3.2)]{Mac95}.

\begin{customthm}{{\bf B}}
Let $r\geq1$ and $\mu$ be a partition.  We have
\begin{enumerate}
\item
(Down-up vertical $\imath$Pieri rule)
\begin{align*}
\iV_\mu\cdot \iV_{(1^r)}
&= \sum_{a+b=r} \sum_{\la\stackrel{b}{\downarrow} \nu}  \sum_{\mu \stackrel{a}{\downarrow} \nu}
\vth^a \frac{b_{\nu}(t)}{b_{\la}(t)} \varphi_r(t) f_{\nu,(1^{b})}^\la(t) f_{\nu,(1^a)}^\mu(t) \; \iV_\la;
\end{align*}

\item
(Up-down vertical $\imath$Pieri rule)
\begin{align*}
\iV_\mu\cdot \iV_{(1^r)}
&= \sum_{a+b=r}\sum_{i=0}^{\min(a,b)}\sum_{\xi \stackrel{a-i}{\downarrow}\la}  \sum_{\xi\stackrel{b-i}{\downarrow}\mu}
(-1)^i \vth^a t^{\frac{i(i-1)}{2}} \frac{\varphi_r(t)}{\varphi_i(t)}  \frac{b_{\mu}(t)}{b_{\xi}(t)}    f_{\la,(1^{{a-i}})}^\xi(t) f_{\mu,(1^{b-i})}^\xi (t)  \; \iV_\la.
\end{align*}
\end{enumerate}
\end{customthm}

However, our proofs of the 3 Pieri rules in Theorem~{\bf A(II)} and Theorem~{\bf B} are indirect, and they rely essentially on the computations in $\imath$Hall algebra of the Jordan quiver. It will be very interesting to find direct algebraic/combinatorial proofs for these 3 Pieri rules.
At the $t=0$ limit, our up-down vertical Pieri rule reduces to the counterpart for (universal) characters of type C \cite{Su90} (also cf. \cite{Ok16}).

\begin{proof} [Proof of Theorems~{\bf A(II)} and ~{\bf B}]
Thanks to the algebra isomorphism in Theorem~{\bf C} below (or in Theorem~\ref{thm:iso}), it suffices to establish the counterparts of these 3 $\imath$Pieri rules in the $\imath$Hall algebra $\iRH(\bfk \QJ)$. These are given in Proposition~\ref{prop:PieriHall-Ho2}, Proposition \ref{prop:haPieri-v}, and Proposition~ \ref{prop:haPieri-vv}, respectively.
\end{proof}

\subsection{The ring isomorphism}

For any partition $\la =(\la_1, \la_2, \ldots)$, denote
\begin{align}
\label{eq:nla}
n(\la) =\sum_{i\geq1}(i-1)\la_i,\qquad\quad
|\la|=\sum_{i\geq1}\la_i.
\end{align}
We replace elements $[S^{(\la)}]$, $[K_S]$ and $q^{-1}$ in the setting of $\iRH(\bfk\QJ)$ by $\fu_{\la}$, $K_\delta$ and by $t$ in the setting of the generic $\imath$Hall algebra $\iRH(\QJ)$.
The following is a generic version of Theorem~\ref{thm:iso} (also see Remark~\ref{rem:PQ} if one compares with the algebra isomorphism in \cite[(3.4)]{Mac95}).

\begin{customthm}{{\bf C}}
There exists a $\Q(t)$-algebra isomorphism $\Phi:  \iRH(\QJ) \longrightarrow \La_{t,\vth}$ such that
\begin{align*}
\Phi (K_\de) =t^{-1} \vth,\qquad
\Phi (\fu_{(r)}) =t^{-r} v_r \quad (r\ge 1).
\end{align*}
Moreover,
$\Phi (\fu_\la)= t^{-|\lambda| -n(\lambda)}\iV_\la$,
for any partition $\la$.
\end{customthm}

The ring structure of $\La_{t,\vth}$ in the basis $\iV_\la$ is determined by the Pieri rule in Theorem~{\bf A(I)} (or any other Pieri rule). Theorem~{\bf C} (or Theorem~\ref{thm:iso}) follows readily if we  establish the counterpart of this Pieri rule in the $\imath$Hall algebra $\iRH(\bfk \QJ)$; this is accomplished in Proposition~ \ref{prop:PieriHall-Ho}.

\subsection{Transitions among several generating sets}

To formulate some generating functions below, we need to enlarge the field $\Q$ to $\Q(\sqq)$, where $\sqq =\sqrt{q}$. Let $z$ be an indeterminate. Recall $\varphi(t)$ and $\exp_q \big(\frac{x}{1-q} \big)$ from \eqref{eq:varphi:r} and \eqref{eq:Euler}. We introduce the following generating functions of several natural generating sets in the $\imath$Hall algebra $\iRH(\bfk \QJ)$ (compare \cite{BKa01, Sch06, LRW20}):
\begin{align}
\label{def:haE}
\haE(z) = &\sum_{r\geq0} \sqq^{r(r-1)} \frac{[S^{(1^r)}]}{|\Aut(S^{(1^r)})|}z^r, 
\\
\label{def:haH}
\haH(z) =&\sum_{r\geq0} \sum\limits_{\la\vdash r} \frac{[S^{(\la)}]}{|\Aut(S^{(\la)})|}z^r,
\\
\label{def:haT}
\haT(z) =&\sum_{r\geq0}\, [S^{(r)}]z^r,
\\
\label{def:haP}
\haP(z) =&\sum_{r\geq1}\haP_rz^{r-1},
\quad \text{where }
\haP_r =\sum\limits_{\la\vdash r} \varphi_{l(\la)-1}(q)\frac{[S^{(\la)}]}{|\Aut(S^{(\la)})|}.
\end{align}

\begin{customthm}{{\bf D}}
[Propositions~\ref{prop:HE}, \ref{prop:QE}, \ref{prop:QP}, \ref{prop:EHP}]
The following identities hold in the $\imath$Hall algebra $\iRH(\bfk \QJ)$:
\begin{align}
\haH(z)\haE(-z)&= \exp_q \Big(\frac{ [K_S] z^2}{1-q} \Big),
 \label{HE}
\\
\haT(z) &= \big(1-q[K_S] z^2 \big) \frac{\haE(-z)}{\haE(- qz)} = \big(1-[K_S] z^2 \big)^{-1} \frac{\haH(qz)}{\haH(z)},
\label{TE}
\\
\haT(z)
&= \exp \Big(\sum_{r\geq1} \frac{q^r-1}{r} \haP_r z^{r} \Big)  \exp \Big(\sum_{k \geq1} \frac{1-q^{k}}{2k} [K_S]^k z^{2k} \Big),
\label{TP}
\\
\haE(z)
&= \exp \Big(\sum_{r\geq1} \frac{(-1)^{r-1}}{r} \haP_r z^{r}  \Big) \exp_q \Big(\frac{ [K_S] z^2}{1-q} \Big)^{\frac12},
\label{EP}
\\
\haH(z)
&= \exp \Big( \sum_{r\geq1} \frac{1}{r} \haP_r z^{r} \Big) \exp_q \Big(\frac{ [K_S] z^2}{1-q} \Big)^{\frac12}.
\label{HP}
\end{align}
\end{customthm}
The identities \eqref{HE}--\eqref{TE} are proved by applying Pieri rules and reducing to some combinatorial identities. The proof of the identity \eqref{TP} turns out to be most challenging, and this identity has applications to $\imath$Hall algebra of the projective line \cite{LRW20}. The identities \eqref{EP}--\eqref{HP} follows from \eqref{HE}--\eqref{TP} by some formal manipulations.
We can transport the identities in Theorem~{\bf D} to corresponding identities in the ring of symmetric functions via the algebra isomorphism $\Phi$ in Theorem~{\bf C}. Much remains to be further explored combinatorially.

While the main motivation of this paper is the interaction between Hall algebra and symmetric functions, it will be interesting to ask if the $\imath$Hall-Littlewood polynomials have any connections to mathematical physics; compare \cite{Ts06, WZ16}. 

\subsection{Organization}

In Section~\ref{sec:Pieri}, we establish the down-up horizontal Pieri rule for the Giambelli type polynomials $\iV_\la$, generalizing the approach of Tamvakis. We then make suitable identifications of $v_r$ to formulate the (modified) $\imath$HL functions. We also formulate the precise relations with (deformed) universal characters of type C.

We review and formulate a version of semi-derived Hall algebras and $\imath$Hall algebras in Section~\ref{sec:Hall}.
In Section~\ref{sec:iHall}, we define the $\imath$Hall algebra $\iRH(\bfk \QJ)$ of the Jordan quiver. Then we establish the 4 versions (down-up vs up-down, horizontal vs vertical) Pieri rules in $\iRH(\bfk \QJ)$.

We collect and establish in Section~\ref{sec:comb} several technical identities which we shall need in the next section. In Section~\ref{sec:identities}, we study several distinguished generating sets for the $\imath$Hall algebra $\iRH(\bfk \QJ)$ in the form of generating functions; they correspond to the $(t,\vth)$-deformations of the generating sets of elementary, homogeneous and power sum symmetric functions. We establish explicit formulas for transitions from one to another generating set.

\vspace{2mm}

\noindent{\bf Acknowledgment.}
We thank Jennifer Morse, Harry Tamvakis, Mark Shimozono, and Soichi Okada for stimulating discussions, very helpful correspondences and references. We thank Harry Tamvakis for explaining and clarifying his original approach toward Pieri rule. ML is partially supported by the National Natural Science Foundation of China (No. 12171333).
SR is partially supported by the National Natural Science Foundation of China (No. 11801473).
WW is partially supported by the National Science Foundation grant DMS-2001351.

\section{$\imath$Hall-Littlewood functions and $\imath$Pieri rules}
  \label{sec:Pieri}

In this section, we define Giambelli type polynomials and establish their Pieri rules. We transform the Giambelli type polynomials into $\imath$HL functions, and compare the latter with (modified) HL functions and deformed universal characters of type C.

\subsection{A Giambelli type formula}
Let $\vv_r$ be commuting variables of degree $r$, for $r \in \Z_{\ge 1}$. Let us set $\vv_0=1$ and $\vv_a=0$, for $a<0$. We define
\[
\vv_\alpha =\vv_{\alpha_1}\vv_{\alpha_2} \cdots
\]
multiplicatively for any integer sequence $\alpha =(\alpha_1, \alpha_2, \ldots)$ with finitely many nonzero entries.
Given 2 integer sequences $\alpha, \beta$, we denote $\alpha \ge \beta$ (or $\beta \le \alpha$) if $\alpha_i \ge \beta_i$, for all $i\ge 1$; in particular, $\alpha \ge 0$ denotes that $\alpha$ is a composition. For a composition $\alpha$, we denote
\begin{align}
 \label{eq:number}
\begin{split}
\ell(\alpha) &= \text{number of nonzero parts in } \alpha,
\\
|\alpha| &=\alpha_1 +\alpha_2+\cdots.
\end{split}
\end{align}

Let $1\le i<j$. The raising operator $R_{ij}$ acts on integer sequences by
\[
R_{ij} \alpha =(\ldots, \alpha_i+1, \ldots, \alpha_j-1, \ldots).
\]
Similarly, we define the lowering operator $L_{ij}$ acting on integer sequences by
\[
L_{ij} \alpha =(\ldots, \alpha_i-1, \ldots, \alpha_j-1, \ldots).
\]
While the raising operators were originally introduced by Young (cf. \cite{Mac95}), some of the significant applications have appeared only recently; see \cite{BKT08, Tam11}. The lowering operators above also appeared in \cite{Lec06}.
The actions of raising and lowering operators on $\vv_\alpha$ are given by letting
\[
R_{ij} \vv_\alpha =\vv_{R_{ij}\alpha}, \qquad L_{ij} \vv_\alpha =\vv_{L_{ij}\alpha}.
\]
Note that all raising and lowering operators commute with each other.

Let $t, \vth$ be indeterminates. Denote by $\cp_n$ the set of partitions $\la$ of $n$ with dominance order $\unrhd$, and we sometime write $\la \vdash n$.
Set $\cp =\cup_{n\ge 0}\cp_n$. The dominance order $\unrhd$ also makes sense on the set of all compositions of $n$.
The ring
\[
\La_{t,\vth} = \Q(t)[\vth] [v_1, v_2, \ldots]
\]
admits a basis $\{v_\la \mid \la \in \cp \}$.

\begin{definition}
 \label{def:V}
For any integer vector $\alpha$, we define $V_\alpha, \iV_\alpha \in \La_{t,\vth}$ by:
\begin{align}  \label{eq:LRq}
V_\alpha &= \prod_{1\le i<j} \frac{1 -R_{ij}}{1 -tR_{ij}} \vv_\alpha,
\qquad
\iV_\alpha = \prod_{1\le i<j} \frac{1 - \vth L_{ij}}{1 - \vth t L_{ij}} \frac{1 -R_{ij}}{1 -tR_{ij}} \vv_\alpha.
\end{align}
Note that $V_r =\iV_{r} =\vv_r$.
\end{definition}
Justification of applying raising (and lowering) operators properly can be found in \cite[\S2.2]{Tam12} and \cite[(4.1) and below]{BMPS19}. Let
\begin{align}
 \label{eq:F}
F(u)=\frac{1-u}{1-tu}=1+ (1-t^{-1}) \sum_{r\geq1} t^r u^r.
\end{align}Denote
\begin{align}  \label{eq:L}
\L_\ell &:= \prod_{1\le i<\ell}  \frac{1 -\vth L_{il}}{1 - \vth t L_{il}}
=\prod_{1\le i< \ell} \Big(1 + (1-t^{-1}) \sum_{s\ge 1} \vth^{s} t^s L_{i\ell}^s \Big),
\\
\qquad
\R_\ell &:= \prod_{1\le i<\ell}  \frac{1 -R_{il}}{1 -tR_{il}}
= \Big(1 + (1-t^{-1}) \sum_{p\ge 1} t^{p} R_{i\ell}^p \Big).
\end{align}

Denote
\begin{align}   \label{eq:D}
D_\ell : =\prod_{1\le i<j \le \ell}  \frac{1 -R_{ij}}{1 -tR_{ij}}\frac{1 -\vth L_{ij}}{1 - \vth t L_{ij}}
=\prod_{k=1}^\ell \L_k \cdot \prod_{k=1}^\ell \R_k.
\end{align}

Below we follow \cite[\S1-2]{Tam11} to manipulate the operators $D_\ell$. We have
\begin{align*}
D_\ell
&= D_{\ell-1} \prod_{1\le i<\ell} \frac{1 -R_{i \ell}}{1 -tR_{i \ell}} \frac{1 -\vth L_{i \ell}}{1 - \vth t L_{i \ell}}
\\
&= D_{\ell-1} \prod_{1\le i< \ell}
\Big(1 + (1-t^{-1}) \sum_{p\ge 1} t^{p} R_{i\ell}^p \Big)
\cdot
\prod_{1\le i< \ell} \Big(1 + (1-t^{-1}) \sum_{s\ge 1} \vth^{s} t^s L_{i\ell}^s \Big).
\end{align*}
Applying $D_\ell$ to $\vv_\alpha \vv_r$, where $\alpha =(\alpha_1, \ldots, \alpha_{\ell-1})$ is an integer vector, and $r\in \Z$, we obtain
\begin{align}
\iV_{\alpha,r} &=D_{\ell} \vv_\alpha \vv_r
= D_{\ell-1} \prod_{1\le i< \ell} \Big(1 + (1-t^{-1}) \sum_{p\ge 1} t^{p} R_{i\ell}^p \Big) \cdot
\sum_{\beta} \vth^{|\beta|} t^{|\beta|} (1-t^{-1})^{\ell(\beta)} \vv_{\alpha -\beta} \vv_{r-|\beta|}
 \notag \\
&= \sum_{\beta, \gamma} \vth^{|\beta|} t^{|\beta|+|\gamma|} (1-t^{-1})^{\ell(\beta) +\ell(\gamma)}\; \iV_{\alpha -\beta +\gamma} \vv_{r-|\beta| -|\gamma|},
\label{eq:recur}
\end{align}
summed over all compositions $\beta, \gamma \in \N^{\ell-1}$. The formula \eqref{eq:recur} gives a recursive formula for computing $\iV_{\alpha}$. Repeatedly applying \eqref{eq:recur}, we see that
\[
\iV_{\alpha} \in v_{\alpha} + \sum_{\alpha' \rhd \alpha}\Q(t)v_{\alpha'}+ \sum_{|\alpha''| <|\alpha|}\Q(t)[\vth]v_{\alpha''}.
\]
As $\{v_\la \mid \la \in \cp \}$ forms a $\Q(t)[\vth]$-basis for $\La_{t,\vth}$, the following is immediate.
\begin{corollary}
$\{\iV_\la \mid \la \in \cp \}$ forms a $\Q(t)[\vth]$-basis for $\La_{t,\vth}$.
\end{corollary}

\subsection{Basic properties of $\iV_\la$}

\begin{proposition}  \label{prop:HT}
Suppose that an equation (of finite sums) of the form
$\sum_{\alpha} a_\alpha \iV_\alpha =\sum_{\alpha} b_\alpha \iV_\alpha$ holds. Then
$\sum_{\alpha} a_\alpha \iV_{(\mu, \alpha)} =\sum_{\alpha} b_\alpha \iV_{(\mu, \alpha)}$, for any integer vector $\mu$.
\end{proposition}

\begin{proof}
The proof is verbatim the same as for \cite[Proposition 1]{Tam11}.
\end{proof}

\begin{lemma} \label{lem:iQQ}
For $a\ge b \ge 0$, $\iV_{(a,b)}$ has the following expansion:
\begin{align*}
\iV_{(a,b)}
&=  V_{(a,b)} + (t-1) \sum_{s=1}^{b} \vth^s t^{s-1} V_{(a-s,b-s)}.
\end{align*}
\end{lemma}

\begin{proof}
We compute
\begin{align*}
\iV_{(a,b)}
&=  \frac{1 -R_{12}}{1 -tR_{12}} \cdot \frac{1 -\vth L_{12}}{1 - \vth t L_{12}}  \vv_{(a,b)}
\\
&=  \frac{1 -R_{12}}{1 -tR_{12}} \Big(1 + (1-t^{-1}) \sum_{s\ge 1}  \vth^s t^s L_{12}^s \Big) \vv_{(a,b)}
\\
&=  \frac{1 -R_{12}}{1 -tR_{12}} \Big(\vv_{(a,b)} + (1-t^{-1}) \sum_{s=1}^{\min(a,b)} \vth^s t^s \vv_{(a-s,b-s)} \Big)
\\
&=  V_{(a,b)} + (1-t^{-1}) \sum_{s=1}^{\min(a,b)} \vth^s t^s V_{(a-s,b-s)}.
\end{align*}
The lemma is proved.
\end{proof}

\begin{lemma}  \label{lem:straight}
For any $a, b \in \Z$, we have
\[
\iV_{(a,b)} +\iV_{(b-1,a+1)} =t (\iV_{(a+1,b-1)} +\iV_{(b,a)}).
\]
In particular, we have $\iV_{(a,a+1)} =t \iV_{(a+1,a)}$.
\end{lemma}

\begin{proof}
Observe that $(1-R_{12}) (\vv_{(a,b)} +\vv_{(b-1,a+1)})=0$, and
\[
\frac{1 -R_{12}}{1 -tR_{12}}
=  (1 -R_{12}) + t\frac{1 -R_{12}}{1 -tR_{12}} R_{12}.
\]
Hence we have
\begin{align*}
\iV_{(a,b)} +\iV_{(b-1,a+1)}
&= \frac{1 -\vth L_{12}}{1 -\vth t L_{12}} \cdot \frac{1 -R_{12}}{1 -tR_{12}}   (\vv_{(a,b)} +\vv_{(b-1,a+1)})
\\
&= \frac{1 -\vth L_{12}}{1 -\vth t L_{12}} \cdot \left( (1 -R_{12}) + t\frac{1 -R_{12}}{1 -tR_{12}} R_{12} \right)  (\vv_{(a,b)} +\vv_{(b-1,a+1)})
\\
&= \frac{1 -\vth L_{12}}{1 -\vth t L_{12}} \cdot   t\frac{1 -R_{12}}{1 -tR_{12}} R_{12}   (\vv_{(a,b)} +\vv_{(b-1,a+1)})
\\
&= t (\iV_{(a+1,b-1)} +\iV_{(b,a)}).
\end{align*}
The formula $\iV_{(a,a+1)} =t \iV_{(a+1,a)}$ follows from the general formula by setting $b=a+1$.
\end{proof}

%


\begin{lemma}
 \label{lem:cd1}
For $c \in \Z, d\ge 1$, we have
\begin{align}  \label{eq:ccd1}
\iV_{(c, c+d)} + (1-t) \sum_{i=1}^{d-1} \iV_{(c+i, c+d-i)} =t\; \iV_{(c+d, c)}.
\end{align}
\end{lemma}

\begin{proof}
We prove the formula by induction on $d$. The case for $d=1$ follows by Lemma~\ref{lem:straight} with $(a,b)=(c,c+1)$, while the case for $d=2$ follows by Lemma~\ref{lem:straight} with $(a,b)=(c,c+2)$.
Let $d>2$. by the inductive assumption,
$\iV_{(c+1, c+d-1)}  + (1-t) \sum_{i=1}^{d-3} \iV_{(c+1+i, c+d-1-i)}= t \iV_{(c+d-1,c+1)},
$
which can be rewritten as
$$\iV_{(c+1, c+d-1)} - t \iV_{(c+d-1,c+1)} + (1-t) \sum_{i=2}^{d-2} \iV_{(c+i, c+d-i)}=0.
$$
By Lemma~\ref{lem:straight} (with $a=c, b=c+d$), we have
\[
\iV_{(c,c+d)} - t \iV_{(c+1, c+d-1)} + \iV_{(c+d-1,c+1)} = t \iV_{(c+d,c)}.
\]
Adding the above 2 identities we obtain the identity \eqref{eq:ccd}.
\end{proof}

\begin{lemma}
 \label{lem:cd}
For $c \in \Z, d\ge 1$, and any integer vector $\beta$, we have
\begin{align}  \label{eq:ccd}
\iV_{(c, c+d, \beta)} + (1-t) \sum_{i=1}^{d-1} \iV_{(c+i, c+d-i, \beta)} =t\; \iV_{(c+d, c,\beta)}.
\end{align}
\end{lemma}

\begin{proof}
This is a variant of \cite[Corollary~1]{Tam11}, and we provide some detailed argument.
Let $\beta \in \Z^\ell$, and we shall prove the lemma by induction on $\ell$. The case for $\ell=0$ follows by Lemma~\ref{lem:cd1}. For $\ell \ge 1$, write $\beta =(\beta',b)$. It follows by \eqref{eq:recur} that
\[
\iV_{\delta, b}  = \sum_{\rho, \gamma} \vth^{|\rho|} t^{|\rho|+|\gamma|} (1-t^{-1})^{\ell(\rho) +\ell(\gamma)}\; \iV_{\delta -\rho +\gamma} \vv_{b-|\rho| -|\gamma|}
\]
where we take $\delta =(c+i,c+d-i, \beta')$, for $0\le i\le d$. The desired formula \eqref{eq:ccd} is then reduced with help of these formulas to the inductive assumption.
\end{proof}

\subsection{Mirror identities}

Given partitions $\mu, \nu$ and $b\in \N$, we use
\[
\nu \rightarrow \mu,
\qquad (\text{respectively, } \nu \stackrel{b}{\rightarrow} \mu)
\]
to denote that $\nu \leq \mu$ and $\mu -\nu$ is a horizontal strip (respectively, a horizontal $b$-strip). Denote by $m_i(\mu)$ the number of times $i$ occurs as a part of $\mu$. Denote by $\mu' =(\mu_1', \mu_2', \ldots)$ the conjugate partition of $\mu$.

For a horizontal $r$-strip $\sigma=\la-\nu$,  let $I=I_{\la-\nu}$ (respectively, $J=J_{\la-\nu}$) be the set of integers $i\geq 1$ such that $\sigma_i'=1$ and $\sigma_{i+1}'=0$ (respectively, $\sigma_i'=0$ and $\sigma_{i+1}'=1$). Following \cite[III]{Mac95}, we introduce
\begin{align}
  \label{eq:phipsi}
  \begin{split}
\varphi_{\la/\nu}(t) &=\prod\limits_{i\in I_{\la-\nu}}(1-t^{m_i(\la)}),
\\
\psi_{\la/\nu}(t) &=\prod\limits_{j\in J_{\la-\nu}}(1-t^{m_j(\nu)}).
\end{split}
\end{align}

Our horizontal Pieri rule in Theorem~\ref{thm:Pieri-Ho} below is formulated using the polynomials $\varphi, \psi$ above. To prove that, we shall require  2 mirror identities. The first mirror identity below is due to Tamvakis \cite[\S2.2]{Tam11} when $\vth=0$  (i.e., with $\iV$ replaced  by $V$), though it requires a little work by the reader to unravel it. 

Denote by $\N$ the set of non-negative integers. For $\gamma =(\gamma_1, \ldots, \gamma_{\ell+1}) \in \N^{\ell+1}$, we set
\[
\gamma_{\bullet} :=(\gamma_1, \ldots, \gamma_\ell).
\]

\begin{proposition} [Mirror identity I] (also see \cite[\S2.2]{Tam11})
For a partition $\rho \in \N^\ell$ and an integer $b \ge 0$, we have
\begin{align}
\sum_{\gamma  \in \N^{\ell+1}, |\gamma| = b} (1-t)^{\ell( \gamma_{\bullet})}  \;  \iV_{\rho +\gamma}
&= \sum_{\rho \stackrel{\scriptscriptstyle b}{\rightarrow} \la} \psi_{\la/\rho}(t) \iV_\la.
  \label{iQ:psi}
 \end{align}
\end{proposition}

\begin{proof}
The strategy of proof here is entirely due to Tamvakis \cite[\S2.2]{Tam11}  (and by personal communication), who gave an elegant argument of the Pieri rule for Hall-Littlewood functions using raising operators (which corresponds to our special case when setting $\vth=0$). Below we clarify the exposition therein by adding more intermediate steps in more careful notations (for instance, the $\gamma$ {\em loc. cit.} at various places should be replaced by $\gamma_{\bullet}$ here).

Setting $\nu =\rho +\gamma$, we rewrite
\begin{align*}
\sum_{\gamma  \in \N^{\ell+1}, |\gamma|=b} (1-t)^{\ell( \gamma_{\bullet})}  \; \iV_{\rho +\gamma}
&= \sum_{\nu  \in \cN} (1-t)^{\ell(\nu_{\bullet} -\rho)}  \; \iV_{\nu},
\end{align*}
where $\cN :=\cN(\rho, b) =\{ \nu \in \N^{\ell+1} \mid |\nu| =|\rho|+b, \nu \geq \rho \}$.
The identity \eqref{iQ:psi} is then equivalent to the following identity: for a partition $\rho =(\rho_1, \rho_2, \ldots, \rho_\ell) \in \N^\ell$ and $b \ge 0$,
\begin{align}
 \label{eq:main}
\sum_{\nu  \in \cN} (1-t)^{\ell(\nu_{\bullet} -\rho)}  \; \iV_{\nu}
=\sum_{\rho \stackrel{\scriptscriptstyle b}{\rightarrow} \la} \psi_{\la/\rho}(t) \ \iV_\la.
\end{align}

We shall prove \eqref{eq:main} by induction on $\ell$, with the case for $\ell=0$ being trivial.
Set
\[
\rho^*=(\rho_2, \ldots, \rho_\ell) \in \N^{\ell-1}.
\]
By inductive assumption we have, for $0\le b^* \le b$,
\begin{align}
 \label{eq:ind}
\sum_{\nu^*  \in \bar\cN}
(1-t)^{\ell(\nu^*_{\bullet} -\rho^*)}  \; \iV_{\nu^*}
=\sum_{\rho^* \stackrel{\scriptscriptstyle b^*}{\rightarrow} \la^*} \psi_{\la^*/\rho^*}(t) \ \iV_{\la^*},
\end{align}
where $\bar\cN :=\cN(\rho^*, b^*) =\{ \nu^* \in \N^{\ell} \mid |\nu^*| =|\rho^*|+b^*, \nu^* \geq \rho^* \}$, and $\nu^*_{\bullet} \in \N^{\ell-1}$ is obtained from $\nu^*$ with the $\ell$th part removed.

By Proposition~\ref{prop:HT}, we can change simultaneously $\iV_{\nu^*}$ to $\iV_{\nu_1,\nu^*}$ and $\iV_{\la^*}$ to $\iV_{\nu_1,\la^*}$ on  \eqref{eq:ind}:
\begin{align*}
\sum_{\nu^*  \in \bar\cN}
(1-t)^{\ell(\nu^*_{\bullet}  -\rho^*)}  \; \iV_{\nu_1,\nu^*}
=\sum_{\rho^* \stackrel{\scriptscriptstyle b^*}{\rightarrow} \la^*} \psi_{\la^*/\rho^*}(t) \ \iV_{\nu_1,\la^*}.
\end{align*}
Now multiply both sides of this new identity by $(1-t)^{\ell(\nu_1 -\rho_1)}$ with $\nu_1 -\rho_1=b-b^* \ge 0$, and then sum over all $b^*$ with $0\le b^* \le b$. In this way, with $\nu =(\nu_1, \nu^*)$ we obtain the identity
\begin{align}
 \label{eq:ind2}
\sum_{\nu \in \cN}
(1-t)^{\ell (\nu_{\bullet}  -\rho)}  \; \iV_{\nu}
=\sum_{(\nu_1, \la^*) \in \cN'} (1-t)^{\ell (\nu_1 -\rho_1)} \psi_{\la^*/\rho^*}(t) \ \iV_{\nu_1,\la^*},
\end{align}
where
\[
\cN' = \big\{(\nu_1, \la^*) \in \cN(\rho,b) ~\big |~  \rho^* \rightarrow \la^*,  \ \nu_1 -\rho_1=b-|\la^*-\rho^*| \big \}.
\]

Compare \eqref{eq:main} with \eqref{eq:ind2}. To complete the proof of \eqref{eq:main} it remains to prove that, for a partition $\rho \in \N^\ell$ and $b \ge 0$,
\begin{align}
 \label{eq:ind3}
\sum_{(\nu_1, \la^*) \in \cN'} (1-t)^{\ell (\nu_1 -\rho_1)} \psi_{\la^*/\rho^*}(t) \ \iV_{\nu_1,\la^*}
=\sum_{\rho \stackrel{\scriptscriptstyle b}{\rightarrow} \la} \psi_{\la/\rho}(t) \ \iV_\la.
\end{align}
Note that $\rho \stackrel{\scriptscriptstyle b}{\rightarrow} \la$ on RHS\eqref{eq:ind3} requires
 $\la_2 \le \rho_{1}$, while there is no similar constraint on LHS\eqref{eq:ind3} between $\la_2$ and $\rho_1$. We shall divide the sum on LHS\eqref{eq:ind3} into 2 partial sums, depending on whether $\la_2 < \rho_1$ or $\la_2 \geq \rho_1$.

The first partial sum with $\la_2 < \rho_1$ on LHS\eqref{eq:ind3} is easy to handle: with an identification of notation $\nu_1= \la_1$, we have $(1-t)^{\ell (\nu_1 -\rho_1)} \psi_{\la^*/\rho^*}(t) =\psi_{\la/\rho}(t)$ by inspection of $\psi$ in \eqref{eq:phipsi}. (The reader is recommended to draw some diagram for $\la/\rho$ here and below.)
Hence,
\begin{align*}
\sum_{(\nu_1, \la^*) \in \cN',\, \la_2 < \rho_1} (1-t)^{\ell (\nu_1 -\rho_1)} \psi_{\la^*/\rho^*}(t) \ \iV_{\nu_1,\la^*}
=\sum_{\rho \stackrel{\scriptscriptstyle b}{\rightarrow} \la, \, \la_2 < \rho_1} \psi_{\la/\rho}(t) \ \iV_\la.
\end{align*}

To show the remaining partial sums on both sides of \eqref{eq:ind3} are equal, we make a more precise claim. Write $\la^* =(\la_2, \la_3, \ldots)$.

{\bf Claim ($\texttt R$).} For each $d\ge 0$, we have
\begin{align}
 \label{eq:partial2}
\sum_{\stackrel{(\nu_1, \la^*) \in \cN',\, \la_2 \ge \rho_1}{\nu_1+\la_2 =2\rho_1+d}} (1-t)^{\ell (\nu_1 -\rho_1)} \psi_{\la^*/\rho^*}(t) \ \iV_{\nu_1,\la^*}
=\sum_{\stackrel{\rho \stackrel{\scriptscriptstyle b}{\rightarrow} \la}{\la_1=\rho_1+d, \, \la_2 =\rho_1}} \psi_{\la/\rho}(t) \ \iV_\la.
\end{align}
The Claim is clear for $d=0$, and we shall assume $d\ge 1$ below.

More precisely, we shall establish the identity \eqref{eq:partial2} for fixed $(\la_3, \la_4, \ldots)$. This is then reduced  by Lemma~\ref{lem:cd} to proving a version of \eqref{eq:partial2} where $\iV_{\nu_1,\la^*}$ (and respectively, $\iV_{\la}$) is replaced by $\iV_{(\nu_1,\la_2)}$ (and respectively, $\iV_{(\la_1, \la_2)}$).

By unravelling the definition of $\psi$ in \eqref{eq:phipsi}, we see that this ``two 2-row" reduced version of the identity \eqref{eq:partial2} becomes the next 2 identities.
First, we have the identity
\begin{align}
 \label{eq:cd1}
 \iV_{(c,c+d)} +(1-t) \sum_{i=1}^{d-1}  \iV_{(c+i,c+d-i)} +(1-t)  \iV_{(c+d,c)} = \iV_{(c+d,c)},
\end{align}
when $\la -\rho$ has a box in column $c=\rho_1$; otherwise, we have the identity
\begin{align}
 \label{eq:cd2}
 (1-t^m) \iV_{(c,c+d)} +(1-t)(1-t^m) \sum_{i=1}^{d-1}  \iV_{(c+i,c+d-i)} +(1-t)  \iV_{(c+d,c)} = (1-t^{m+1})\iV_{(c+d,c)},
\end{align}
for $m\ge 1$. (We have set $m=m_{\rho_1} (\rho)$, the multiplicity of parts in $\rho$ equal to $\rho_1$.) These 2 identities \eqref{eq:cd1}--\eqref{eq:cd2} follow directly from Lemma~\ref{lem:cd}.

This completes the proof of Claim ($\texttt R$), the identity \eqref{eq:ind3}, and whence \eqref{iQ:psi}.
\end{proof}

The second mirror identity below at $\vth=0$ (i.e., with $\iV$ replaced by $V$)  and $b=0$ appeared first in \cite[Theorem 1]{Tam11}, with a very different proof. In our context, this identity is a lowering operator counterpart to the mirror identity \eqref{iQ:psi}; this shall become clear shortly in the proof of Theorem~\ref{thm:Pieri-Ho}.

Recall $\varphi_{\la/\nu}(t)$ from \eqref{eq:phipsi} and $\ell(\alpha)$ from \eqref{eq:number}.

\begin{proposition} [Mirror identity II] (see \cite[Theorem 1]{Tam11})
For a partition $\mu \in \N^\ell$ and integers $a,b \ge 0$, we have
\begin{align}
\sum_{\beta \in \N^\ell, |\beta|=a} (1-t)^{\ell( \beta)}  \;   \iV_{\mu -\beta,b}
&= \sum_{\nu \stackrel{a}{\rightarrow} \mu}  \varphi_{\mu/\nu}(t)   \iV_{\nu,b}.
 \label{iQ:phi}
 \end{align}
\end{proposition}

\begin{proof}
 The arguments here are analogous and ``dual" to those for \eqref{iQ:psi}; besides we have adjoined $b$ to the end of indices for $\iV$ throughout.

 Setting $\eta =\mu-\beta$, we rewrite
\begin{align*}
\sum_{\beta \in \N^\ell, |\beta|=a} (1-t)^{\ell( \beta)}  \;   \iV_{\mu -\beta,b}
&=  \sum_{\eta  \in \cm} (1-t)^{\ell (\mu -\eta)}  \; \iV_{\eta,b},
 \end{align*}
where $\cm :=\cm(\mu, a) =\{ \eta \in \N^{\ell} \mid |\eta| =|\mu| -a, \mu \geq \eta \}.$
The identity \eqref{iQ:phi} is then equivalent to the following identity: for a partition $\mu =(\mu_1, \mu_2, \ldots, \mu_\ell) \in \N^\ell$ and $a \ge 0$,
\begin{align}
 \label{eq:mainII}
\sum_{\eta  \in \cm} (1-t)^{\ell(\mu -\eta)}  \; \iV_{\eta,b}
&= \sum_{\nu \stackrel{a}{\rightarrow} \mu}  \varphi_{\mu/\nu}(t)   \iV_{\nu,b}.
\end{align}
We shall prove \eqref{eq:mainII} by induction on $\ell$, with the case for $\ell=0$ being trivial.

Set
\[
\mu^*=(\mu_2, \ldots, \mu_\ell)  \in \N^{\ell-1}.
\]
By inductive assumption we have, for $0\le a^* \le a$,
\begin{align}
 \label{eq:indII}
\sum_{\eta^*  \in \bar\cm}
(1-t)^{\ell (\mu^* -\eta^*)}  \; \iV_{\eta^*,b}
=\sum_{\nu^* \stackrel{a^*}{\rightarrow} \mu^*}  \varphi_{\mu^*/\nu^*}(t)   \iV_{\nu^*,b},
\end{align}
where $\bar\cm :=\cm(\mu^*, a^*) =\{ \eta^* \in \N^{\ell-1} \mid |\eta^*| =|\mu^*| -a^*, \mu^* \geq \eta^* \}.$

By Proposition~\ref{prop:HT}, we can change simultaneously $\iV_{\eta^*,b}$ to $\iV_{\eta_1,\eta^*,b}$ and $\iV_{\nu^*,b}$ to $\iV_{\eta_1,\nu^*,b}$ on the two sides of \eqref{eq:indII}:
\begin{align*}
\sum_{\eta^*  \in \bar\cm}
(1-t)^{\ell (\mu^* -\eta^*)}  \; \iV_{\eta_1,\eta^*,b}
=\sum_{\nu^* \stackrel{a^*}{\rightarrow} \mu^*}  \varphi_{\mu^*/\nu^*}(t)   \iV_{\eta_1,\nu^*,b}.
\end{align*}
Now multiply  two sides of the above new identity by $(1-t)^{\ell(\mu_1 -\eta_1)}$ with $\mu_1 -\eta_1 =a-a^* \ge 0$, and then sum over all $0\le a^* \le a$. In this way, with $\eta = (\eta_1,\eta^*)$, we obtain the identity
\begin{align}
 \label{eq:ind2II}
\sum_{\eta  \in \cm} (1-t)^{\ell (\mu -\eta)}  \; \iV_{\eta,b}
&= \sum_{(\eta_1, \nu^*) \in \cm'} (1-t)^{\ell(\mu_1 -\eta_1)} \varphi_{\mu^*/\nu^*}(t)   \iV_{\eta_1,\nu^*,b},
\end{align}
where
\[
\cm' =\big\{ (\eta_1, \nu^*) \in \cm(\mu,a) ~\big |~ \nu^* \rightarrow \mu^*, \mu_1 -\eta_1= a-| \mu^* -\eta^*|
\big \}.
\]

Compare \eqref{eq:mainII} with \eqref{eq:ind2II}. To complete the proof of \eqref{eq:mainII} it remains to prove that, for a partition $\mu \in \N^\ell$ and $a \ge 0$,
\begin{align}
 \label{eq:ind3II}
\sum_{(\eta_1, \nu^*) \in \cm'} (1-t)^{\ell(\mu_1 -\eta_1)} \varphi_{\mu^*/\nu^*}(t)   \iV_{\eta_1,\nu^*,b}
&= \sum_{\nu \stackrel{a}{\rightarrow} \mu}  \varphi_{\mu/\nu}(t)   \iV_{\nu,b}.
\end{align}
Note that $\nu \stackrel{a}{\rightarrow} \mu$ on RHS\eqref{eq:ind3II} requires
$\mu_2 \le \nu_1$, while there is no similar constraint on LHS\eqref{eq:ind3II} between $\mu_2$ and $\eta_1$. 
We shall divide the sum on LHS\eqref{eq:ind3II} into 2 partial sums, depending on whether $\mu_2 <\eta_1$ or $\mu_2 \geq \eta_1$.

The first partial sum with $\mu_2 <\eta_1$ on LHS\eqref{eq:ind3II} is easy to handle: with an identification of notation $\eta_1 =\nu_1$, we have $(1-t)^{\ell (\mu_1 -\nu_1)} \varphi_{\mu^*/\nu^*}(t) =\varphi_{\mu/\nu}(t)$ by inspection of $\varphi$ in \eqref{eq:phipsi}, and thus
\begin{align*}
 \sum_{(\eta_1, \nu^*) \in \cm',\, \mu_2 <\eta_1} (1-t)^{\ell(\mu_1 -\eta_1)} \varphi_{\mu^*/\nu^*}(t)   \iV_{\eta_1,\nu^*,b}
&= \sum_{\nu \stackrel{a}{\rightarrow} \mu,\, \mu_2 <\eta_1}  \varphi_{\mu/\nu}(t)   \iV_{\nu,b}.
 \end{align*}

To show the remaining partial sums on both sides of \eqref{eq:ind3II} are equal, we make a more precise claim. We write $\nu^*=(\nu_2, \nu_3, \ldots)$.

{\bf Claim ($\texttt L$).} For each $d\ge 0$, we have
\begin{align}
 \label{eq:partial2II}
 \sum_{\stackrel{(\eta_1, \nu^*) \in \cm', \, \mu_2 \geq \eta_1}{\eta_1+\nu_2 =2\mu_2-d}} (1-t)^{\ell(\mu_1 -\eta_1)} \varphi_{\mu^*/\nu^*}(t)   \iV_{\eta_1,\nu^*,b}
&= \sum_{\stackrel{\nu \stackrel{a}{\rightarrow} \mu}{\nu_1=\mu_2, \, \nu_2=\mu_2-d}  }
\varphi_{\mu/\nu}(t)   \iV_{\nu,b}.
\end{align}
 The Claim is clear for $d=0$, and we shall assume $d\ge 1$ below.

More precisely, we shall establish the identity \eqref{eq:partial2II} for fixed $(\nu_3, \nu_4, \ldots)$. This is then reduced  by Lemma~\ref{lem:cd} to proving a version of \eqref{eq:partial2II} where $\iV_{\eta_1,\nu^*,b}$ (and respectively, $\iV_{\nu,b}$) is replaced by $\iV_{(\eta_1,\nu_2,b)}$ (and respectively, $\iV_{(\nu_1, \nu_2,b)}$).

By unravelling the definition of $\varphi$ in \eqref{eq:phipsi}, we see that this ``top 2-row" reduced version of the identity \eqref{eq:partial2II} becomes the next 2 identities:
\begin{equation}
\begin{cases}
 \scriptstyle
 (1-t)\iV_{(c-d,c)} +(1-t)^2 \sum_{i=1}^{d-1}  \iV_{(c-i,c+i-d)} +(1-t)  \iV_{(c,c-d)} = (1-t^2) \iV_{(c,c-d)},
& \text{ if } c=\mu_2=\mu_1,
\\
\scriptstyle (1-t) \iV_{(c-d,c)} +(1-t)^2 \sum_{i=1}^{d-1}  \iV_{(c-i,c+i-d)} +(1-t)^2  \iV_{(c,c-d)} = (1-t)  \iV_{(c,c-d)},
 & \text{ if } c= \mu_2 < \mu_1.
 \end{cases}
\end{equation}
These 2 identities follow directly from Lemma~\ref{lem:cd}.

This completes the proof of Claim ($\texttt L$), the identity \eqref{eq:ind3II}, and whence \eqref{iQ:phi}.
\end{proof}

\subsection{Horizontal $\imath$Pieri rule}


As we shall see, the following Pieri rule corresponds to the multiplication formula \eqref{Pieri:h} in the $\imath$Hall algebra setting.

\begin{theorem}  [Down-up horizontal $\imath$Pieri rule] 
\label{thm:Pieri-Ho}
For $r\ge 1$ and any partition $\mu$, we have
\begin{align}
\iV_\mu  \cdot \vv_r
&=\sum_{a+b=r} \sum_{\nu \stackrel{a}{\rightarrow} \mu ,\nu\stackrel{b}{\rightarrow} \la} \vth^{a}\, \varphi_{\mu/\nu}(t) \psi_{\la/\nu}(t) \;  \iV_\la.
 \label{Pieri:hSF3}
\end{align}
\end{theorem}

\begin{proof}
Set $\ell =\ell(\mu)$. Recall from \eqref{eq:L}--\eqref{eq:D} that $\L_{\ell+1}^{-1} =1 +(1-t) \sum_{k\ge 1} \vth^k L_{i,\ell+1}^k$ and $D_\ell =\R_{\ell+1}^{-1} \cdot D_{\ell+1} \cdot \L_{\ell+1}^{-1}$. Then we have
\begin{align}
\iV_\mu  \cdot \vv_r
&= D_\ell \vv_{\mu,r}
 \label{iQ:r1} \\
&= \R_{\ell+1}^{-1}
\cdot D_{\ell+1}
\prod_{1\le i \le \ell} (1 +(1-t) \sum_{k\ge 1} \vth^k L_{i,\ell+1}^k) \vv_{\mu,r}
\notag \\
&= \R_{\ell+1}^{-1}
\cdot D_{\ell+1}\cdot
 \sum_{\beta \in \N^\ell} \vth^{|\beta|}  (1-t)^{\ell( \beta)}  \;   \vv_{\mu -\beta, r-|\beta|}
\notag
\\
&= \R_{\ell+1}^{-1}
\cdot
 \sum_{\beta \in \N^\ell} \vth^{|\beta|}  (1-t)^{\ell( \beta)}  \;   \iV_{\mu -\beta, r-|\beta|}
\notag
\\
&= \sum_a \sum_{\nu \stackrel{a}{\rightarrow} \mu} \vth^a \varphi_{\mu/\nu}(t)
\R_{\ell+1}^{-1}  \iV_{\nu, r-a},
 \notag
\end{align}
where we have used the mirror identity \eqref{iQ:phi} with $b=r-a$ in the last step.

Then,  a simple but key next step is to convert $\iV_{\nu, r-a}$ back to $\vv_{\nu, r-a}$ and work with the raising/lowering operators. As $D_{\ell+1}$ and $\R_{\ell+1}^{-1}$ commute, we rewrite the RHS of \eqref{iQ:r1} using $\iV_{\nu, r-a} = D_{\ell+1}\vv_{\nu, r-a}$ as
\begin{align}
  \label{iQ:r2}
\iV_\mu  \cdot \vv_r
&= \sum_a \sum_{\nu \stackrel{a}{\rightarrow} \mu} \vth^a \varphi_{\mu/\nu}(t) \,
D_{\ell+1} \cdot \R_{\ell+1}^{-1}   \vv_{\nu, r-a}
 \\
&= \sum_a \sum_{\nu \stackrel{a}{\rightarrow} \mu} \vth^a \varphi_{\mu/\nu}(t) \,
D_{\ell+1} \cdot \prod_{1\le i \le \ell} (1 +(1-t) \sum_{k\ge 1} R_{i,\ell+1}^k)  \vv_{\nu, r-a} %
\notag \\
&= \sum_a \sum_{\nu \stackrel{a}{\rightarrow} \mu} \vth^a \varphi_{\mu/\nu}(t)
D_{\ell+1}
\sum_{\gamma  \in \N^{\ell+1}, |\gamma| =r-a} (1-t)^{\ell( \gamma_{\bullet})} \vv_{\nu +\gamma}
\notag \\
&= \sum_a \sum_{\nu \stackrel{a}{\rightarrow} \mu} \vth^a \varphi_{\mu/\nu}(t)
\sum_{\gamma  \in \N^{\ell+1}, |\gamma| =r-a} (1-t)^{\ell( \gamma_{\bullet})} \iV_{\nu +\gamma}.
\notag
\end{align}

The theorem follows now by applying the mirror identity \eqref{iQ:psi} (where $\rho$ and $b$ are replaced by $\nu$ and $r-a$ here) to the RHS of  \eqref{iQ:r2}.
\end{proof}

\begin{remark}
 The Pieri rule in \eqref{Pieri:hSF3} was derived from the Giambelli formula \eqref{eq:LRq} for $\iV_\la$. Conversely (by a simple induction), one can derive the Giambelli formula from the Pieri rule.
 \end{remark}

\subsection{$\imath$HL functions}

We recall the Hall-Littlewood (HL) functions in variables $x=(x_1, x_2, \ldots)$, following \cite[Pages 208--211]{Mac95}. Denote by $h_r$ the $r$th homogeneous symmetric function, for $r\ge 0$.
For an indeterminate $u$, denote
\begin{align*}
H(u)&= \sum_{r=0}^\infty h_r u^r = \prod_{i\ge 1} \frac1{1-ux_i},
\qquad
Q(u) = \sum_{r=0}^\infty Q_{r} u^r=H(u)/H(tu).
\end{align*}

We set $q_r =Q_r$, and $q_\la =q_{\la_1} q_{\la_2} \cdots$, for any composition $\la$.
Let $u_1,u_2,\cdots$ be independent indeterminates. Recall $F(u)=(1-u)/(1-tu)$ from \eqref{eq:F}. Then the HL symmetric function $Q_\la$ is the coefficient of $u^\la=u_1^{\la_1} u_2^{\la}\cdots$ in
\begin{align}
\label{eq:HL}
Q(u_1,u_2,\dots)= \prod_{i\geq1}Q(u_i) \prod_{i<j} F(u_iu_j^{-1}).
\end{align}
According to \cite[p.212]{Mac95}, the generating function \eqref{eq:HL} for HL functions $Q_\la$ can be restated that
\begin{align}  \label{eq:Rq}
Q_\la = \prod_{i<j} \frac{1 -R_{ij}}{1 -tR_{ij}} q_\la
=\Big(1 + (t-1) \sum_{r\ge 1} t^{r-1}  R_{ij}^r  \Big) q_\la.
\end{align}

%

In this subsection, we shall set $\vv_i = Q_{r}$ in the definition \eqref{eq:LRq} of $V_\la, \iV_\la$, for any composition $\la$. In this way, we identify the ring $\La_{t,\vth} = \Q(t)[\vth] [q_1, q_2, \ldots]$ with the ring of symmetric functions.

Then we see from  \eqref{eq:LRq} and \eqref{eq:Rq} that $V_\la$ coincides with $Q_\la$.
On the other hand, $\iV_\la$ gives us another class of symmetric functions, which are called {\em $\imath$HL functions} and denoted by $\iQ_\la$. Alternatively, in the spirit of the equivalence between \eqref{eq:HL} and \eqref{eq:Rq}, we can define $\iQ_\la$ as follows.

\begin{definition} [$\imath$HL functions]
 \label{def:iQ}
For any composition $\la$, $\iQ_\la$ is the coefficient of $u^\la=u_1^{\la_1}u_2^{\la_2}\cdots$ in
\begin{align}
\label{HL functions}
\iQ(u_1,u_2,\dots)=&\prod_{i\geq1} Q(u_i)\prod_{i<j} F(u_i^{-1}u_j)F(\vth  u_iu_j)
\\\notag
=&\prod_{i\geq1} Q(u_i) \prod_{i<j} \frac{(1-u_i^{-1}u_j)(1- \vth  u_iu_j)}{(1- tu_i^{-1}u_j)(1-t \vth  u_iu_j)}.
\end{align}
\end{definition}

\begin{remark}
 \label{rem:Okada}
Okada \cite[\S3.1]{Ok21} introduced the universal symplectic Q-functions $\mathcal Q^C_\lambda$, which are parametrized by strict partitions $\lambda$. We learned from Okada  that a specialization of the iHL functions $Q^\imath_\lambda$ (for $\la$ strict) coincides with the universal symplectic Q-functions (compare \eqref{HL functions} with \cite[(2.14), (3.7)-(3.8)]{Ok21}):
\[
Q^\imath_\lambda|_{t=-1,\theta=1} = \mathcal Q^C_\lambda.
\]
\end{remark}

\subsection{Modified $\imath$HL functions}

Let $x=(x_1, x_2, \ldots)$. Recall the modified Hall-Littlewood functions $H_\la(x; t)$ can be defined using plethysm, $H_\la(x; t) =Q(x/(1-t); t)$, cf. \cite[Ex.~7, p.234]{Mac95}, where it was denoted by $Q_\la'$. It can also be defined via Garsia's version of Jing's vertex operators. 
Via raising operators the modified HL functions are defined to be
\begin{align}  \label{eq:HRh}
H_\mu = \prod_{1\le i<j}  \frac{1 -R_{ij}}{1 -tR_{ij}} h_\mu.
\end{align}
Note $H_{(r)} =h_r$.

In this subsection, we shall take $\vv_r =h_r$, the $r$th homogenous symmetric function, in the definition \eqref{eq:LRq}, and then redenote $\iV_\la$ as $\iH_\la$.
In this way, we identify $\La_{t,\vth} = \Q(t)[\vth] [h_1, h_2, \ldots]$ with the ring of symmetric functions.

\begin{definition}
The modified $\imath$HL function $\iH_\la \in \La_{t,\vth}$, for any composition $\la$, is given by
\begin{align}  \label{eq:LRh2}
\iH_\la = \prod_{1\le i<j} \frac{1 - \vth L_{ij}}{1 - \vth t L_{ij}} \frac{1 -R_{ij}}{1 -tR_{ij}} h_\la.
\end{align}
\end{definition}

\begin{proposition}
$\{\iH_\la \}_{\la \in \mathcal P}$ forms a $\Q(t)[\vth]$-basis for the ring $\La_{t,\vth}$.
\end{proposition}

\begin{proof}
Clearly $\iH_\la \in \La_{t,\vth}$.
Note $\iH_{(r)} = h_r$, and $\iH_\la |_{\vth=0} =H_\la$ by \eqref{eq:HRh}; in other words, $\iH_\la $ is equal to $H_\la$ plus lower degree symmetric functions (or equivalently, higher degree in terms of $\vth$). So the transition matrix expressing $\{\iH_\la\}_{|\la|\le n}$ via $\{H_\la\}_{|\la|\le n}$, for each $n>0$, is uni-triangular.
As $\{H_\la \}_{\la \in \mathcal P}$ forms a $\Q(t)[\vth]$-basis for $\La_{t,\vth}$, so does $\{\iH_\la \}_{\la \in \mathcal P}$.
\end{proof}

Now let us set $t=0$, and consider the ring of symmetric functions $\La_{\vth} := \Q [\vth] [h_1, h_2, \ldots]$. We set $t=0$ in $\iH_\la$ to define the {\em $\imath$Schur functions} $\is_\la \in \Lambda_\vth$ by a Giambelli type formula (compare Definition~ \ref{def:V}):
\begin{align}  \label{eq:LRh}
\is_\la = \prod_{1\le i<j} (1 - \vth L_{ij}) (1 -R_{ij}) h_\la.
\end{align}

The $\is_\la$ is homogeneous of degree $|\la|$ if we set the degree $|\vth | =2$. It we regard $\vth$ as a usual parameter of degree 0, then $\is_\la$ is
inhomogeneous with leading term being the Schur function $s_\la$, i.e., $\is_\la |_{\vth=0} =s_\la$.

We denote
\begin{align}
N_{\mu,r}^\la &= \Big | \big \{\nu \in \cp ~  |~ \nu \rightarrow \mu, \mu \rightarrow \la, \text{ and } |\mu-\nu| +|\la-\nu| =r \} \Big |.
\label{eq:N}
\end{align}
Note that $N_{\mu,r}^\la =0$ unless $ -r \le |\la| -|\mu| \le r$ and $|\la| -|\mu| \equiv r \pmod 2$.
The following Pieri rule for $\imath$Schur functions is an immediate corollary by setting $t=0$ in Theorem~\ref{thm:Pieri-Ho}. (It can also be proved directly as for Theorem~\ref{thm:Pieri-Ho}, using a well-known and simpler mirror identities for $\is_\la$, cf., e.g., \cite[Lemma 2]{Tam12}.)

\begin{corollary}
 \label{cor:Pieri0}
For $r\ge 1$ and any partition $\mu$, we have
\begin{align*}
\is_\mu  \cdot h_r
&= \sum_{a=0}^r \sum_{\la \vdash (|\mu|+r-2a)} \vth^{a}\,  N_{\mu,r}^\la  \is_\la.
\end{align*}
\end{corollary}

\subsection{(Deformed) universal characters}

It turns out that that the modified $\imath$HL function $\iH_\la$ appeared earlier in a very different form. The deformed universal characters of type C, denoted by $\mathbb B_\la^{\dmno}$, was introduced in \cite[\S6.2]{SZ06} in plethystic notation via creation operators.

\begin{proposition} \label{prop:B=H}
We have $\iH_\la \big |_{\vth=1} =\mathbb B_\la^{\dmno}$, for any partition $\la$.
\end{proposition}

\begin{proof}
(Mark Shimozono)
Let $\la \in \N^\ell$ and $Z=(z_1, \ldots, z_\ell)$.
Following \cite{SZ06}, we first set up some notations, with products below running over $1\le i<j \le \ell$:
\begin{align*}
R(Z) &=\prod_{i<j} (1-z_j/z_i),
\qquad
\Omega[-s_{1,1}[Z]] =\prod_{i<j} (1-z_iz_j),
\\
\Omega[ZX] &= \sum_{\beta  \in \N^\ell} z^\beta h_\beta[X],
\qquad
\Pi  
= \prod_{i<j} \frac{1}{(1-tz_j/z_i)(1-tz_iz_j)}.
\end{align*}
Then the definition of $\mathbb B_\la^{\dmno}$ via creation operators in \cite[(78)-(79)]{SZ06} can be quickly converted to
\begin{align*}
\mathbb B_\la^{\dmno} &= R(Z) \cdot \Omega[-s_{1,1}[Z]] \cdot \Pi \cdot \Omega[ZX] \big|_{z^\la}
\\
&= \prod_{i<j} \frac{(1-z_j/z_i)(1-z_iz_j)}{(1-tz_j/z_i)(1-tz_iz_j)} \sum_{\beta  \in \N^\ell} z^\beta h_\beta[X] \big|_{z^\la},
\end{align*}
which is equal to $\iH_\la \big |_{\vth=1}$ in the raising/lowering operator form.
\end{proof}

The universal classical characters of type $C$ in \cite{KT87} are defined for any partition $\la$ (here we follow the notation $s^{\dmno}_\la$ used in \cite[(21)]{SZ06}). Then  by a variant of Proposition~\ref{prop:B=H}, we have
\begin{align}
s^{\dmno}_\la = \is_\la \big |_{\vth=1},
\end{align}
where $s^{\dmno}_\la$ is defined via creation operators \cite[(21)]{SZ06} and   $\is_\la \big |_{\vth=1}$ is defined via raising/lowering operator form  in \eqref{eq:LRh}.

Define the polynomials $d_{\la\mu}(t, \vth)$ and $d_{\la\mu}(t)$ via the expansions
\[
\iH_\mu =\sum_{\la} d_{\la\mu}(t, \vth) \, \is_\la,\qquad
\iH_\mu \big |_{\vth=1} =\sum_{\la} d_{\la\mu}(t) \, \is_\la \big |_{\vth=1}.
\]
We have a remarkable positivity result from \cite[Theorem 18]{SZ06}:
\[
d_{\la\mu}(t) \in \N [t].
\]
By assigning degree $|\vth|=2$, $\iH_\mu$ is homogeneous of degree $|\mu|$ and $\is_\la$  is homogeneous of degree $|\la|$. As $\vth$ is exactly used to make up the difference of degrees, we conclude that
\[
d_{\la\mu}(t, \vth) = \vth^{\frac{|\la|-|\mu|}2} d_{\la\mu}(t).
\]

\begin{remark}
Pieri rules for characters of finite type C were computed in \cite{Su90} and  for the universal characters of classical type can be found in \cite{Ok16}. 
 The universal Pieri rule of type C coincides with the formula in Corollary~\ref{cor:Pieri0} by setting $\vth=1$. Four Pieri rules for universal characters can be obtained from specializations at $t=0$ and $\vth=1$ of Theorems~{\bf A(I)}, {\bf A(II)} and {\bf B}, two of which were well known \cite{Su90, Ok16}. It is curious to see if the remaining two can lead to interesting new forms of Pieri rules for irreducible characters of finite classical type after applying modification rules.
\end{remark}

\section{Generalities on semi-derived Hall algebras}
  \label{sec:Hall}

  In this short section, we review and present a variant of semi-derived Hall algebras and $\imath$Hall algebras, following \cite{Lu19, LW20, LRW20} (also cf. \cite{Br13, Gor2}).

\subsection{Hall algebras}
   \label{subsec:HA}
Let $\ca$ be an essentially small abelian category, linear over a finite field $\bfk=\F_\bq$. Assume that $\ca$ is Hom-finite and $\Ext^1$-finite. 
Given objects $M,N,L\in\ca$, let $\Ext^1(M,N)_L\subseteq \Ext^1(M,N)$ be the subset parameterizing extensions whose middle term is isomorphic to $L$. The {\em Hall algebra} (or {\em Ringel-Hall algebra}) $\ch(\ca)$ is defined to be the $\Q$-vector space with the isoclasses $[M]$ of objects $M \in \ca$ as a basis and multiplication given by (cf., e.g., \cite{Br13})
\begin{align}
\label{eq:mult}
[M]\diamond [N]=\sum_{[L]\in \Iso(\ca)}\frac{|\Ext^1(M,N)_L|}{|\Hom(M,N)|}[L].
\end{align}

Given three objects $X,Y,Z$, the Hall number is defined to be
$$G_{XY}^Z:= |\{L\subseteq Z\mid L \cong Y\text{ and }Z/L\cong X\}|.$$
Denote by $\aut(X)$ the automorphism group of $X$. The Riedtmann-Peng formula reads
\begin{align}
 \label{eq:RP}
G_{XY}^Z= \frac{|\Ext^1(X,Y)_Z|}{|\Hom(X,Y)|} \cdot \frac{|\aut(Z)|}{|\aut(X)|\cdot |\aut(Y)|}.
\end{align}

\subsection{Category of $1$-periodic complexes}

\label{subsec:periodic}

Let $\ca$ be a hereditary abelian category  which is essentially small with finite-dimensional homomorphism and extension spaces. 
 A $1$-periodic complex $X^\bullet$ in $\ca$ is a pair $(X,d)$ with $X\in\ca$ and a differential $d:X\rightarrow X$. A morphism $(X,d) \rightarrow (Y,e)$ is given by a morphism $f:X\rightarrow Y$ in $\ca$ satisfying $f\circ d=e\circ f$. The category of all $1$-periodic complexes in $\ca$, denoted by $\cc_1(\ca)$, is an abelian category. Denote by $\cc_{1,ac}(\ca)$ the full subcategory of $\cc_1(\ca)$ consisting of acyclic complexes $X^\bullet=(X,d)$ (here acyclic means $\ker d=\Im d$).
Denote by $H(X^\bullet)\in\ca$ the cohomology group of $X^\bullet$, i.e., $H(X^\bullet)=\ker d/\Im d$.

For any $X\in\ca$, denote the stalk complex by
\[
C_X =(X,0)
\]
 (or just by $X$ when there is no confusion), and denote by $K_X$ the following acylic complex:
\[
K_X:=(X\oplus X, d),
\qquad \text{ where }
d=\left(\begin{array}{cc} 0&\Id \\ 0&0\end{array}\right).
\]

\begin{lemma}
[\text{\cite[Lemma 2.2]{LRW20}}]
\label{lem:pd acyclic}
For any acyclic complex $K^\bullet$ and $p \ge 2$, we have
\begin{align}
\label{Ext2vanish}
\Ext^p_{\cc_1(\ca)}(K^\bullet,-)=0=\Ext^p_{\cc_1(\ca)}(-,K^\bullet).
\end{align}
\end{lemma}

For any $K^\bullet\in\cc_{1,ac}(\ca)$ and $M^\bullet\in\cc_{1}(\ca)$, by \cite[Corollary 2.4]{LRW20}, define
\begin{align*}
\langle K^\bullet, M^\bullet\rangle=\dim_{\bfk} \Hom_{\cc_{1}(\ca)}( K^\bullet, M^\bullet)-\dim_{\bfk}\Ext^1_{\cc_{1}(\ca)}(K^\bullet, M^\bullet),
\\
\langle  M^\bullet,K^\bullet\rangle=\dim_{\bfk} \Hom_{\cc_{1}(\ca)}( M^\bullet, K^\bullet)-\dim_{\bfk}\Ext^1_{\cc_{1}(\ca)}( M^\bullet,K^\bullet).
\end{align*}
These formulas give rise to well-defined bilinear forms (called {\em Euler forms}), again denoted by $\langle \cdot, \cdot \rangle$, on the Grothendieck groups $K_0(\cc_{1,ac}(\ca))$
and $K_0(\cc_{1}(\ca))$.

Denote by $\langle \cdot,\cdot\rangle_\ca$ 
 the Euler form of $\ca$,  i.e.,
\[
\langle M,N\rangle_\ca=\dim_\bfk\Hom_{\ca}(M,N)-\dim_\bfk\Ext^1_{\ca}(M,N).
\]
Let $\res: \cc_1(\ca)\rightarrow\ca$ be the restriction functor.
Let $K_0(\ca)$ be the Grothendieck group of $\ca$. For any object $X\in\ca$, we denote by $\widehat{X}$ its class in $K_0(\ca)$.

\begin{lemma}
[\text{\cite[Lemma 2.7]{LRW20}}]
\label{lemma compatible of Euler form}
We have
\begin{itemize}
\item[(1)]
$\langle K_X, M^\bullet\rangle = \langle X,\res (M^\bullet) \rangle_\ca$,\; $\langle M^\bullet,K_X\rangle =\langle \res(M^\bullet), X \rangle_\ca$, for $X\in \ca$, $M^\bullet\in\cc_1(\ca)$;
\item[(2)] $\langle M^\bullet,N^\bullet\rangle=\frac{1}{2}\langle \res(M^\bullet),\res(N^\bullet)\rangle_\ca$, for  $M^\bullet,N^\bullet\in\cc_{1,ac}(\ca)$.
\end{itemize}
\end{lemma}

\subsection{Semi-derived Ringel-Hall algebras}
\label{sec:iHall}

We continue to work with a hereditary abelian category $\ca$ as in \S\ref{subsec:periodic}.
Let $\ch(\cc_1(\ca))$ be the Ringel-Hall algebra of $\cc_1(\ca)$ over $\Q$, i.e., $\ch(\cc_1(\ca))=\bigoplus_{[X^\bullet]\in \Iso(\cc_1(\ca))} \Q[X^\bullet]$, with multiplication
\begin{align*}
[M^\bullet]\diamond[N^\bullet]=\sum_{[L^\bullet]\in\Iso(\cc_1(\ca)) } \frac{|\Ext^1(M^\bullet,N^\bullet)_{L^\bullet}|}{|\Hom(M^\bullet,N^\bullet)|}[L^\bullet].
\end{align*}

Following \cite{LP21, Lu19, LW20, LRW20},
we consider the ideal $\cI$ of $\ch(\cc_1(\ca))$ generated by
\begin{align}
  \label{eq:ideal}
&\{[K_1^\bullet]-[K_2^\bullet] \mid K_1^\bullet,K_2^\bullet\in\calc_{1,ac}(\ca) \text{ with }\widehat{\Im d_{K_1^\bullet}}=\widehat{\Im d_{K_2^\bullet}}\} \bigcup
\\\notag
&\{[L^\bullet]-[K^\bullet\oplus M^\bullet]\mid \exists \text{ exact sequence } 0 \rightarrow K^\bullet \rightarrow L^\bullet \rightarrow M^\bullet \rightarrow 0 \text{ with }K^\bullet \text{ acyclic}\}.
\end{align}
We denote
\[
\cs:=\{ a[K^\bullet] \in \ch(\cc_1(\ca))/\cI \mid a\in \Q(\sqq)^\times, K^\bullet \text{ acyclic}\},
\]
a multiplicatively closed subset of $\ch(\cc_1(\ca))/ \cI$ with the identity $[0]$.

\begin{lemma}
[\text{\cite[Proposition A.5]{Lu19}}]
The multiplicatively closed subset $\cs$ is a right Ore, right reversible subset of $\ch(\cc_1(\ca))/\cI$. Equivalently, there exists a right localization of
$\ch(\cc_1(\ca))/\cI$ with respect to $\cs$, denoted by $(\ch(\cc_1(\ca))/\cI)[\cs^{-1}]$.
\end{lemma}

The algebra $(\ch(\cc_1(\ca))/\cI)[\cs^{-1}]$ is the {\em semi-derived Ringel-Hall algebra} of $\cc_1(\ca)$ in the sense of \cite{LP21, LW19}, and will be denoted by $\cs\cd\ch(\cc_1(\ca))$.

For any $\alpha\in K_0(\ca)$,  there exist $X,Y\in\ca$ such that $\alpha=\widehat{X}-\widehat{Y}$. Define $[K_\alpha]:=\bq^{-\langle \alpha,\widehat{Y}\rangle_\ca}[K_X]\diamond [K_Y]^{-1}$. This is well defined, see, e.g., \cite[\S 3.2]{LP21}. Let $\ct(\ca)$  be the subalgebra of $\cs\cd\ch(\cc_1(\ca))$ generated by $[K_\alpha]$, for $\alpha\in K_0(\ca)$.

\begin{proposition}
[\text{\cite[Proposition 2.9]{LRW20}}]
\label{prop:hallbasis}
$\cs\cd\ch(\cc_1(\ca))$ is a free left (respectively, right) $\ct(\ca)$-modules with a basis $\{[M]\mid [M]\in\Iso(\ca)\}$. In particular, $\cs\cd\ch(\cc_1(\ca))$ has a basis
\begin{align*}
\{[M]\diamond[K_\alpha]\mid [M]\in\Iso(\ca), \alpha\in K_0(\ca)\}.
\end{align*}
\end{proposition}

For any $M^\bullet=(M,d)$, we have
\begin{align}
\label{eq:decom}
[M^\bullet]=\bq^{\langle H(M^\bullet), \Im d\rangle_\ca}[H(M^\bullet)]\diamond[K_{\widehat{\Im d}}]
\end{align}
 in $\cs\cd\ch(\cc_1(\ca))$; see \cite[Lemma 2.10]{LRW20}.


The following multiplication formula in the semi-derived Ringel-Hall algebra will be useful later.

\begin{proposition}
[\text{\cite[Proposition 3.10]{LW20}}]
\label{prop:iHallmult}
Let $\ca$ be a hereditary abelian category over $\bfk$, and let $A,B\in\ca\subset \cc_1(\ca)$. Then we have in $\cs\cd\ch(\cc_1(\ca))$
\begin{align}
\label{Hallmult}
[A]\diamond[B]=&
\sum_{[M]\in\Iso(\ca)}\sum_{[L],[N]\in\Iso(\ca)}  \bq^{\langle N,L\rangle_\ca -\langle A,B\rangle_\ca}\frac{|\Ext^1(N, L)_{M}|}{|\Hom(N,L)| }
\\
\notag
&\qquad \times \big|\{f:A\rightarrow B\mid \ker f\cong N, \coker f\cong L
\} \big |\cdot [M]\diamond[K_{\widehat{A}-\widehat{N}}].
\end{align}
\end{proposition}

\begin{corollary}
\label{cor:Hallmult}
Let $\ca$ be a hereditary abelian category over $\bfk$, and let $A,B\in\ca\subset \cc_1(\ca)$. Then we have in $\cs\cd\ch(\cc_1(\ca))$
\begin{align}
[A]\diamond[B]=&
\sum_{[M],[L],[N]\in\Iso(\ca)} \bq^{\langle N,L\rangle_\ca -\langle A,B\rangle_\ca} G_{N,L}^M\cdot\frac{|\aut(N)|\cdot |\aut(L)|}{|\aut(M)|}
\\\notag
&\qquad\qquad \times  \Big(\sum_{[X]\in\Iso(\ca)}  G_{X,N}^A G_{L,X}^B|\aut(X)|\Big) [M]\diamond[K_{X}].
\end{align}
\end{corollary}

\begin{proof}
From \eqref{Hallmult}, one obtains that
\begin{align*}
[A]\diamond[B]=&
\sum_{[M],[L],[N],[X]\in\Iso(\ca)}  \bq^{\langle N,L\rangle_\ca -\langle A,B\rangle_\ca}\frac{|\Ext^1(N, L)_{M}|}{|\Hom(N,L)| }
\\
\notag
&\qquad \times \big |\{f:A\rightarrow B\mid \ker f\cong N, \Im f\cong X,\coker f\cong L \} \big |\cdot [M]\diamond[K_{X}].
\end{align*}
A direct computation shows that
\begin{align*}
&\big |\{f:A\rightarrow B\mid \ker f\cong N, \Im f\cong X,\coker f\cong L
\} \big |
\\
=&\big |\{f_1:A\rightarrow X\mid f_1 \text{ is epic},\ker f_1\cong N\} \big | \cdot
\big | \{f_2:X\rightarrow B \mid f_2 \text{ is monic}, \coker f\cong L\} \big |/|\aut(X)|
\\
=&G_{X,N}^AG_{L,X}^B|\aut(X)|.
\end{align*}
The desired formula follows by using the Riedtmann-Peng formula \eqref{eq:RP}.
\end{proof}

\section{$\imath$Hall algebra of the Jordan quiver}
  \label{sec:iHall}

In this section, we formulate the basic properties including 4 Pieri rules for the $\imath$Hall algebra of the Jordan quiver. We also establish an algebra isomorphism from the $\imath$Hall algebra to the ring of symmetric functions.

\subsection{Commutativity}

Let $ \QJ$ be the Jordan quiver, i.e., a quiver with a single vertex $1$ and a single loop $\alpha:1\rightarrow 1$.
Let $\rep^{\rm nil}_\bfk(\QJ)$ be the category of the finite-dimensional nilpotent representations of $\QJ$ over $\bfk$, and $\calc_1(\rep^{\rm nil}_\bfk(\QJ))$ the category of $1$-periodic complexes in $\rep^{\rm nil}_\bfk(\QJ)$.
Let $\Lambda_{\texttt{J}}^\imath$ be the $\imath$quiver algebra of the $\imath$quiver $\QJ$ equipped with trivial involution \cite{LW19, LW20}, which can be identified with
\[
\Lambda_{\texttt{J}}^\imath=\bfk\ov{\QJ} /( \alpha\varepsilon-\varepsilon \alpha, \varepsilon^2 ),
\]
where $\ov{\QJ}$ is enhanced from the Jordan quiver by a second arrow $\varepsilon$ as follows:
\vspace{-4mm}
\begin{equation}
\begin{picture}(100,30)(0,10)
%
\begin{tikzpicture}[scale=1.1]
\draw[color=purple] (7,0) .. controls (7.2,0.35) and (7.6,0.35) .. (7.6,-0.05);
\draw[-latex,color=purple] (7.6,-0.05) .. controls (7.6,-0.45) and (7.2,-0.45) .. (7,-0.1);
\draw[-] (6.7,0) .. controls (6.5,0.35) and (6.1,0.35) .. (6.1,-0.05);
\draw[-latex] (6.1,-0.05) .. controls (6.1,-0.45) and (6.5,-0.45) .. (6.7,-0.1);
\node at (6.85,0 ){ \small$1$};
\node at (5.9,0 ){ \small$\alpha$};
\node at (7.8,0 ){ \small\color{purple}{$\varepsilon$}};
\end{tikzpicture}
\end{picture}
\end{equation}
\vspace{-2mm}

Then $\Lambda_{\texttt{J}}^\imath$ is a commutative $\bfk$-algebra. Clearly, $\rep^{\rm nil}(\Lambda_{\texttt{J}}^\imath)\cong \calc_1(\rep^{\rm nil}_\bfk (\QJ))$, and we shall identify these categories below. We can view $\rep^{\rm nil}_\bfk(\QJ)$ naturally as a full subcategory of $\rep^{\rm nil}(\Lambda_{\texttt{J}}^\imath)$.

It is well known that $\rep_\bfk^{\rm nil}(\QJ)$ is a uniserial category. Let $S$ be the simple object in $\rep_\bfk^{\rm nil}(\QJ)$.
Then any indecomposable object of $\rep_\bfk^{\rm nil}(\QJ)$  (up to isomorphisms) is of the form $S^{(n)}$ of length $n\geq1$.
Thus the set of isomorphism classes $\Iso(\rep_\bfk^{\rm nil}(\QJ))$ is canonically isomorphic to the set $\cp$ of all partitions, via the assignment
\begin{align}
\la =(\la_1,\la_2,\cdots,\la_r)\mapsto S^{(\la)}= S^{(\la_1)}\oplus \cdots \oplus S^{(\la_r)}.
\end{align}

Consider the Hall algebra $\ch\big(\calc_1(\rep^{\rm nil}_\bfk (\QJ))\big)$. Let $\cJ$ be the two-sided ideal of  $\ch\big(\calc_1(\rep^{\rm nil}_\bfk (\QJ))\big)$ generated by
\begin{align}
  \label{eq:idealJ}
&\Big\{ [M^\bullet]-[N^\bullet]\mid H(M^\bullet)\cong H(N^\bullet), \quad \widehat{\Im d_{M^\bullet}}=\widehat{\Im d_{N^\bullet}} \Big\}.
\end{align}
The quotient algebra
\begin{align}
\tMHLJ:=\ch\big(\calc_1(\rep^{\rm nil}_\bfk (\QJ))\big)/\cJ
\end{align}
is called the {\em$\imath$Hall algebra} of $\Lambda_{\texttt{J}}^\imath$; compare \cite{LW19, LW20} (see also Section~\ref{sec:Hall}). We use $\ast$ to denote the multiplication of $\tMHLJ$ (this is compatible with the twisted multiplication of Hall algebras since the Euler form of $\calc_1(\rep^{\rm nil}_\bfk (\QJ))$ is trivial). The version of the $\imath$Hall algebra defined here does not require $[K_S]$ to be invertible.

For any acyclic complex $K^\bullet=(K,d)$, by definition, $[K^\bullet]=[K_{\Im d}]$ in  $\tMHLJ$. In fact, if $\widehat{X}=\widehat{Y}$ in $K_0(\rep^{\rm nil}_\bfk (\QJ))$, we have $[K_X]=[K_Y]$. So we sometimes denote $[K_X]$ by $[K_{\widehat{X}}]$ in $\tMHLJ$ below. For any complex $M^\bullet=(M,d)$, by definition we have
\begin{align}
\label{eq:decom1}
[M^\bullet]=[H(M^\bullet)\oplus K_{\Im d}]=[H(M^\bullet)]\ast[K_{\widehat{\Im d}}].
\end{align}

\begin{lemma}
\label{lem:iHallbasis}
The identity morphism of $\ch(\cc_1(\rep^{\rm nil}_\bfk (\QJ) ))$ induces an algebra embedding
$\Upsilon: \tMHLJ\rightarrow \cs\cd\ch(\cc_1(\rep^{\rm nil}_\bfk (\QJ) ))$.
\end{lemma}

\begin{proof}
For any $M^\bullet,N^\bullet$, if $H(M^\bullet)\cong H(N^\bullet)$ and $\widehat{\Im d_{M^\bullet}}=\widehat{\Im d_{N^\bullet}}$, then $[M^\bullet]=[N^\bullet]$ in $\cs\cd\ch(\cc_1(\rep^{\rm nil}_\bfk (\QJ) ))$ by \eqref{eq:decom}. So the identity morphism of $\ch(\cc_1(\rep^{\rm nil}_\bfk (\QJ) ))$ induces an algebra morphism
$\Upsilon: \tMHLJ\rightarrow \cs\cd\ch(\cc_1(\rep^{\rm nil}_\bfk (\QJ) ))$.

The injectivity of $\Upsilon$ follows from \eqref{eq:decom1} and Proposition \ref{prop:hallbasis}.
\end{proof}

\begin{proposition}
\label{prop:iHallbasis}
The following hold in $\tMHLJ$:
\begin{enumerate}
\item
For any acyclic complex $K$, $[K]$ is a central element of  $\tMHLJ$.
\item
$\tMHLJ$ has a basis (called $\imath$Hall basis) given by
\begin{align*}
\{[M]*[K_S]^{a}\mid [M]\in\Iso(\rep^{\rm nil}_\bfk (\QJ)), a\in\N\}.
\end{align*}
\end{enumerate}
\end{proposition}

\begin{proof}
The first statement is trivial.
The second one follows  from the proof of Lemma \ref{lem:iHallbasis}, since the rank of $K_0(\rep^{\rm nil}_\bfk (\QJ))$ is one.
\end{proof}

The multiplication formula \eqref{Hallmult} holds in $\tMHLJ$.

\begin{lemma}
\label{lem:comm Jordan}
The algebra $\tMHLJ$ is commutative.
\end{lemma}

\begin{proof}
Since $\cs\cd\ch(\cc_1(\rep^{\rm nil}_\bfk (\QJ) ))$ is commutative by \cite[Lemma 4.5]{LRW20}, the result follows from Lemma \ref{lem:iHallbasis}.
\end{proof}

\begin{proposition}
\label{prop:iHJordan}
The $\imath$Hall algebra $\tMHLJ$ of the Jordan quiver is isomorphic to
\begin{enumerate}
\item
the polynomial $\Q$-algebra in the infinitely many generators
$[K_S], [S^{(1^r)}],$ for $r\ge 1$;
\item
the polynomial $\Q$-algebra in the infinitely many generators
$[K_S], [S^{(r)}],$ for $r\ge 1$.
\end{enumerate}
\end{proposition}

\begin{proof}

We only prove the first statement, and the second one can be proved in the same way.

Let $\ch(\bfk \QJ)$ be the (twisted) Hall algebra of $\rep_\bfk^{\rm nil}(\QJ)$. From \cite[Chapter II, (2.3)]{Mac95}, $\ch(\bfk \QJ)$ is isomorphic to the polynomial algebra in the infinitely many generators $[S^{(1^r)}],$ for $r\ge 1$.

By Lemma \ref{lem:comm Jordan}, the algebra $\tMHLJ$ is commutative.
By Proposition \ref{prop:iHallbasis}(2), $\tMHLJ$ has a Hall basis given by
$[M] * [K_S]^a$, where $[M]\in\Iso(\rep^{\rm nil}_\bfk (\QJ))$ and $a\in\N$. By \cite[Lemma~5.4]{LW19}, $\tMHLJ$ is a filtered algebra by setting $\deg ([S^{(\la)}]*[K_S]^a)=|\la|$, whose  associated graded algebra is isomorphic to the tensor algebra $\ch(\bfk \QJ)\otimes \Q\big[[K_S]\big]$.
Then the desired result follows by a standard filtered algebra argument.
\end{proof}

\subsection{Generic $\imath$Hall algebra}
\label{subsec:generic}

According to Hall and Steinitz (cf. \cite[III]{Mac95}), there exists a unique polynomial $G_{\mu,\nu}^\lambda(T)\in\Z[T]$ such that
$$G_{S^{(\mu)},S^{(\nu)}}^{S^{(\lambda)}}=G_{\mu,\nu}^\lambda(q).$$

For any partition $\lambda$, recall from \eqref{eq:nla} that $n(\la)=\sum_i(i-1)\la_i$ and $|\la|=\sum_{i\geq1}\la_i$. The following is well known; see, e.g., \cite[Lemma~ 2.8]{Sch06}.
\begin{align}
\label{eq:aut}
|\Aut(S^{(\lambda)})|=&\bq^{|\lambda|+2n(\lambda)}\prod_i(1-\bq^{-1})(1-\bq^{-2})\cdots (1-\bq^{-m_i(\lambda)})
\\\notag
=&\bq^{|\lambda|+2n(\lambda)}\prod_{i\geq1}\varphi_{m_i(\lambda)}(\bq^{-1}).
\end{align}

\begin{lemma}
\label{lem:Hallpoly}
For any partitions $\lambda,\mu,\nu$, there exists a Laurent polynomial $\mathcal G_{\mu,\nu}^{\lambda}(T)\in \Z[T,T^{-1}]$ such that
\begin{align}
[S^{(\mu)}]*[S^{(\nu)}]=&
\sum_{\lambda} \mathcal G_{\mu,\nu}^{\lambda}(\bq)[S^{(\lambda)}]*[K_{S}]^{\frac{|\mu|+|\nu|-|\lambda|}{2}}.
\end{align}
\end{lemma}

\begin{proof}
By Corollary \ref{cor:Hallmult} and the above, we can write
\begin{align}
[S^{(\mu)}]*[S^{(\nu)}]=&
\sum_{\lambda} \mathcal G_{\mu,\nu}^{\lambda}(\bq)[S^{(\lambda)}]*[K_{S}]^{\frac{|\mu|+|\nu|-|\lambda|}{2}},
\end{align}
where $\mathcal G_{\mu,\nu}^{\lambda}(T)=\frac{F_1(T)}{F_2(T)}$ for  $F_1(T),F_2(T)\in\Z[T]$, and $F_2(\bq)=|\aut(S^{(\la)})|$.
A comparison with \eqref{Hallmult} shows that $\bq^m \mathcal G_{\mu,\nu}^{\lambda}(q)=\frac{\bq^m F_1(q)}{F_2(q)}\in\Z$ for some $m\in \N$.  Therefore, $F_2(T)$ divides $T^m F_1(T)$, and thus we have $\mathcal G_{\mu,\nu}^{\lambda}(T)\in \Z[T,T^{-1}]$ by noting that $\bq$ can be any prime power.
\end{proof}


We now define the {\em generic $\imath$Hall algebra} of the Jordan quiver as the generic (twisted) modified Ringel-Hall algebra of $\LaJ^\imath$  over $\Q(t)$, and denote it by $\tMHg$. The algebra $\tMHg$ is the free $\Q(t)[K_\de]$-module with a basis $\{\fu_\lambda\mid \lambda\in\cp\}$ and multiplication given by
\begin{align}
\label{eq:g}
\fu_\mu*\fu_\nu=\sum_{\lambda\in\cp}\mathcal G_{\mu,\nu}^\lambda(t^{-1})\fu_\lambda*K_\de^{\frac{|\mu|+|\nu|-|\lambda|}{2}}.
\end{align}

\begin{remark}
 \label{rem:filter}
The filtered algebra structure on $\tMHg$ (or $\tMHLJ$) can be made more precise. The multiplication formula in the associated graded of $\tMHg$ is modified from \eqref{eq:g} by setting $K_\de=0$. 
In particular, when $|\la|=|\mu|+|\nu|$, $\mathcal G_{\mu,\nu}^\lambda(t)$ is equal to  $ g_{\mu,\nu}^\lambda(t)$ in \cite[III, (4.2)]{Mac95} multiplied by
$t^{2(n(\mu)+n(\nu)-n(\la))} \prod_{i\geq1} {\varphi_{m_i(\mu)}(t^{-1})\varphi_{m_i(\nu)}(t^{-1})}/{\varphi_{m_i(\lambda)}(t^{-1})}.$ (The difference here is due to the different normalizations of the Hall products.)
\end{remark}

\subsection{The down-up horizontal $\imath$Pieri rule} 

We now establish the first $\imath$Pieri rule in $\tMHLJ$.

\begin{lemma}
\label{lem:HallfomrulaJordan}
For $r\geq 1$ and any partition $\mu$, the following identities hold in $\tMHLJ$:
\begin{align}
\label{Pieri:h}
[S^{(\mu)}]*[S^{(r)}]
&= \sum_{a=0}^r\sum_{\la,\nu} G_{\nu,(r-a)}^\la(\bq)G_{(a),\nu}^\mu(\bq)
\\\notag
&\qquad\qquad \times
\frac{|\Aut(S^{(\nu)})|\cdot |\Aut(S^{(a)})|\cdot|\Aut(S^{(r-a)})| }{|\Aut(S^{(\la)})|}   [S^{(\la)}]*[K_S]^a,
 \\
\label{Pieri:v}
[S^{(\mu)}]*[S^{(1^r)}]
&= \sum_{a=0}^r\sum_{\la,\nu} G_{\nu,(1^{r-a})}^\la(\bq) G_{(1^a),\nu}^\mu(\bq) \frac{|\Aut(S^{(\nu)})| \cdot|\Aut(S^{(1^{r-a})})| }{|\Aut(S^{(\la)})|}
\\\notag
&\qquad\qquad \times
(\bq^r-1)\cdots(\bq^r-\bq^{a-1}) [S^{(\la)}]*[K_S]^a.
\end{align}
\end{lemma}

\begin{proof}
Since $S^{(r)}$ is uniserial, its sub-object has the form $S^{(a)}$ for some $0\leq a\leq r$, and then $G_{S^{(r-a)},S^{(a)}}^{S^{(r)}}=1$
(equivalently, $G_{{(r-a)},{(a)}}^{{(r)}}(\bq)=1$). Then \eqref{Pieri:h} follows immediately from Corollary~ \ref{cor:Hallmult}.

Similarly, since $S^{(1^r)}$ is semisimple, its sub-object has the form $S^{(1^a)}$ for some $0\leq a\leq r$, and then
$$G_{S^{(1^{r-a})},S^{(1^a)}}^{S^{(1^r)}}\cdot |\aut(S^{(a)})|=(\bq^r-1)\cdots(\bq^r-\bq^{a-1}).$$
Then \eqref{Pieri:v} follows from Corollary~ \ref{cor:Hallmult}.
\end{proof}

We shall now renormalize the formulas \eqref{Pieri:h}--\eqref{Pieri:v} to make them compatible with the symmetric function side.

For a horizontal $r$-strip $\sigma=\nu-\mu$, let $I=I_{\nu-\mu}$ be the set of integers $i\geq 1$ such that $\sigma_i'=1$ and $\sigma_{i+1}'=0$.
Recall from \cite[II,(4.13)]{Mac95} that
\begin{align*}
G_{\mu,(r)}^{\nu}(T)
=& \frac{T^{n(\nu)-n(\mu)}}{1-T^{-1}}\prod\limits_{i\in I_{\nu-\mu}}(1-T^{-m_i(\nu)})
\end{align*}
if $r>0$; and $G_{\mu,(0)}^{\nu}(T)=1$. So we have
\begin{align}
\label{eq:Gmurnu}
G_{\mu,(r)}^{\nu}(\bq)\cdot|\aut(S^{(r)})|=\bq^{n(\nu)-n(\mu)+r}\prod\limits_{i\in I_{\nu-\mu}}(1-\bq^{-m_i(\nu)}),\quad \forall r\geq0.
\end{align}

For any partition $\lambda$, denote by
\begin{align}
\label{eq:haV}
\iVh_\lambda:= \bq^{-|\lambda|-n(\lambda)}[S^{(\lambda)}]\in \iRH(\bfk\QJ).
\end{align}
In particular,
\begin{align}
\label{eq:haV-hv}
\iVh_{(r)}= \bq^{-r}[S^{(r)}], \qquad \quad \iVh_{(1^r)}= \bq^{-\binom{r+1}{2}}[S^{(1^r)}],\quad \forall r>0.
\end{align}

Introduce the following notations, for $r\ge 0$ and any partition $\la$:
\begin{align}
\label{eq:varphi:r}
\varphi_r(t) &= (1-t)(1-t^2)\cdots (1-t^r);
\\
\label{def:bla}
b_\la(t) &= \prod_{i\geq1} \varphi_{m_i(\la)}(t).
\end{align}

\begin{proposition}
\label{prop:PieriHall-Ho}
For $r\ge 1$ and any partition $\mu$, we have
\begin{align}
\label{eq:haPH-ho}
\iVh_\mu* \iVh_{(r)}= &\sum_{a+b =r}
\sum_{\nu \stackrel{a}{\rightarrow} \mu}\sum_{\nu \stackrel{b}{\rightarrow} \la
} \bq^{-a}\, \varphi_{\mu/\nu}(\bq^{-1}) \psi_{\la/\nu}(\bq^{-1}) \;  \iVh_\la*[K_S]^{a}.
\end{align}
\end{proposition}

\begin{proof}
By \eqref{eq:haV}--\eqref{eq:haV-hv}, the identity \eqref{Pieri:h} translates into
\begin{align}   \label{Pieri:hSF}
\begin{split}
\iVh_\mu  * \iVh_{(r)}
&= \sum_{a=0}^r\sum_{\la,\nu} \bq^{-(r+|\mu|-|\la|+n(\mu)-n(\la))} G_{\nu,(r-a)}^\la(\bq) G_{(a),\nu}^\mu(\bq)
 \\
&\qquad\qquad  \times \frac{|\Aut(S^{(\nu)})| \cdot|\Aut(S^{(a)})|\cdot|\Aut(S^{(r-a)})| }{|\Aut(S^{(\la)})|}   \iVh_\la*[K_S]^a.
\end{split}
 \end{align}
Plugging \eqref{eq:aut} and \eqref{eq:Gmurnu} into \eqref{Pieri:hSF}, we have
\begin{align*}
\iVh_\mu * \iVh_{(r)}
&= \sum_{a=0}^r\sum_{\la,\nu} \bq^{-(r+|\mu|-|\la|+n(\mu)-n(\la))} \bq^{-n(\nu)+n(\la)} \bq^{-n(\nu)+n(\mu)} \bq^{|\nu|+2n(\nu)} \bq^{r} \bq^{-|\la|-2n(\la)}
\\
&\qquad\qquad \times
\prod_{i\in I_{\la-\nu}} (1-\bq^{-m_i(\la)}) \prod_{i\in I_{\mu-\nu}} (1-\bq^{-m_i(\mu)}) \frac{\prod_i \varphi_{m_i(\nu)}(\bq^{-1})}{\prod_{i} \varphi_{m_i(\la)}(\bq^{-1})} \iVh_\la*[K_S]^a
\\
&=\sum_{a+b =r}
\sum_{\nu \stackrel{a}{\rightarrow} \mu}\sum_{\nu \stackrel{b}{\rightarrow} \la
}  \bq^{-a}\varphi_{\la/\nu}(\bq^{-1})\varphi_{\mu/\nu}(\bq^{-1}) \frac{b_{\nu}(\bq^{-1})}{b_\la(\bq^{-1})}\;\iVh_\la*[K_S]^a.
\end{align*}
The above formula can be converted to \eqref{eq:haPH-ho} using the following identity  \cite[III (5.12)]{Mac95}
\begin{align}
\label{eq:bphi}
\psi_{\la/\nu}(\bq^{-1})=\varphi_{\la/\nu}(\bq^{-1}) \frac{b_{\nu}(\bq^{-1})}{b_\la(\bq^{-1})}.
\end{align}
The proposition follows.
\end{proof}

\subsection{$\imath$Hall algebra and the ring $\La_{\vth}$}

Recall the  ring $\La_{t,\vth} = \Q(t)[\vth] [v_1, v_2, \ldots]$. Denote $\La_{\vth} :=\La_{t,\vth} \otimes_{\Q(t)} \Q$, where $t$ acts on $\Q$ by $\bq^{-1}$. Recall $\iVh_\la$ from \eqref{eq:haV}.

\begin{theorem}
\label{thm:iso}
There exists a $\Q$-algebra isomorphism $\Phi_{(\bq)}: \iRH(\bfk \QJ) \longrightarrow \La_{\vth}$ such that
\begin{align*}
\Phi_{(\bq)} ([K_S]) = q\vth,\qquad
\Phi_{(\bq)} ([S^{(r)}]) =q^r v_r  \quad (r\ge 1).
\end{align*}
Moreover, for any partition $\la$, we have
\begin{align}
\Phi_{(\bq)}([S^{(\la)}])=\bq^{|\la| +n(\la)} \iV_\la,
\qquad
\Phi_{(\bq)}(\iVh_\la)= \iV_\la.
\end{align}
\end{theorem}

\begin{proof}
The first statement follows from Proposition \ref{prop:iHJordan}.

The second statement will be proved using Proposition \ref{prop:PieriHall-Ho} and Theorem~ \ref{thm:Pieri-Ho}, by inductions on $n=|\lambda|$ and the dominance order on $\cp_n$.
The case for $n=1$ is clear.

The case for $\lambda$ with $\ell(\lambda)=1$ follows by definition. If $\ell(\lambda)>1$, let $\mu$ be the partition obtained from $\lambda$ by deleting the last row. Suppose that this last row has $r$ elements. Then by Proposition~ \ref{prop:PieriHall-Ho},
we have
\begin{align*}
\iVh_\mu* \iVh_{(r)}= &\iVh_\la+
\sum_{\mu \stackrel{r}{\rightarrow} \gamma,\la\rhd\gamma}   \psi_{\la/\mu}(\bq^{-1}) \;  \iVh_\gamma
+\sum_{a+b =r,a\neq0}
\sum_{\nu \stackrel{a}{\rightarrow} \mu}\sum_{\nu \stackrel{b}{\rightarrow} \gamma
} \bq^{-a}\, \varphi_{\mu/\nu}(\bq^{-1}) \psi_{\gamma/\nu}(\bq^{-1}) \;  \iVh_\gamma*[K_S]^{a}.
\end{align*}
By Theorem \ref{thm:Pieri-Ho}, we have
\begin{align*}
\iV_\mu  \cdot \vv_r
=& \iV_\la+
\sum_{\mu \stackrel{r}{\rightarrow} \gamma,\la\rhd\gamma
}  \psi_{\gamma/\mu}(t) \;  \iV_\gamma
+ \sum_{a+b =r,a\neq0}
\sum_{\nu \stackrel{a}{\rightarrow} \mu}\sum_{\nu \stackrel{b}{\rightarrow} \gamma
} \vth^{a}\, \varphi_{\mu/\nu}(t) \psi_{\gamma/\nu}(t) \;  \iV_\gamma.
\end{align*}
By the inductive assumption and comparing the above 2 identities, we have $\Phi_{(\bq)}(\iVh_\la)= \iV_\la$.
\end{proof}

\begin{remark}
 \label{rem:PQ}
The Hall products in this paper and in \cite{Mac95} use different normalizations (though the resulting Hall algebras of the Jordan quiver are isomorphic via $[S^\la] \mapsto [S^\la]/ |\Aut(S^\la)|$); cf. Remark~\ref{rem:filter}. Under our Hall product, the isomorphism in \cite[(3.4)]{Mac95} would become $u_\la \mapsto q^{|\la| +n(\la)} Q_\la(x;q^{-1})$.
\end{remark}

\subsection{The up-down horizontal $\imath$Pieri rule}

We shall convert the down-up horizontal $\imath$Pieri rule to a new up-down version. The two formulas are related by Green's formula \cite{Gr95}.

\begin{lemma}\label{recursion formula of phi and psi}
Assume $\cu_{a,b}$ and $\cv_{a,b}$ (for $a,b\in\N$) satisfy the relations
\begin{align}
  \label{eq:UV}
\cu_{a,b}=\cv_{a,b}+\sum\limits_{i=1}^{\min(a,b)}\bq^{i-1}(\bq-1)\cv_{a-i,b-i},
\quad \text{ for all } a,b.
\end{align}
Then we have
\[
\cv_{a,b}=\cu_{a,b}-(\bq-1)\sum\limits_{i=1}^{\min(a,b)}\cu_{a-i,b-i},
\quad \text{ for all } a,b.
\]
\end{lemma}

\begin{proof}
Set $\cv_{a,b}=\cu_{a,b} =0$, if either $a$ or $b$ is negative.
If $a=0$ or $b=0$, then we have $\cv_{a,b}=\cu_{a,b}$.
By induction on $\min(a,b)$ we have
\begin{align*}
\cv_{a,b}&= \cu_{a,b}-\sum\limits_{i\ge 1}\bq^{i-1}(\bq-1)\cv_{a-i,b-i}
\\
&=\cu_{a,b}-\sum\limits_{i\ge 1}\bq^{i-1}(\bq-1)\cu_{a-i,b-i}
+\sum\limits_{i, j \ge 1}\bq^{i-1}(\bq-1)^2 \cu_{a-i-j,b-i-j} \quad (\text{set }k=i+j)
\\
&=\cu_{a,b}-\sum\limits_{i\ge 1}\bq^{i-1}(\bq-1)\cu_{a-i,b-i}
+\sum\limits_{k\ge 2} (\bq^{k-1}-1)(\bq-1)\cu_{a-k,b-k},
\end{align*}
which can be readily converted to the formula in the lemma.
\end{proof}

\begin{lemma}
For any $a,b\geq 0$, we have
\begin{align}
\label{inverse of Green's formula2}
\sum_{\nu}&G_{S^{(\nu)},S^{(b)}}^{S^{(\lambda)}}
G_{S^{(\nu)}, S^{(a)}}^{S^{(\mu)}}\frac{|\Aut(S^{(\nu)})| \cdot |\aut(S^{(a)})| \cdot |\aut(S^{(b)})|}{|\Aut(S^{(\lambda)})|}
\\\notag
=&\sum_{\xi} G_{S^{(\lambda)},S^{(a)}}^{S^{(\xi)}}  G_{S^{(\mu)}, S^{(b)}}^{S^{(\xi)}} \frac{|\Aut(S^{(\mu)})| \cdot |\aut(S^{(a)})| \cdot |\aut(S^{(b)})|}{|\Aut(S^{(\xi)})|}
\\\notag
& -\sum\limits_{i=1}^{\min(a,b)}(\bq-1) \sum_{\xi} G_{S^{(\lambda)},S^{(a-i)}}^{S^{(\xi)}}  G_{S^{(\mu)}, S^{(b-i)}}^{S^{(\xi)}} \frac{|\Aut(S^{(\mu)})| \cdot |\aut(S^{(a-i)})| \cdot |\aut(S^{(b-i)})|}{|\Aut(S^{(\xi)})|}.
\end{align}
\end{lemma}

\begin{proof}

Observe that the Euler form is trivial for the Jordan quiver. By Green's formula \cite{Gr95},
\begin{align}\label{Green's formula}
\sum_{\xi} G_{S^{(\lambda)},S^{(a)}}^{S^{(\xi)}}  &G_{S^{(\mu)}, S^{(b)}}^{S^{(\xi)}} \frac{1}{|\Aut(S^{(\xi)})|}
=\sum\limits_{i=0}^{\min(a,b)}\sum_{\nu}G_{S^{(\nu)},S^{(b-i)}}^{S^{(\lambda)}}
G_{S^{(\nu)}, S^{(a-i)}}^{S^{(\mu)}}
\\\notag
&\times \frac{|\Aut(S^{(i)})|\cdot |\Aut(S^{(a-i)})|\cdot |\Aut(S^{(b-i)})|\cdot |\Aut(S^{(\nu)})|}{|\Aut(S^{(a)})|\cdot |\Aut(S^{(b)})|\cdot |\Aut(S^{(\lambda)})|\cdot |\Aut(S^{(\mu)})|}.
\end{align}
%
Denote by
\[
{\cu}_{a,b}=\sum_{\xi} G_{S^{(\lambda)},S^{(a)}}^{S^{(\xi)}}  G_{S^{(\mu)}, S^{(b)}}^{S^{(\xi)}} \frac{|\Aut(S^{(\mu)})| \cdot |\aut(S^{(a)})| \cdot |\aut(S^{(b)})|}{|\Aut(S^{(\xi)})|}
\]
and
\[
{\cv}_{a,b}=\sum_{\nu}G_{S^{(\nu)},S^{(b)}}^{S^{(\lambda)}}
G_{S^{(\nu)}, S^{(a)}}^{S^{(\mu)}}\frac{|\Aut(S^{(\nu)})| \cdot |\aut(S^{(a)})| \cdot |\aut(S^{(b)})|}{|\Aut(S^{(\lambda)})|}.
\]
Recall that $|\aut(S^{(i)})|=\bq^{i-1}(\bq-1)$ for $i \geq 1$.
The identity \eqref{Green's formula} can be reformulated as
the identity \eqref{eq:UV}.
Then by Lemma \ref{recursion formula of phi and psi} we have ${\cv}_{a,b}={\cu}_{a,b}-(\bq-1)\sum\limits_{i=1}^{\min(a,b)}{\cu}_{a-i,b-i}$ for any $a,b\geq 0$, and whence \eqref{inverse of Green's formula2}.
\end{proof}


\begin{proposition}
\label{prop:PieriHall-Ho2}
For $r\ge 1$ and any partition $\mu$, we have
\begin{align}
\iVh_\mu* \iVh_{(r)}=&\sum_{a+b =r}
\sum_{\la \stackrel{a}{\rightarrow} \xi}\sum_{\mu \stackrel{b}{\rightarrow} \xi} \bq^{-a}\, \varphi_{\xi/\la}(\bq^{-1}) \psi_{\xi/\mu}(\bq^{-1}) \;  \iVh_\la*[K_S]^{a}
\\\notag
&+\sum_{a+b =r}\sum_{i=1}^{\min(a,b)}
\sum_{\la \stackrel{a-i}{\rightarrow} \xi}\sum_{\mu \stackrel{b-i}{\rightarrow} \xi}\bq^{-a}(\bq^{-i}-\bq^{-i+1})  \varphi_{\xi/\la}(\bq^{-1}) \psi_{\xi/\mu}(\bq^{-1}) \;  \iVh_\la*[K_S]^{a}.
\end{align}
\end{proposition}

\begin{proof}
Combining with \eqref{eq:aut}, \eqref{eq:Gmurnu} and \eqref{eq:bphi}, we obtain that, for any $0\leq i\leq \min(a,b)$, $\la \stackrel{a-i}{\rightarrow} \xi$ and  $\mu \stackrel{b-i}{\rightarrow} \xi$,
\begin{align}
 \label{the coefficient of multiple two Hall numbers}
&G_{S^{(\lambda)},S^{(a-i)}}^{S^{(\xi)}}  G_{S^{(\mu)}, S^{(b-i)}}^{S^{(\xi)}}
 \frac{|\Aut(S^{(\mu)})| \cdot  |\Aut(S^{(a-i)})| \cdot |\Aut(S^{(b-i)})| }{|\Aut(S^{(\xi)})|}\\\notag
 =&\bq^{n(\xi)-n(\la)+a-i}
\prod\limits_{j\in I_{\xi-\la}}(1-\bq^{-m_j(\xi)})\cdot\bq^{n(\xi)-n(\mu)+b-i}\prod\limits_{j\in I_{\xi-\mu}}(1-\bq^{-m_j(\xi)})  \\\notag
&\qquad\qquad \times \bq^{|\mu|+2n(\mu)}\prod_{j\geq1}
\varphi_{m_j(\mu)}(\bq^{-1})\cdot \Big(\bq^{|\xi|+2n(\xi)}\prod_{j\geq1}\varphi_{m_j(\xi)}(\bq^{-1})\Big)^{-1}\\\notag
=&\bq^{n(\mu)-n(\la)+a-i} \varphi_{\xi/\la}(\bq^{-1}) \varphi_{\xi/\mu}(\bq^{-1}) \frac{b_{\mu}(\bq^{-1})}{b_{\xi}(\bq^{-1})}\\\notag
=&\bq^{n(\mu)-n(\la)+a-i} \varphi_{\xi/\la}(\bq^{-1}) \psi_{\xi/\mu}(\bq^{-1}).
\end{align}

Plugging \eqref{inverse of Green's formula2} into \eqref{Pieri:h} and using \eqref{the coefficient of multiple two Hall numbers}, we have
\begin{align*}
[S^{(\mu)}]*[S^{(r)}]
&=\sum_{a+b=r}\sum_{\la \stackrel{a}{\rightarrow} \xi}\sum_{\mu \stackrel{b}{\rightarrow} \xi} \bq^{n(\mu)-n(\la)+a} \varphi_{\xi/\la}(\bq^{-1}) \psi_{\xi/\mu}(\bq^{-1})
[S^{(\la)}]*[K_S]^a
\\\notag
& +(1-\bq)\sum_{a+b=r}\sum\limits_{i=1}^{\min(a,b)} \sum_{\la \stackrel{a-i}{\rightarrow} \xi}\sum_{\mu \stackrel{b-i}{\rightarrow} \xi} \bq^{n(\mu)-n(\la)+a-i} \varphi_{\xi/\la}(\bq^{-1}) \psi_{\xi/\mu}(\bq^{-1})
[S^{(\la)}]*[K_S]^a.
\end{align*}
This identity can then easily be converted to the identity in the proposition by using
$\iVh_\lambda= \bq^{-|\lambda|-n(\lambda)}[S^{(\lambda)}]$ for any partition $\lambda$, cf. \eqref{eq:haV}.
\end{proof}

\subsection{The down-up vertical $\imath$Pieri rule}

Let $n, r \in \Z$. Define
\begin{align*}\begin{bmatrix} n \\ r \end{bmatrix}_+(t)=\frac{\varphi_n(t)}{\varphi_r(t)\varphi_{n-r}(t)}
=\frac{(1-t^n)(1-t^{n-1})\cdots(1-t^{n-r+1})}{(1-t)(1-t^2)\cdots(1-t^r)}
\end{align*}
for $r\ge 0$ and
$\begin{bmatrix} n \\ r \end{bmatrix}_+(t)=0$
for $r<0$.

For any partitions $\la,\mu$ and $m\geq0$, define (cf. \cite[III, (3.2)]{Mac95})
\begin{align}
 \label{eq:fmu}
f_{\mu,(1^m)}^\la(t)=&\begin{cases}
\prod_{i\geq1} \qbinom{\la'_i-\la'_{i+1}}{\la'_i-\mu'_i}_+(t),& \text{ if }\lambda-\mu \text{ is a vertical }m\text{-strip},
\\
0,& \text{ otherwise.}
\end{cases}
\end{align}

Recall $n(\la)$ from \eqref{eq:nla}. By \cite[III, (3.3)]{Mac95}, we have
\begin{align}
\label{eq:Gf}
G^\la_{\mu,(1^m)}(q)=q^{n(\la)-n(\mu)-n(1^m)} f^\la_{\mu,(1^m)}(q^{-1}).
\end{align}
Given partitions $\mu, \nu$ and $a \in \N$, we use
\[
\mu {\downarrow}\nu,
\qquad (\text{respectively, } \mu \stackrel{a}{\downarrow}\nu)
\]
to denote that $\nu \leq \mu$ and $\mu -\nu$ is a vertical strip (respectively, a vertical $a$-strip).

\begin{proposition}
\label{prop:haPieri-v}
For $r\geq1$ and any partition $\mu$, we have
\begin{align}
\label{eq:haPieri-v1}
\iVh_\mu*\iVh_{(1^r)}
=&\sum_{a+b=r} \sum_{\la\stackrel{b}{\downarrow} \nu} \sum_{\mu \stackrel{a}{\downarrow}\nu}
\bq^{-a} \frac{b_{\nu}(\bq^{-1})}{b_\la(\bq^{-1})}\varphi_r(\bq^{-1})  f_{\nu,(1^{b})}^\la(\bq^{-1}) f_{\nu,(1^a)}^\mu(\bq^{-1}) \cdot \iVh_\la*[K_S]^a.
\end{align}
\end{proposition}

\begin{proof}
In this proof, we shall use the shorthand notation $\sum_{\la,\nu}$ to stand for $\sum_{\la\stackrel{b}{\downarrow} \nu} \sum_{\mu \stackrel{a}{\downarrow}\nu}$.
Starting from \eqref{eq:haV} and \eqref{Pieri:v} and using also \eqref{eq:Gf} and \eqref{eq:aut}, we have
\begin{align*}
\iVh_\mu*\iVh_{(1^r)}
&= \bq^{-|\mu|-n(\mu)-r-n(1^r)}[S^{(\mu)}] *[S^{(1^r)}]\\\notag
&= \bq^{-|\mu|-n(\mu)-r-n(1^r)}\sum_{a+b=r}\sum_{\la,\nu} G_{\nu,(1^{b})}^\la G_{(1^a),\nu}^\mu \frac{|\Aut(S^{(\nu)})| \cdot|\Aut(S^{(1^{b})})| }{|\Aut(S^{(\la)})|}
\\\notag
&\qquad\qquad\qquad\qquad\qquad
 \times (\bq^r-1)\cdots(\bq^r-\bq^{a-1})\cdot[S^{(\la)}]*[K_S]^a
\\\notag
&= \bq^{-|\mu|-n(\mu)-r-n(1^r)}\sum_{a+b=r}\sum_{\la,\nu} \bq^{n(\la)-n(\nu)-n(1^{b})}f_{\nu,(1^{b})}^\la(\bq^{-1})\cdot \bq^{n(\mu)-n(\nu)-n(1^a)}f_{\nu,(1^a)}^\mu (\bq^{-1})
\\\notag
&\qquad \times \frac{\bq^{|\nu|+2n(\nu)} b_\nu(\bq^{-1})\cdot \bq^{b^2}\varphi_{b}(\bq^{-1}) }{\bq^{|\la|+2n(\la)} b_\la(\bq^{-1})} (\bq^r-1)\cdots(\bq^r-\bq^{a-1}) \cdot \bq^{|\la|+n(\la)}\cdot\iVh_\la*[K_S]^a.
\end{align*}
The RHS above can be further simplified to be
\begin{align*}
=&\sum_{a+b=r}\sum_{\la,\nu} \bq^{-|\mu|-n(\mu)-r-n(1^r)} \bq^{n(\la)-n(\nu)-n(1^{b})} \bq^{n(\mu)-n(\nu)-n(1^a)}\bq^{|\la|+n(\la)} \frac{\bq^{|\nu|+2n(\nu)} \bq^{b^2} }{\bq^{|\la|+2n(\la)} }
\\\notag
&\times \frac{b_{\nu}(\bq^{-1})}{b_\la(\bq^{-1})} \varphi_{b}(\bq^{-1})\cdot \bq^{ra} (1-\bq^{a-r-1}) \cdots (1-\bq^{-r}) \cdot f_{\nu,(1^{b})}^\la(\bq^{-1}) f_{\nu,(1^a)}^\mu (\bq^{-1}) \cdot\iVh_\la*[K_S]^a
\\\notag
=&\sum_{a+b=r}\sum_{\la,\nu} \bq^{-|\mu|+|\nu|-r}\bq^{-n(1^r)-n(1^{b})-n(1^a)} \bq^{ra+b^2}
\frac{b_{\nu}(\bq^{-1})}{b_\la(\bq^{-1})} \varphi_r(\bq^{-1})
\\
&\qquad\qquad \times
f_{\nu,(1^{b})}^\la(\bq^{-1}) f_{\nu,(1^a)}^\mu (\bq^{-1}) \cdot \iVh_\la*[K_S]^a
\\\notag
=&\sum_{a+b=r}\sum_{\la,\nu}  \bq^{-a}\frac{b_{\nu}(\bq^{-1})}{b_\la(\bq^{-1})}\varphi_r(\bq^{-1})   f_{\nu,(1^{b})}^\la(\bq^{-1}) f_{\nu,(1^a)}^\mu (\bq^{-1}) \cdot \iVh_\la*[K_S]^a.
\end{align*}
The proof is completed.
\end{proof}

\begin{remark}
The coefficients of the RHS \eqref{eq:haPieri-v1} are in $\Z[q^{-1}]$, or equivalently, the coefficients in the vertical $\imath$Pieri rules (see Theorem {\bf C}) are in $\Z[t]$.
In fact, it follows first from Lemma \ref{lem:Hallpoly} and Theorem~ \ref{thm:iso} that the coefficients in the vertical $\imath$Pieri rules are in $\Z[t,t^{-1}]$, i.e.,
\begin{align*}
\sum_{a+b=r}\sum_{\la\stackrel{b}{\downarrow} \nu}  \sum_{\mu \stackrel{a}{\downarrow} \nu}\frac{b_{\nu}(t)}{b_{\la}(t)} \varphi_r(t) f_{\nu,(1^{b})}^\la(t) f_{\nu,(1^a)}^\mu(t)\in\Z[t,t^{-1}],
\end{align*}
for any $\mu,r,\la$.
They are actually in $\Z[t]$ by noting that $b_\nu(t),b_{\la}(t),\varphi_r(t),f_{\nu,(1^a)}^\mu(t),f_{\nu,(1^a)}^\mu(t)\in\Z[t]$, and $b_\la(0)=1$.
\end{remark}

\subsection{The up-down vertical $\imath$Pieri rule}

\begin{lemma}\label{recursion formula of phi and psi2}
Assume $\cu_{a,b}$ and $\cv_{a,b}$ (for $a,b\in\N$) satisfy the relations
\begin{align}
\label{eq:assump}
\cu_{a,b}=\sum\limits_{i=0}^{\min(a,b)}q^{-(a+b-i)i}\frac{1}{\varphi_i(q^{-1})}\cv_{a-i,b-i},
\quad \text{ for all } a,b.
\end{align}
Then we have
\[
\cv_{a,b}=\sum\limits_{i=0}^{\min(a,b)}(-1)^{i}\frac{q^{-i(a+b)+\binom{i+1}{2}}}{\varphi_i(q^{-1})}  \cu_{a-i,b-i}
\quad \text{ for all } a,b.
\]
\end{lemma}

\begin{proof}
Set $\cv_{a,b}=\cu_{a,b} =0$, if either $a$ or $b$ is negative.
If $a=0$ or $b=0$, then $\cv_{a,b}=\cu_{a,b}$. By induction on $\min(a,b)$, we have
\begin{align*}
\cv_{a,b}=&\cu_{a,b}-\sum\limits_{i\geq1}q^{-(a+b-i)i}\frac{1}{\varphi_i(q^{-1})}\cv_{a-i,b-i}
\\
=&\cu_{a,b}-\sum\limits_{i\geq1}\sum\limits_{j\geq0}q^{-(a+b-i)i}\frac{1}{\varphi_i(q^{-1})}(-1)^{j}\frac{q^{-j(a+b-2i)+\binom{j+1}{2}}}{\varphi_j(q^{-1})}  \cu_{a-i-j,b-i-j} \quad (\text{set }k=i+j)
\\
=&\cu_{a,b}-\sum_{k\geq1} (-1)^k q^{-k(a+b)+\binom{k+1}{2}} \Big(\sum_{i=1}^k(-1)^iq^{-\binom{i+1}{2}+ik} \frac{1}{\varphi_i(q^{-1})\varphi_{k-i}(q^{-1})} \Big)\cu_{a-k,b-k}
\\
=&\sum\limits_{k=0}^{\min(a,b)}(-1)^{k}\frac{q^{-k(a+b)+\binom{k+1}{2}}}{\varphi_k(q^{-1})}  \cu_{a-k,b-k}.
\end{align*}
The last equality used the following formula (cf. \cite[1.3.1(c)]{Lus93})
\begin{align*}
&\sum\limits_{i=1}^{k}(-1)^{i}\frac{q^{-\binom{i+1}{2}+ik}}{\varphi_i(q^{-1})\varphi_{k-i}(q^{-1})} =\frac{1}{\varphi_k(q^{-1})}\sum_{i=1}^{k}(-1)^{i}q^{-\binom{i+1}{2}+ik} \qbinom{k}{i}_+(q^{-1})=-\frac{1}{\varphi_k(q^{-1})}.
\end{align*}
\end{proof}

\begin{lemma}
\label{lem:Hallnum2}
For any $a,b\geq0$ and partitions $\la,\mu$, we have
\begin{align*}
\sum_{\nu}& G_{S^{(\nu)},S^{(1^{b})}}^{S^{(\lambda)}}
G_{S^{(\nu)}, S^{(1^{a})}}^{S^{(\mu)}}\frac{ |\Aut(S^{(\nu)})|}{|\Aut(S^{(\lambda)})|}
\\
&=\sum\limits_{i=0}^{\min(a,b)}\sum_{\xi} (-1)^i \frac{\bq^{-i(a+b)+\binom{i+1}{2}}}{\varphi_i(\bq^{-1})} G_{S^{(\lambda)},S^{(1^{a-i})}}^{S^{(\xi)}}  G_{S^{(\mu)}, S^{(1^{b-i})}}^{S^{(\xi)}}\frac{|\Aut(S^{(\mu)})|}{|\Aut(S^{(\xi)})|}.
\end{align*}
\end{lemma}

\begin{proof}
By Green's formula \cite{Gr95}, we have
\begin{align}
\sum_{\xi} & G_{S^{(\lambda)},S^{(1^a)}}^{S^{(\xi)}}  G_{S^{(\mu)}, S^{(1^b)}}^{S^{(\xi)}} \frac{1}{|\Aut(S^{(\xi)})|}\\\notag
& =\sum\limits_{i=0}^{\min(a,b)}\sum_{\nu}G_{S^{(\nu)},S^{(1^{b-i})}}^{S^{(\lambda)}}
G_{S^{(\nu)}, S^{(1^{a-i})}}^{S^{(\mu)}}G_{S^{(1^{i})}, S^{(1^{a-i})}}^{S^{(1^{a})}} G_{S^{(1^{i})}, S^{(1^{b-i})}}^{S^{(1^{b})}}
\\\notag
&\qquad\qquad \times \frac{|\Aut(S^{(1^i)})|\cdot |\Aut(S^{(1^{a-i})})|\cdot |\Aut(S^{(1^{b-i})})|\cdot |\Aut(S^{(\nu)})|}{|\Aut(S^{(1^a)})|\cdot |\Aut(S^{(1^b)})|\cdot |\Aut(S^{(\lambda)})|\cdot |\Aut(S^{(\mu)})|}.
\end{align}
For any $a,b,i\geq0$ with $i\leq \min(a,b)$, we have
\begin{align*}
G_{S^{(1^{i})}, S^{(1^{a-i})}}^{S^{(1^{a})}} G_{S^{(1^{i})}, S^{(1^{b-i})}}^{S^{(1^{b})}}
\frac{|\Aut(S^{(1^i)})|\cdot |\Aut(S^{(1^{a-i})})|\cdot |\Aut(S^{(1^{b-i})})|}{|\Aut(S^{(1^a)})|\cdot |\Aut(S^{(1^b)})|}
=&\bq^{-(a+b-i)i}\frac{1}{\varphi_i(\bq^{-1})}.
\end{align*}
So we have
\begin{align}
\label{eq:recursive}
\sum_{\xi} & G_{S^{(\lambda)},S^{(1^a)}}^{S^{(\xi)}}  G_{S^{(\mu)}, S^{(1^b)}}^{S^{(\xi)}}\frac{|\Aut(S^{(\mu)})|}{|\Aut(S^{(\xi)})|}\\\notag
&= \sum\limits_{i=0}^{\min(a,b)}\sum_{\nu}G_{S^{(\nu)},S^{(1^{b-i})}}^{S^{(\lambda)}}
G_{S^{(\nu)}, S^{(1^{a-i})}}^{S^{(\mu)}}\frac{ |\Aut(S^{(\nu)})|}{|\Aut(S^{(\lambda)})|}\cdot\bq^{-(a+b-i)i}\frac{1}{\varphi_i(\bq^{-1})}.
\end{align}

Denote by
\begin{align*}
\cu_{a,b}=&\sum_{\xi} G_{S^{(\lambda)},S^{(1^a)}}^{S^{(\xi)}}  G_{S^{(\mu)}, S^{(1^b)}}^{S^{(\xi)}}\frac{|\Aut(S^{(\mu)})|}{|\Aut(S^{(\xi)})|},
\quad
\cv_{a,b}= \sum_{\nu}G_{S^{(\nu)},S^{(1^{b})}}^{S^{(\lambda)}}
G_{S^{(\nu)}, S^{(1^{a})}}^{S^{(\mu)}}\frac{ |\Aut(S^{(\nu)})|}{|\Aut(S^{(\lambda)})|}.
\end{align*}
Then the identity \eqref{eq:recursive} can be converted to the identity \eqref{eq:assump}, 
and the result follows from Lemma \ref{recursion formula of phi and psi2}.
\end{proof}

\begin{proposition}
\label{prop:haPieri-vv}
For $r\geq1$ and any partition $\mu$, we have
\begin{align}
\label{eq:haPieri-vv1}
\iVh_\mu*\iVh_{(1^r)}
=&\sum_{a+b=r} \sum_{i=0}^{\min(a,b)} \sum_{\xi \stackrel{a-i}{\downarrow}\la}  \sum_{\xi\stackrel{b-i}{\downarrow}\mu} (-1)^i\bq^{-a-\frac{i(i-1)}{2}} \frac{b_{\mu}(\bq^{-1})}{b_{\xi}(\bq^{-1})} \frac{\varphi_r(\bq^{-1})}{\varphi_i(\bq^{-1})}
\\\notag
&\qquad\qquad\qquad \times  f_{\la,(1^{a-i})}^\xi(\bq^{-1}) f_{(1^{b-i}),\mu}^\xi (\bq^{-1}) \cdot \iVh_\la*[K_S]^{a}.
\end{align}
\end{proposition}

\begin{proof}
From \eqref{eq:haV}, 
\eqref{Pieri:v}  and  \eqref{eq:aut}, we have
\begin{align*}
\iVh_\mu*\iVh_{(1^r)}
=&\bq^{-|\mu|-n(\mu)-r-n(1^r)}[S^{(\mu)}] *[S^{(1^r)}]\\\notag
=&\bq^{-|\mu|-n(\mu)-r-n(1^r)}\sum_{a+b=r}\sum_{\la,\nu} G_{\nu,(1^{b})}^\la(\bq) G_{(1^a),\nu}^\mu(\bq) \frac{|\Aut(S^{(\nu)})|\cdot |\Aut(S^{(1^{b})})| }{|\Aut(S^{(\la)})|}
\\\notag
&\qquad\qquad\qquad\qquad \times
(\bq^r-1)\cdots(\bq^r-\bq^{a-1}) \cdot [S^{(\la)}]*[K_S]^a
\\
=&\bq^{-|\mu|-n(\mu)-r-n(1^r)}\sum_{a+b=r}\sum_{\la,\nu} \bq^{-ar+r^2+a^2}\varphi_r(\bq^{-1})
 \\
 &\qquad\qquad\qquad \times
 G_{\nu,(1^{b})}^\la(\bq) G_{(1^a),\nu}^\mu(\bq) \frac{|\Aut(S^{(\nu)})|  }{|\Aut(S^{(\la)})|} \cdot [S^{(\la)}]*[K_S]^a.
\end{align*}
By Lemma \ref{lem:Hallnum2}, we have
\begin{align*}
\iVh_\mu*\iVh_{(1^r)}
=&\bq^{-|\mu|-n(\mu)-r-n(1^r)}\sum_{a+b=r} \sum\limits_{i=0}^{\min(a,b)}\sum_{\la,\xi} \bq^{-ar+a^2+r^2}\varphi_r(\bq^{-1}) (-1)^i \frac{\bq^{-ir+\binom{i+1}{2}}}{\varphi_i(\bq^{-1})}
\\
&\qquad\qquad \times
G_{S^{(\lambda)},S^{(1^{a-i})}}^{S^{(\xi)}}  G_{S^{(\mu)}, S^{(1^{b-i})}}^{S^{(\xi)}} \frac{|\aut(S^{(\mu)})|}{|\Aut(S^{(\xi)})|}
 \cdot [S^{(\la)}]*[K_S]^a.
\end{align*}
Using \eqref{eq:Gf} and \eqref{eq:aut}, we then obtain
\begin{align*}
&\iVh_\mu*\iVh_{(1^r)}\\
=&\bq^{-|\mu|-n(\mu)-r-n(1^r)}\sum_{a+b=r}\sum\limits_{i=0}^{\min(a,b)}\sum_{\xi \stackrel{a-i}{\downarrow}\la}  \sum_{\xi\stackrel{b-i}{\downarrow}\mu} \bq^{-ar+a^2+r^2}\varphi_r(\bq^{-1}) (-1)^i \frac{\bq^{-ir+\binom{i+1}{2}}}{\varphi_i(\bq^{-1})}
\\
&\qquad\qquad\qquad \times
\bq^{n(\xi)-n(\la)-n(1^{a-i})}f_{\la,(1^{a-i})}^\xi(\bq^{-1})
\bq^{n(\xi)-n(\mu)-n(1^{b-i})}f_{(1^{b-i}),\mu}^\xi (\bq^{-1})
\\
&\qquad\qquad\qquad \times
\frac{\bq^{|\mu|+2n(\mu)} b_\mu(\bq^{-1}) }{\bq^{|\xi|+2n(\xi)} b_\xi(\bq^{-1})}   \bq^{|\la|+n(\la)}\cdot\iVh_\la*[K_S]^{a}
\\\notag
=&\sum_{a+b=r}\sum_{i=0}^{\min(a,b)}\sum_{\xi \stackrel{a-i}{\downarrow}\la}  \sum_{\xi\stackrel{b-i}{\downarrow}\mu} \bq^{-|\mu|-n(\mu)-r-n(1^r)} \bq^{n(\xi)-n(\la)-n(1^{a-i})} \bq^{n(\xi)-n(\mu)-n(1^{b-i})}\bq^{|\la|+n(\la)} \frac{\bq^{|\mu|+2n(\mu)}  }{\bq^{|\xi|+2n(\xi)} }
\\\notag
&\qquad \times \bq^{-ar+a^2+r^2}\varphi_r(\bq^{-1})(-1)^i \frac{\bq^{-ir+\binom{i+1}{2}}}{\varphi_i(\bq^{-1})} \frac{b_{\mu}(\bq^{-1})}{b_{\xi}(\bq^{-1})}  f_{\la,(1^{a-i})}^\xi(\bq^{-1})  f_{(1^{b-i}),\mu}^\xi (\bq^{-1}) \cdot\iVh_\la*[K_S]^{a}
\\\notag
=&\sum_{a+b=r} \sum_{i=0}^{\min(a,b)}\sum_{\xi \stackrel{a-i}{\downarrow}\la}  \sum_{\xi\stackrel{b-i}{\downarrow}\mu}(-1)^i\bq^{-a+i-r-\binom{r}{2} -\binom{a-i}{2}-\binom{b-i}{2}-ar+a^2+r^2-ir+\binom{i+1}{2}}\frac{\varphi_r(\bq^{-1})}{\varphi_i(\bq^{-1})}
\\\notag
& \qquad\qquad\qquad \times  \frac{b_{\mu}(\bq^{-1})}{b_{\xi}(\bq^{-1})}  f_{\la,(1^{a-i})}^\xi(\bq^{-1})  f_{(1^{b-i}),\mu}^\xi (\bq^{-1}) \cdot\iVh_\la*[K_S]^{a}
\\\notag
=&\sum_{a+b=r}\sum_{i=0}^{\min(a,b)}\sum_{\xi \stackrel{a-i}{\downarrow}\la}  \sum_{\xi\stackrel{b-i}{\downarrow}\mu}  (-1)^i\bq^{-a-\frac{i(i-1)}{2}} \frac{b_{\mu}(\bq^{-1})}{b_{\xi}(\bq^{-1})} \frac{\varphi_r(\bq^{-1})}{\varphi_i(\bq^{-1})}   f_{\la,(1^{a-i})}^\xi(\bq^{-1}) f_{(1^{b-i}),\mu}^\xi (\bq^{-1}) \cdot \iVh_\la*[K_S]^{a}.
\end{align*}
The proof is completed.
\end{proof}

\section{Identities arising from $\imath$Hall computations}
  \label{sec:comb}

This section plays a supportive role. We establish several identities which will be used in the $\imath$Hall algebra computations in Section~\ref{sec:identities}.

\subsection{Counting extensions}

\begin{lemma}
\label{extensions}
For any $r\geq1$ and any partition $\nu=(\nu_1, \nu_2, \ldots)$, we have
\begin{align}
\label{eq:Ext1}
\sum_{\mu:\, \ell(\mu)=\ell(\nu)+1}  \big|\Ext^1(S^{(\nu)},S^{(r)})_{S^{(\mu)}}\big|
=&q^{ \sum_{i}\min(r,\nu_i)-\ell(\nu)},
\\
\label{eq:Ext2}
\sum_{\mu:\, \ell(\mu)=\ell(\nu)} \big|\Ext^1(S^{(\nu)},S^{(r)})_{S^{(\mu)}}\big|
 =&q^{ \sum_{i} \min(r,\nu_i)}(1-q^{- \ell(\nu)}).
\end{align}
\end{lemma}

\begin{proof}
If $|\Ext^1(S^{(\nu)},S^{(r)})_{S^{(\mu)}}\big|\neq0$, then $\ell(\mu)=\ell(\nu)$ or $\ell(\mu)=\ell(\nu)+1$.
Now the formula \eqref{eq:Ext2} follows from \eqref{eq:Ext1}, since $\big|\Ext^1(S^{(\nu)},S^{(r)}) \big|
=\big|\Hom(S^{(\nu)},S^{(r)})|=q^{ \sum_{i} \min(r,\nu_i)}$.

We shall prove \eqref{eq:Ext1} by induction on $r$.

Let $r=1$. Then the only $[\xi]\in \Ext^1(S^{(\nu)}, S)_{S^{(\mu)}}$ such that $\ell(\mu)=\ell(\nu)+1$ is $[\xi]=0$.  Hence $\sum_{\mu:\, \ell(\mu)=\ell(\nu)+1}\big|\Ext^1(S^{(\nu)},S)_{S^{(\mu)}}\big|=1.$

For $r>1$, fix a short exact sequence
$0\rightarrow S\xrightarrow{f_1} S^{(r)} \xrightarrow{f_2} S^{(r-1)}\rightarrow0$.
By applying $\Hom(S^{(\nu)},-)$, we obtain the following long exact sequence
\begin{align}
\label{eqn:long exact}
0 & \longrightarrow \Hom(S^{(\nu)},S)\longrightarrow \Hom(S^{(\nu)},S^{(r)}) \longrightarrow \Hom(S^{(\nu)},S^{(r-1)})
\\
\notag
& \longrightarrow \Ext^1(S^{(\nu)},S) {\longrightarrow} \Ext^1(S^{(\nu)},S^{(r)}) \stackrel{\phi}{\longrightarrow} \Ext^1(S^{(\nu)},S^{(r-1)})\longrightarrow0.
\end{align}
Let us describe the action of $\phi$.
For any short exact sequence $0\rightarrow S^{(r)} \xrightarrow{h_1} S^{(\mu)} \xrightarrow{h_2} S^{(\nu)}\rightarrow0$ representing $[\eta]\in \Ext^1(S^{(\nu)},S^{(r)})_{S^{(\mu)}}$, we have a pushout diagram of short exact sequences
\[
\xymatrix{ S\ar@{=}[d] \ar[r]^{f_1} & S^{(r)} \ar[r]^{f_2} \ar[d]^{h_1} & S^{(r-1)} \ar[d]^{g_1} \\
S \ar[r] & S^{(\mu)}  \ar[r] \ar[d]^{h_2} & S^{(\la)} \ar[d]^{g_2} \\
& S^{(\nu)} \ar@{=}[r] & S^{(\nu)}  }
\]
with the rightmost column representing $\phi([\eta])=[\xi]$.

We know that the second row does not split since the first row does not. It follows that $\ell(\la)=\ell(\mu)$. Hence,
\begin{align}
 \label{eq:ext3}
\sum_{\mu:\, \ell(\mu)=\ell(\nu)+1} \big|\Ext^1(S^{(\nu)},S^{(r)})_{S^{(\mu)}}\big|
&= \sum_{\la:\, \ell(\la)=\ell(\nu)+1}\big|\phi^{-1}(\Ext^1(S^{(\nu)},S^{(r-1)})_{S^{(\la)}}) \big|\\
&= \sum_{\la:\, \ell(\la)=\ell(\nu)+1}\big|\Ext^1(S^{(\nu)},S^{(r-1)})_{S^{(\la)}}\big|\cdot|\ker\phi|.
\notag
\end{align}

By the inductive assumption, we have
\begin{align}\label{induction on r-1}
\sum_{\la:\, \ell(\la)
=\ell(\nu)+1}\big|\Ext^1(S^{(\nu)},S^{(r-1)})_{S^{(\la)}}\big|=q^{\sum_{i}\min(r-1,\nu_i)- \ell(\nu)}.
\end{align}
It follows from \eqref{eqn:long exact} that
 \begin{align}\label{order of ker phi}
|\ker\phi|=\frac{\big|\Hom(S^{(\nu)},S^{(r)})\big|\cdot\big|\Ext^1(S^{(\nu)},S)\big|}{\big|\Hom(S^{(\nu)},S)\big|\cdot \big|\Hom(S^{(\nu)},S^{(r-1)})\big|}
= q^{\sum_{i}\min(r,\nu_i)-\sum_{i}\min(r-1,\nu_i)}.
\end{align}
Plugging \eqref{induction on r-1} and \eqref{order of ker phi} into \eqref{eq:ext3}, we obtain
\eqref{eq:Ext1}.
The lemma is proved.
\end{proof}

We record the following curious identity. Recall $\varphi_r (t)$ from \eqref{eq:varphi:r}.
\begin{corollary}
\label{cor:sumext0}
For any $r\geq 1$ and any partition $\nu\neq \emptyset$, we have
\begin{align*}
\sum\limits_{\mu}\varphi_{\ell(\mu)-1}(q)\big|\Ext^1(S^{(\nu)}, S^{(r)})_{S^{(\mu)}}\big|=0.
\end{align*}
\end{corollary}

\begin{proof}
By Lemma \ref{extensions}, we have
\begin{align*}
 \sum\limits_{\mu} & \varphi_{\ell(\mu)-1}(q)\big|\Ext^1(S^{(\nu)}, S^{(r)})_{S^{(\mu)}}\big|\\
&= \sum\limits_{\mu:\, \ell(\mu)=\ell(\nu)+1}\varphi_{\ell(\nu)}(q)\big|\Ext^1(S^{(\nu)},S^{(r)})_{S^{(\mu)}}\big|
+\sum\limits_{\mu:\, \ell(\mu)=\ell(\nu)}\varphi_{\ell(\nu)-1}(q)\big|\Ext^1(S^{(\nu)},S^{(r)})_{S^{(\mu)}}\big|\\
&= \varphi_{\ell(\nu)}(q) \cdot q^{ \sum_{i} \min(r,\nu_i)- \ell(\nu)}
+\varphi_{\ell(\nu)-1}(q)\cdot q^{ \sum_{i} \min(r,\nu_i)}(1-q^{- \ell(\nu)})\\
&= q^{ \sum_{i}\min(r,\nu_i)- \ell(\nu)}
\Big(\varphi_{\ell(\nu)}(q)+\varphi_{\ell(\nu)-1}(q)(q^{\ell(\nu)}-1) \Big)
=0.
\end{align*}
\end{proof}

\subsection{An identity of Hall numbers}

For $r>0$ and a partition $\omega =(i^{m_i})_{i\ge 1}$, where $m_i$ denotes the number of parts in $\omega$ equal to $i$, we introduce new partitions
\[
\omega_+^r =(i^{m_i})_{i>r},
\quad
\omega_-^r =(i^{m_i})_{i<r},
\quad
\omega^r=(r^{m_r}).
\]
It follows that $\sum\limits_{i}\min(r, \omega_i)=|\omega^r_-|+|\omega^r|+r\ell(\omega^r_+)$.

\begin{lemma}
For $r\geq1$ and any partition $\omega=(\omega_1, \omega_2, \ldots)$, we have
\begin{align}
\label{eq:Gaut1}
\sum_{\nu:\, \ell(\nu)=\ell(\omega)}G_{S^{(\nu)}, S^{(r)}}^{S^{(\omega)}} |\Aut(S^{(r)})|
& =q^{\sum\limits_{i}\min(r, \omega_i)-\ell(\omega)}\Big(q^{\ell(\omega^r_+)}-1\Big),
\\
\label{eq:Gaut2}
\sum_{\nu:\, \ell(\nu)=\ell(\omega)-1}G_{S^{(\nu)}, S^{(r)}}^{S^{(\omega)}} |\Aut(S^{(r)})|
& =q^{\sum\limits_{i}\min(r, \omega_i)-\ell(\omega)}\Big(q^{\ell(\omega)}-q^{\ell(\omega^r_+)}-q^{\ell(\omega^r_-)}+1\Big).
\end{align}
\end{lemma}

\begin{proof}
Consider the short exact sequence
\begin{align}
  \label{eq:ses1}
0\longrightarrow S^{(r)}\stackrel{\phi}{\longrightarrow} S^{(\omega)} \stackrel{\psi}{\longrightarrow} S^{(\nu)}\longrightarrow0.
\end{align}
Then $\psi$ is determined by $\phi$ up to an automorphism in $\Aut(S^{(\nu)})$. Denote by  $\phi^{(r)}: S^{(r)}\rightarrow S^{(\omega^r)}$ and $\phi^{(r)}_\pm: S^{(r)}\rightarrow S^{(\omega^r_\pm)}$ the compositions of $\phi$ with the corresponding canonical projections.

Observe that $\phi$ is injective if and only if there exists some component of $\phi^{(r)}_+$ or $\phi^{(r)}$ which is injective, hence the number of such $\phi$'s is given by
\begin{align}  \label{sum1}
\sum_{\nu}G_{S^{(\nu)}, S^{(r)}}^{S^{(\omega)}} |\Aut(S^{(r)})|
&= q^{|\omega^r_-|}(q^{r(\ell(\omega^r_+)+\ell(\omega^r))}-q^{(r-1)(\ell(\omega^r_+)+\ell(\omega^r))})
\\
& =q^{\sum\limits_{i}\min(r, \omega_i)}(1-q^{\ell(\omega^r_-)-\ell(\omega)})
\notag \\
& =q^{\sum\limits_{i}\min(r, \omega_i)-\ell(\omega)}(q^{\ell(\omega)}-q^{\ell(\omega^r_-)}).
\notag
\end{align}

Note $\ell(\nu)=\ell(\omega)$ if and only if $\phi^{(r)}_+$ admits an injective component and $\phi^{(r)}_-, \phi^{(r)}$ contain no surjective components. Thus,
\begin{align*}
\sum_{\nu: \ell(\nu)=\ell(\omega)} G_{S^{(\nu)}, S^{(r)}}^{S^{(\omega)}} |\Aut(S^{(r)})|
&=q^{\sum\limits_{\omega_i\leq r}(\omega_i-1)}(q^{r\ell(\omega^r_+)}-q^{(r-1)\ell(\omega^r_+)})\\
&=q^{r\ell(\omega^r_+)+|\omega^r|+|\omega^r_-|-\ell(\omega^r)-\ell(\omega^r_-)}(1-q^{-\ell(\omega^r_+)})
 \notag \\
&=q^{\sum\limits_{i}\min(r, \omega_i)-\ell(\omega)}(q^{\ell(\omega^r_+)}-1).
\notag
\end{align*}
This proves \eqref{eq:Gaut1}.

Observe that $\ell(\nu)=\ell(\omega)$ or $\ell(\omega)-1$. The second identity in the lemma follows from \eqref{sum1} and \eqref{eq:Gaut1}.
\end{proof}

\subsection{A combinatorial identity}

Set
\[
\sqq =\sqrt{q},
\qquad
[n]_\sqq =\frac{\sqq^n -\sqq^{-n}}{\sqq -\sqq^{-1}}.
\]

The next lemma will be used in the proof of Proposition \ref{prop:QP}.
\begin{lemma}
\label{lem:omega<m}
Let $n\geq1$. For any partition $\mu$ of $n$, we have
\begin{align*}
 &\sum_{1\leq r\leq n-1}\sqq^{n-r}[n-r]_\sqq \Big(q^{\ell(\mu)}+q^{\ell(\mu)-1}
-q^{\ell(\mu)-1+\ell(\mu^r_+)}-q^{\ell(\mu^r_-)}\Big)\cdot q^{\sum\limits_{i}\min(r, \mu_i)-\ell(\mu)}\\
&\qquad\qquad=\sqq^{n}[n]_\sqq(q^{\ell(\mu)-1}-1).
\end{align*}
\end{lemma}

\begin{proof}
Denote $\mu=(i^{m_i})_{i\ge 1}$ and $\ell =\ell(\mu)$. Then $\sum\limits_{i}\min(r, \mu_i)=n -\sum\limits_{i>r} (i-r) m_i$. Hence the formula in the lemma is converted to be ($\sum_{1\leq r\leq n-1}$ is replaced by $\sum_{1\leq r\leq n}$ thanks to $[n-n]_\sqq=0$):
\begin{align}
  \label{eq:id}
& \sum_{1\leq r\leq n}
\sqq^{n-r}[n-r]_\sqq
\Big(q^{\ell}+q^{\ell-1}
-q^{\ell-1+\sum\limits_{i> r}m_i}-q^{\ell -\sum\limits_{i\ge  r}m_i}\Big) q^{n-\ell -\sum\limits_{i>r} (i-r) m_i}
=\sqq^{n}[n]_\sqq(q^{\ell -1}-1).
\end{align}
Since $q^{\ell}+q^{\ell-1}
-q^{\ell-1+\sum\limits_{i> r}m_i}-q^{\ell -\sum\limits_{i\ge  r}m_i}=0$ whenever $r> n$, the summation $\sum_{1\leq r\leq n}$ in \eqref{eq:id} can be replaced by $\sum_{r\ge 1}$.

Denote by $\mu' =(\mu_1', \mu_2', \ldots)$ the conjugate partition of $\mu$. Note that $\mu_r' =\sum_{i\ge r} m_i$.
Denote by $\mu^{(r)}$ the partition obtained from the Young diagram of $\mu$ by removing the first $r$ columns, for $r\ge 0$. Then it follows by definition that $\mu^{(r)} =\emptyset$, for $r\gg 0$, and




\begin{align}  \label{eq:mua}
|\mu^{(r-1)}| = |\mu^{(r)}| +\mu_r', \quad \text{ for }r \ge 1.
\end{align}

In particular, we have $|\mu^{(r)}| =\sum_{i>r} (i-r) m_i$. The identity \eqref{eq:id} admits the following equivalent reformulation:

\begin{align}
  \label{eq:id2}
& \sum_{r\ge 1}
\sqq^{n-r}[n-r]_\sqq q^{n -|\mu^{(r)}|}
\big(1+q^{-1}
-q^{-1+\mu_{r+1}'}-q^{-\mu_r'} \big)
=\sqq^{n}[n]_\sqq(q^{\ell -1}-1).
\end{align}

Let us prove \eqref{eq:id2}. Recall $q=\sqq^2$. Noting $|\mu^{(0)}|=n$ and $|\mu^{(1)}|=n-\ell$, applying \eqref{eq:mua} twice we have

\begin{align*}
\text{LHS } \eqref{eq:id2}
& =\sum_{r\ge 1} \sqq^{n-r}[n-r]_\sqq q^{n}
\big( (1+q^{-1}) q^{-|\mu^{(r)}|}
-q^{-1-|\mu^{(r+1)}|}-q^{-|\mu^{(r-1)}|} \big)
\\
& =\sum_{r\ge 1} \sqq^{n-r}[n-r]_\sqq q^{n}
\big( q^{ -|\mu^{(r)}|} -q^{-|\mu^{(r-1)}|} \big)
\\
&\qquad - \sum_{r\ge 1} \sqq^{n-r-2}[n-r]_\sqq q^{n}
\big( q^{-|\mu^{(r+1)}|}  - q^{-|\mu^{(r)}|}\big).
\end{align*}

By shifting the index $r\mapsto r+1$ for the partial sum $\sum_{r\ge 2}$ in the first summand above, we rewrite the above as

\begin{align*}
\text{LHS } \eqref{eq:id2}
& =
\sqq^{n-1}[n-1]_\sqq q^{n}
\big( q^{-(n-\ell)} -q^{- n} \big)
\\
& \quad +\sum_{r\ge 1} \big(\sqq^{n-r-1}[n-r-1]_\sqq -\sqq^{n-r-2}[n-r]_\sqq \big) q^{n}
\big( q^{-|\mu^{(r+1)}|}  - q^{-|\mu^{(r)}|}\big)
\\%
& = \sqq^{n-1}[n-1]_\sqq q^{n}
\big( q^{\ell-n} -q^{- n} \big)
-\sum_{r\ge 1} \sqq^{-1} q^{n} \big( q^{-|\mu^{(r+1)}|}  - q^{-|\mu^{(r)}|}\big)
\\%
& = \sqq^{n-1}[n-1]_\sqq q^{n}
\big( q^{\ell -n} -q^{- n} \big)
- \sqq^{-1} q^{n} ( 1 - q^{-(n-\ell)} )
\\
&= \sqq^{n}[n]_\sqq(q^{\ell -1}-1).
\end{align*}

This proves \eqref{eq:id2} and hence the lemma.
\end{proof}


\section{Transition relations in the $\imath$Hall algebra}
 \label{sec:identities}

In this section, we work with the $\imath$Hall algebra of the Jordan quiver $\iRH(\bfk \QJ) \otimes_{\Q} \Q(\sqq)$, where $\bfk =\mathbb F_q$ and $\sqq=\sqrt{q}.$ We shall determine the relations among the generating functions $\haE(z), \haH(z), \haT(z)$ and $\haP(z)$ (see \eqref{def:haE}, \eqref{def:haH}, \eqref{def:haT}, and \eqref{def:haP}) of several generating sets for the $\imath$Hall algebra.

\subsection{Euler's identity}

(In this subsection, $q$ can also be interpreted as a formal variable.)

For $n\in \N \cup \{\infty \}$, we denote by
\[
(x;q)_n =(1-x)(1-qx)\ldots (1-q^{n-1} x).
\]
We have the following identity due to Euler:
\begin{align}
 \label{eq:Euler}
\exp_q\Big(\frac{x}{1-q} \Big):= \sum_{r\geq0} \frac{x^r}{(1-q)(1-q^2)\ldots (1-q^r)}
= \frac{1}{(x; q)_\infty}.
\end{align}
It is elementary to rewrite $\frac{1}{(x; q)_\infty}$ in terms of the usual exponential function:
\[
\frac{1}{(x; q)_\infty} = \exp \Big (\sum_{k\ge 1} \frac{x^k}{k (1-q^k)} \Big).
\]
Hence we have
\begin{align}  \label{exp=exp}
\exp_q \big(x/(1-q) \big) = \exp \Big (\sum_{k\ge 1} \frac{x^k}{k (1-q^k)} \Big).
\end{align}

\subsection{$\haH(z)$ vs $\haE(z)$}
We prepare some lemmas before formulating the relation between $\haH(z)$ vs $\haE(-z)$.

\begin{lemma}
\label{lem:sumsemisimple}
For $r\ge 1$ and any partition $\mu$ with $\ell(\mu)\geq r$, we have
\begin{align*}
\sum_{\la} G_{(1^{r}),\la}^\mu(\bq)  |\Aut(S^{(1^r)})| =(q^{\ell(\mu)}-1) (q^{\ell(\mu)}-q)\cdots (q^{\ell(\mu)}-q^{r-1}).
\end{align*}
\end{lemma}

\begin{proof}
By definition of $G_{(1^{r}),\la}^\mu(\bq) $ we have
\begin{align*}
\sum_{\la} G_{(1^{r}),\la}^\mu(\bq)  |\Aut(S^{(1^r)})|=& \big|\{d\in\Hom(S^{(\mu)},S^{(1^r)})\mid d \text{ is surjective}\} \big|.
\end{align*}
Let $f:S^{(\mu)}\rightarrow {\rm top}(S^{(\mu)})=S^{(1^{\ell(\mu)})}$ be the natural projection. For any epimorphism $d: S^{(\mu)}\rightarrow S^{(1^r)}$, there exists a unique epimorphism
$g:S^{(1^{\ell(\mu)})}\rightarrow S^{(1^r)}$ such that $d=gf$. Thus,
\begin{align*}
\big|\{d\in\Hom(S^{(\mu)},S^{(1^r)})\mid d \text{ is surjective}\} \big|
=& \big|\{g\in\Hom(S^{(1^{\ell(\mu)})},S^{(1^r)})\mid g \text{ is surjective}\} \big|
\\
=&(q^{\ell(\mu)}-1) (q^{\ell(\mu)}-q) \cdots (q^{\ell(\mu)}-q^{r-1}).
\end{align*}
\end{proof}

\begin{lemma}
(cf., e.g., \cite[Lemma 5.4]{LW20})
\label{lem:qbinom1}
Let $p \in \Z_{\geq1}$.
Let $d\in \Z$ be such that $|d| \le p-1$ and $d\equiv p-1\pmod 2$. Then
$\sum\limits_{n=0}^{p}(-1)^nv^{-dn}\qbinom{p}{n}=0$.
\end{lemma}
Recall $\haE(z), \haH(z)$ from \eqref{def:haE} and \eqref{def:haH}. Recall $\exp_q (\frac{x}{1-q} )$ from \eqref{eq:Euler}.

\begin{proposition}
 \label{prop:HE}
We have
\begin{align*}
\haH(z)\haE(-z) =& \exp_q \Big(\frac{ [K_S] z^2}{1-q} \Big).
\end{align*}
\end{proposition}

\begin{proof}
Recall from \eqref{Pieri:v} that
\begin{align*}
[S^{(\la)}]*[S^{(1^a)}]
&= \sum_{r=0}^a\sum_{\omega,\nu} G_{\nu,(1^{a-r})}^\omega(\bq) G_{(1^r),\nu}^\la(\bq) \frac{|\Aut(S^{(\nu)})| \cdot|\Aut(S^{(1^{a-r})})| }{|\Aut(S^{(\omega)})|}
\\\notag
&\qquad\qquad\qquad \times (\bq^a-1)\cdots(\bq^a-\bq^{r-1}) [S^{(\omega)}]*[K_S]^r.
\end{align*}
Hence we have
\begin{align*}
&\haH(z)\haE(-z)
\\
= &\sum_{a\geq0} \sum_{b\geq0} \sum\limits_{\la\vdash b}(-1)^a\sqq^{a(a-1)} \frac{[S^{(\la)}]}{|\Aut(S^{(\la)})|}* \frac{[S^{(1^a)}]}{|\Aut(S^{(1^a)})|}z^{a+b}
\\
=&\sum_{a\geq0} \sum_{b\geq0} \sum\limits_{\la\vdash b}(-1)^a\sqq^{a(a-1)} \sum_{\omega,r\geq0}\sum_{\nu}G_{\nu,(1^{a-r})}^\omega(\bq) G_{(1^r),\nu}^\la(\bq) \frac{|\Aut(S^{(\nu)})| \cdot|\Aut(S^{(1^{a-r})})| }{|\Aut(S^{(\omega)})|\cdot|\aut(S^{(\la)})|\cdot|\aut(S^{(1^a)})|}
\\\notag
&\qquad\qquad\qquad\qquad \times
(\bq^a-1)\cdots(\bq^a-\bq^{r-1}) [S^{(\omega)}]*[K_S]^rz^{a+b}
\\
=&\sum_{a\geq0} \sum_{b\geq0} (-1)^a\sqq^{a(a-1)} \sum_{\omega,r\geq0}\sum_{\nu\vdash (b-r)}G_{\nu,(1^{a-r})}^\omega(\bq)  \frac{|\Aut(S^{(1^{a-r})})| }{|\Aut(S^{(\omega)})|\cdot|\aut(S^{(1^a)})|\cdot |\aut(S^{(1^r)})|}
\\\notag
&\qquad\qquad\qquad\qquad \times
(\bq^a-1)\cdots(\bq^a-\bq^{r-1}) [S^{(\omega)}]*[K_S]^rz^{a+b}.
\end{align*}
Here the last equality follows by the Riedtmann-Peng formula
\begin{align*}
\sum_{\la \vdash b} \frac{ G_{(1^r),\nu}^\la(\bq)   |\aut(S^{(\nu)})|}{|\aut(S^{(\la)})|}=&\frac{|\Ext^1(S^{(1^{r})},S^{(\nu)} )|}{|\Hom(S^{(1^{r})},S^{(\nu)})|\cdot|\aut(S^{(1^{r})})|}
=\frac{1}{|\aut(S^{(1^{r})})|}.
\end{align*}
Observe that $(\bq^a-1)\cdots (\bq^a-\bq^{r-1}) \frac{1}{|\aut(S^{(1^r)})|}=\sqq^{r(a-r)}\qbinom{a}{r}_\sqq$. Hence
\begin{align*}
\haH(z)\haE(-z)
=&\sum_{a\geq0} \sum_{b\geq0} (-1)^a\sqq^{a(a-1)} \sum_{\omega,r\geq0}\sum_{\nu\vdash (b-r)}
G_{\nu,(1^{a-r})}^\omega(\bq)  \frac{|\Aut(S^{(1^{a-r})})| }{|\Aut(S^{(\omega)})|\cdot|\aut(S^{(1^a)})|}
 \\
 &\qquad\qquad\qquad\qquad\qquad\qquad \times
 \sqq^{r(a-r)}\qbinom{a}{r}_\sqq [S^{(\omega)}]*[K_S]^rz^{a+b}.
\end{align*}

Fix $\omega$ and $r$ in the following. If $\omega=\emptyset$, then $a=r=b$ and $\nu=\emptyset$. The coefficient of $[K_S]^rz^{2r}$ is equal to
\begin{align*}
(-1)^r\frac{\sqq^{r(r-1)}}{|\aut(S^{(1^r)})|}= \frac{1}{(1-q^r)(1-q^{r-1})\cdots (1-q)}.
\end{align*}
Now assume that $\omega\neq \emptyset$. Then $G_{\nu,(1^{a-r})}^\omega(\bq)\neq 0$ implies that $a+b=|\omega|+2r$ and $r\leq a\leq \ell(\omega)+r$. Hence, by Lemma~\ref{lem:sumsemisimple}, the coefficient of $[S^{(\omega)}]*[K_S]^rz^{a+b}$ (with the factor $\frac{1}{|\Aut(S^{(\omega)})|}$ ignored) is equal to
\begin{align*}
&\sum_{r\leq a\leq \ell(\omega)+r} (-1)^a\sqq^{a(a-1)}\sqq^{r(a-r)}\qbinom{a}{r}_\sqq \frac{1}{|\Aut(S^{(1^a)})|} \big(\sum_\nu G^\omega_{\nu,(1^{a-r})}(\bq)|\Aut(S^{(1^{a-r})})|\big)
\\
=&(-1)^r\sqq^{r(r-1)} \frac{1}{(q^r-1)\cdots (q^r-q^{r-1})} \\
& \quad +\sum_{r< a\leq \ell(\omega)+r} (-1)^a\sqq^{a(a-1)}\sqq^{r(a-r)}\qbinom{a}{r}_\sqq \frac{(q^{\ell(\omega)}-1)\cdots (q^{\ell(\omega)}-q^{a-r-1})}{(q^a-1)\cdots (q^a-q^{a-1}) }
\\
=&\frac{(-1)^r\sqq^{-\binom{r+1}{2}}}{[r]_\sqq^!(\sqq-\sqq^{-1})^r} +  \sum_{r< a\leq \ell(\omega)+r} (-1)^a \frac{\sqq^{r(a-r)+(a-r)\ell(\omega)+\binom{r+1}{2}-a(r+1)}}{[r]_\sqq^!(\sqq-\sqq^{-1})^r} \frac{[\ell(\omega)]_\sqq^!}{[a-r]_\sqq^![\ell(\omega)-a+r]_\sqq^!}
\\
=& \frac{(-1)^r\sqq^{-\binom{r+1}{2}}}{[r]_\sqq^!(\sqq-\sqq^{-1})^r} \Big( 1+ \sum_{r< a\leq \ell(\omega)+r} (-1)^{a-r} \sqq^{(a-r)(\ell(\omega)-1)} \frac{[\ell(\omega)]_\sqq^!}{[a-r]_\sqq^![\ell(\omega)-a+r]_\sqq^!}\Big)
\\
=& \frac{(-1)^r\sqq^{-\binom{r+1}{2}}}{[r]_\sqq^!(\sqq-\sqq^{-1})^r} \sum_{0\leq c\leq \ell(\omega)} (-1)^c \sqq^{c(\ell(\omega)-1)} \qbinom{\ell(\omega)}{c}_\sqq
=0,
\end{align*}
where the last 2 equalities used a substitution $c=a-r$ and Lemma \ref{lem:qbinom1}. Then the desired formula follows.
\end{proof}

\subsection{$\haT(z)$ vs $\haE(z)$}

\begin{lemma}
\label{lem:semiind}
Assume $a,b>0$. The following identities hold in $\iRH(\bfk \QJ)$:
for $a>1$,
\begin{align}
&[S^{(a)}]*[S^{(1^b)}]=q^{-b}[ S^{(a)}\oplus S^{(1^b)}]+(1-q^{-b}) [S^{(a+1)}\oplus S^{(1^{b-1})}]
\\\notag
&\qquad+(q-q^{1-b})[S^{(a-1)}\oplus S^{(1^{b-1})}]*[K_S]+(q-q^{1-b})(q^{b-1}-1)[S^{(a)}\oplus S^{(1^{b-2})}]*[K_S];
\end{align}
for $a=1$,
\begin{align}
[S]*[S^{(1^b)}]=&q^{-b}[ S^{(1^{b+1})}]+(1-q^{-b}) [S^{(1^{b-1})}\oplus S^{(2)}]+(q^b-1)[S^{(1^{b-1})}]*[K_S].
\end{align}

\end{lemma}

\begin{proof}
It follows by a direct computation, or by Lemma \ref{lem:HallfomrulaJordan}.
\end{proof}

Inspired by \cite{LRW20}, we define
$$\haT_{m}=\frac{[S^{(m)}]}{\sqq-\sqq^{-1}}.$$
Recall $\haT(z)$ from \eqref{def:haT}. We further denote
\begin{align}
 \label{def:haT2}
\haT(z)=\sum_{r\geq0} (\sqq-\sqq^{-1})\haT_r z^r
=\sum_{r\geq0}[S^{(r)}]z^r.
\end{align}

\begin{proposition}
  \label{prop:QE}
  We have
\begin{align*}
\haT(z)=& \big(1-q[K_S] z^2 \big) \frac{\haE(-z)}{\haE(- qz)},
\\
\haT(z)=& \big(1-[K_S] z^2 \big)^{-1} \frac{\haH(qz)}{\haH(z)}.
\end{align*}
\end{proposition}

\begin{proof}
The 2 formulas in the proposition are equivalent by Proposition~\ref{prop:HE} with the help of Euler's identity \eqref{eq:Euler}. Let us prove the first formula.

To that end, we shall prove the following equivalent formula:
\begin{align} \label{QE=E}
\haT(z)\haE(- qz)=&(1-q[K_S] z^2)\haE(-z).
\end{align}

By definition we have
\begin{align}\label{haT times haE}
\haT(z)\haE(- qz)
=\sum_{a\geq0}\sum_{b\geq0} (-q)^b\sqq^{b(b-1)}   \frac{[S^{(a)}] *[S^{(1^b)}]}{|\Aut(S^{(1^b)})|}z^{a+b}.
\end{align}
By Lemma \ref{lem:semiind}, only the terms $[S^{(a)}\oplus S^{(1^{b})}]z^{a+b}$ and $[S^{(a)}\oplus S^{(1^{b})}]*[K_S]z^{a+b+2}$ can appear on the right-hand side of \eqref{haT times haE}.

For $a>1$, the term $[S^{(a)}\oplus S^{(1^{b})}]z^{a+b}$ comes from $[S^{(a)}] *[S^{(1^b)}]$ or $[S^{(a-1)}] *[S^{(1^{b+1})}]$, and hence its coefficient is equal to
\begin{align*}
&(-q)^b \sqq^{b(b-1)} \frac{q^{-b}}{|\Aut(S^{(1^b)})|} +(-q)^{b+1} \sqq^{b(b+1)} \frac{1-q^{-(b+1)}}{|\Aut(S^{(1^{b+1})})|}=0,
\end{align*}
thanks to $|\Aut(S^{(1^{b+1})})|=|\Aut(S^{(1^{b})})|\cdot q^{b}(q^{b+1}-1).$

Similarly, the term $[S^{(a)}\oplus S^{(1^{b})}]*[K_S]z^{a+b+2}$ comes from $[S^{(a+1)}] *[S^{(1^{b+1})}]$ or $[S^{(a)}] *[S^{(1^{b+2})}]$, and hence its coefficient is equal to
\begin{align*}
&(-q)^{b+1}\sqq^{b(b+1)}  \frac{q-q^{1-(b+1)}}{|\Aut(S^{(1^{b+1})})|} +(-q)^{b+2} \sqq^{(b+1)(b+2)} \cdot \frac{(q-q^{1-(b+2)})(q^{(b+2)-1}-1)}{|\Aut(S^{(1^{b+2})})|}=0.
\end{align*}
Hence, only the terms $[S^{(1^{b})}]z^{b}$ and $[S^{(1^{b})}]*[K_S]z^{b+2}$ can appear nontrivially on RHS of \eqref{haT times haE}.

For any $b\geq 0$, the term $[S^{(1^{b})}]z^{b}$ comes from $[S^{(0)}]*[S^{(1^b)}]$ or $[S] *[S^{(1^{b-1})}]$, and hence its coefficient is equal to
\begin{align*}
(-q)^b \sqq^{b(b-1)}  \frac{1}{|\Aut(S^{(1^b)})|} +(-q)^{b-1} \sqq^{(b-1)(b-2)} \frac{q^{-(b-1)}}{|\Aut(S^{(1^{b-1})})|}
= (-1)^b \frac{\sqq^{b(b-1)}}{|\Aut(S^{(1^b)})|}.
\end{align*}
The term $[S^{(1^{b})}]*[K_S]z^{b+2}$ comes from $[S] *[S^{(1^{b+1})}]$ and $[S^{(2)}] *[S^{(1^{b})}]$, and hence its coefficient is equal to
\begin{align*}
(-q)^{b+1}\sqq^{b(b+1)} \frac{q^{b+1}-1}{|\Aut(S^{(1^{b+1})})|}+(-q)^{b}\sqq^{b(b-1)} \frac{q-q^{1-b}}{|\Aut(S^{(1^{b})})|}
= (-1)^{b+1} \frac{q \sqq^{b(b-1)}}{|\Aut(S^{(1^{b})})|}.
\end{align*}
Therefore we have
\begin{align*}
\haT(z)\haE(- qz)
&= \sum_{b\geq0}(-1)^b \frac{\sqq^{b(b-1)}}{|\Aut(S^{(1^b)})|}[S^{(1^{b})}]z^{b}+\sum_{b\geq0}(-1)^{b+1} \frac{q \sqq^{b(b-1)}}{|\Aut(S^{(1^{b})})|}[S^{(1^{b})}]*[K_S]z^{b+2}\\
&= \widehat{E}(-z) -q[K_S] z^2\widehat{E}(-z)\\
&= (1-q[K_S] z^2)\haE(-z).
\end{align*}
This proves the proposition.
\end{proof}

\subsection{$\haT(z)$ vs $\haP(z)$}

Recall $\haT(z), \haT_r$ from \eqref{def:haT} and \eqref{def:haT2}. Define $\widehat{T}_m\in \iRH(\bfk \QJ)$ by
\begin{align}\label{eq:ThT}
\haT(z)  =\exp\big( (\sqq-\sqq^{-1}) \sum_{m\ge 1} \widehat{T}_m z^m \big).
\end{align}

\begin{lemma}
The identity \eqref{eq:ThT} is equivalent to
\begin{align}
 \label{eq:Tr2}
m\haT_m=(\sqq-\sqq^{-1})\sum_{1\leq l\leq m}l \widehat{T}_l*\haT_{m-l},\quad  \forall m\geq1,
\end{align}
where $\haT_0=\frac{1}{\sqq-\sqq^{-1}}$.
\end{lemma}

\begin{proof}
First assume \eqref{eq:ThT} holds. Applying the differentiation $\frac{d}{dz}$ to \eqref{eq:ThT}, we have
\begin{align}\label{DerivationFormula}
\sum_{m\geq 1} (\sqq-\sqq^{-1}) m \haT_{m} z^{m-1}
&= \Big((\sqq-\sqq^{-1})\sum_{l\ge 1} l \widehat{T}_l z^{l-1} \Big)  \Big( 1+ \sum_{m\geq 1} (\sqq-\sqq^{-1})\haT_{m} z^m\Big).
\end{align}
Now \eqref{eq:Tr2} follows by comparing coefficients of $z^{m-1}$ on both sides.

Conversely, assume \eqref{eq:Tr2} holds. That is, \eqref{DerivationFormula} holds. 
We obtain
$$\frac{d}{dz} \haT (z) =\big((\sqq-\sqq^{-1})\sum_{l\geq1} l \widehat{T}_l z^{l-1} \big)  \haT(z),$$
and hence, $\frac{d}{dz} \big(\ln \haT (z) \big) = \frac{d}{dz} \big((\sqq-\sqq^{-1})\sum_{l\geq1}  \widehat{T}_l z^{l} \big)$. Then, $\ln \haT (z) = (\sqq-\sqq^{-1})\sum_{l\geq1}  \widehat{T}_l z^{l}$ (as both sides have constant term 0), whence \eqref{eq:ThT}.
%
\end{proof}

Recall $\haP(z)$ and $\haP_r$ from \eqref{def:haP}, for $r\ge 1$.


\begin{proposition}
\label{prop:QP}
We have
\begin{align}\label{eq:haTandhaP}
\widehat{T}_m &=\sqq^m \frac{[m]_\sqq}{m} \haP_m -\delta_{m, ev} \sqq^{m/2} \frac{[m/2]_\sqq}{m} [K_S]^{m/2}  \qquad (m\geq1),
\\
\label{eq:QP}
\haT(z)
&=  \exp \Big(\sum_{r\geq1} \frac{q^{r}-1}{r} \haP_r z^{r} \Big)  \exp \Big( \sum_{k \geq1} \frac{1-q^{k}}{2k} [K_S]^k z^{2k} \Big).
\end{align}
\end{proposition}

\begin{proof}
The identity \eqref{eq:QP} follows from \eqref{eq:ThT} and \eqref{eq:haTandhaP} by a formal manipulation. So it remains to prove \eqref{eq:haTandhaP}. We proceed by induction on $m$.
It is trivial for $m=1$.

Let  $m\geq 2$. By applying \eqref{eq:Tr2}, the inductive assumption and then Lemma \ref{lem:HallfomrulaJordan}, we have
\begin{align*}
m\widehat{T}_m &= m\haT_{m}-(\sqq-\sqq^{-1})\sum_{1\leq l\leq m-1}l \widehat{T}_l*\haT_{m-l}
\\
&= m\haT_{m}-(\sqq-\sqq^{-1})\sum_{1\leq l\leq m-1}(\sqq^l [l]_\sqq \haP_l -\delta_{l, ev} \sqq^{l/2} [l/2]_\sqq [K_S]^{l/2})*\haT_{m-l}
\\
&= m\haT_{m}-\sum_{1\leq l\leq m-1}\sqq^l [l]_\sqq \sum\limits_{\la\vdash l} \varphi_{\ell(\la)-1}(q)\frac{[S^{(\la)}]*[S^{(m-l)}]}{|\Aut(S^{(\la)})|}
 \\
 &\qquad\qquad +\sum_{1\leq l\leq m-1}\delta_{l, ev} \sqq^{l/2} [l/2]_\sqq [K_S]^{l/2}*[S^{(m-l)}]
\\
&= m\haT_{m}-\sum_{1\leq l\leq m-1}\sqq^l [l]_\sqq \sum\limits_{\la\vdash l} \varphi_{\ell(\la)-1}(q)\frac{1}{|\Aut(S^{(\la)})|}\\
&\quad \times \sum_{\omega,\nu,r\geq0}
\frac{|\Ext^1(S^{(\nu)}, S^{(m-l-r)})_{S^{(\omega)}}|\cdot|\Ext^1(S^{(r)},S^{(\nu)})_{S^{(\la)}}| \cdot |\Aut(S^{(\la)})| } {|\Hom(S^{(\nu)},S^{(m-l-r)})|\cdot|\Hom(S^{(r)},S^{(\nu)})|\cdot|\Aut(S^{(\nu)})|} [S^{(\omega)}]*[K_S]^r\\
&\qquad\qquad +\sum_{1\leq l\leq m-1}\delta_{l, ev} \sqq^{l/2} [l/2]_\sqq [K_S]^{l/2}*[S^{(m-l)}].
\end{align*}
By Corollary \ref{cor:sumext0}, the second summand of the above formula is zero unless either $r=0$ or $\nu=\emptyset$. In particular, if $\nu=\emptyset$, then $\la=(r), r=l$ and $\omega=(m-2l)$. So we have
\begin{align}
\label{eq:HmThetam}
m\widehat{T}_m
&= m\haT_{m}-\sum_{1\leq l\leq m-1}\sqq^l [l]_\sqq \sum\limits_{\la\vdash l} \varphi_{\ell(\la)-1}(q)
 \sum_{\omega}
\frac{|\Ext^1(S^{(\la)}, S^{(m-l)})_{S^{(\omega)}}|} {|\Hom(S^{(\la)},S^{(m-l)})|\cdot |\Aut(S^{(\la)})|} [S^{(\omega)}]\\\notag
&\qquad
-\sum_{1\leq l\leq m-1}\sqq^l [l]_\sqq [S^{(m-2l)}]*[K_S]^l +\sum_{1\leq l\leq m-1}\delta_{l, ev} \sqq^{l/2} [l/2]_\sqq [K_S]^{l/2}*[S^{(m-l)}].
\end{align}

Observe that
\begin{align*}
&-\sum_{1\leq l\leq m-1}\sqq^l [l]_\sqq [S^{(m-2l)}]*[K_S]^l +\sum_{1\leq l\leq m-1}\delta_{l, ev} \sqq^{l/2} [l/2]_\sqq [K_S]^{l/2}*[S^{(m-l)}]\\
=&-\sum_{1\leq l\leq \frac{m}{2}}\sqq^l [l]_\sqq [S^{(m-2l)}]*[K_S]^l +\sum_{1\leq l\leq \frac{m-1}{2}} \sqq^{l} [l]_\sqq [K_S]^{l}*[S^{(m-2l)}]\\
=&-\delta_{m, ev} \sqq^{m/2} \frac{[m/2]_\sqq}{m} [K_S]^{m/2}.
\end{align*}

Hence the proof of \eqref{eq:haTandhaP} reduces to the validity of the following identity:
\begin{align}\label{haH'_m}
\sqq^m[m]_\sqq\haP_m-m\haT_m+
\sum\limits_{1\leq l\leq m-1}\sqq^l [l]_\sqq \sum\limits_{\la\vdash l} \varphi_{ \ell(\la)-1}(q)
 {\textstyle \sum\limits_{\omega}
\frac{|\Ext^1(S^{(\la)}, S^{(m-l)})_{S^{(\omega)}}|} {|\Hom(S^{(\la)},S^{(m-l)})|\cdot |\Aut(S^{(\la)})|}
}
 [S^{(\omega)}]=0.
\end{align}

By Riedtmann-Peng formula, we have
$$\frac{|\Ext^1(S^{(\la)}, S^{(m-l)})_{S^{(\omega)}}|} {|\Hom(S^{(\la)},S^{(m-l)})|\cdot|\Aut(S^{(\la)})|} =\frac{G_{S^{(\la)}, S^{(m-l)}}^{S^{(\omega)}}|\Aut(S^{(m-l)})|
} {|\Aut(S^{(\omega)})|}.$$

Hence the identity (\ref{haH'_m}) is equivalent to the following
\begin{align}  \label{haH'_m2}
\sqq^m[m]_\sqq\haP_m-m\haT_m+\sum\limits_{\omega\vdash m}\sum_{1\leq l\leq m-1}\sqq^l [l]_\sqq\sum\limits_{\la} \varphi_{\ell(\la)-1}(q)
\frac{G_{S^{(\la)}, S^{(m-l)}}^{S^{(\omega)}}|\Aut(S^{(m-l)})|
} {|\Aut(S^{(\omega)})|} [S^{(\omega)}]=0.
\end{align}

We shall prove the identity \eqref{haH'_m2} by considering the coefficients of $[S^{(\omega)}]$ in 2 separated cases:

\emph{Case (a):} $\underline{\ell(\omega)=1}$. Then $\omega= (m)$ and $\la= (l)$. Hence the coefficient of $[S^{(\omega)}]$ is equal to
\begin{align*}
\sqq^m[m]_\sqq & \cdot \frac{1}{|\Aut(S^{(m)})|}-\frac{m}{\sqq-\sqq^{-1}}+\sum_{1\leq l\leq m-1}\sqq^l [l]_\sqq\frac{|\Aut(S^{(m-l)})|}{|\Aut(S^{(m)})|}
\\
=& \frac{q^m-1 }{\sqq-\sqq^{-1}}\cdot\frac{1}{q^{m}(1-q^{-1})}- \frac{m}{\sqq-\sqq^{-1}}+\sum_{1\leq l\leq m-1} \frac{q^l-1 }{\sqq-\sqq^{-1}}\cdot\frac{q^{m-l}(1-q^{-1})}{q^{m}(1-q^{-1})}
\\
=&  \frac{1-q^{-m}}{(\sqq-\sqq^{-1})(1-q^{-1})}- \frac{m}{\sqq-\sqq^{-1}}+\sum_{1\leq l\leq m-1} \frac{1-q^{-l}}{\sqq-\sqq^{-1}} =0.
\end{align*}

\emph{Case (b):} $\underline{\ell(\omega)>1}$. Then the coefficient of $[S^{(\omega)}]$ is equal to
\begin{align*}
\sum_{1\leq l\leq m} & \frac{\sqq^l [l]_\sqq}{{|\Aut(S^{(\omega)})|}} \sum\limits_{\la\vdash l} \varphi_{\ell(\la)-1}(q)G_{S^{(\la)}, S^{(m-l)}}^{S^{(\omega)}}|\Aut(S^{(m-l)})|
\\
=&\sum_{1\leq l\leq m} \frac{\sqq^l [l]_\sqq}{|\Aut(S^{(\omega)})|}\Big(\varphi_{\ell(\omega)-1}(q)\sum\limits_{\la: \ell(\la)=\ell(\omega)}
G_{S^{(\la)}, S^{(m-l)}}^{S^{(\omega)}} {|\Aut(S^{(m-l)})|}
\\&\qquad\qquad\qquad \qquad+{\varphi_{\ell(\omega)-2}(q)}
\sum\limits_{\la: \ell(\la)=\ell(\omega)-1}
G_{S^{(\la)}, S^{(m-l)}}^{S^{(\omega)}} {|\Aut(S^{(m-l)})|} \Big)
\\
\stackrel{(*)}{=}
&\; \frac{\sqq^m [m]_\sqq \varphi_{\ell(\omega)-1}(q)}{|\Aut(S^{(\omega)})|}
\\  \qquad \qquad\qquad\qquad
&+\sum_{1\leq a\leq m-1}\frac{\sqq^{m-a} [m-a]{_\sqq}}{|\Aut(S^{(\omega)})|} \Big({\varphi_{\ell(\omega)-1}(q)}
\big(q^{\ell(\omega^a_+)}-1\big)q^{\sum\limits_{i}\min(a, \omega_i)-\ell(\omega)}
\\& \qquad \qquad  +{\varphi_{\ell(\omega)-2}(q)}\big(q^{\ell(\omega)}-q^{\ell(\omega^a_+)}-q^{\ell(\omega^a_-)}+1\big) q^{\sum\limits_{i}\min(a, \omega_i)-\ell(\omega)}  \Big)
\\
=& \frac{\sqq^m [m]_\sqq  \varphi_{\ell(\omega)-1}(q)}{|\Aut(S^{(\omega)})|}+\frac{ {\varphi_{\ell(\omega)-2}(q)}}{|\Aut(S^{(\omega)})|} \sum_{1\leq a\leq m-1}\sqq^{m-a}[m-a]{_{\sqq}}q^{\sum\limits_{i}\min(a, \omega_i)-\ell(\omega)}
\\&\qquad\qquad\qquad \times \Big((1-q^{\ell(\omega)-1})\big(q^{\ell(\omega^a_+)}-1\big)+\big(q^{\ell(\omega)}-q^{\ell(\omega^a_+)}-q^{\ell(\omega^a_-)}+1\big) \Big)
\\
=& \frac{\sqq^m [m]_\sqq  \varphi_{\ell(\omega)-1}(q)}{|\Aut(S^{(\omega)})|}+\frac{ {\varphi_{\ell(\omega)-2}(q)}}{|\Aut(S^{(\omega)})|} \sum_{1\leq a\leq m-1}\sqq^{m-a}[m-a]{_{\sqq}}q^{\sum\limits_{i}\min(a, \omega_i)-\ell(\omega)}
\\&\qquad\qquad\qquad \times \Big(q^{\ell(\omega)}+q^{\ell(\omega)-1}
-q^{\ell(\omega)-1+\ell(\omega^a_+)}-q^{\ell(\omega^a_-)}\Big)
\\=&\frac{\sqq^m [m]_\sqq \varphi_{\ell(\omega)-1}(q)}{|\Aut(S^{(\omega)})|}+\frac{ {\varphi_{\ell(\omega)-2}(q)}}{|\Aut(S^{(\omega)})|} \cdot \sqq^{m}[m]_\sqq (q^{\ell(\omega)-1}-1)
=0,
\end{align*}
where the second last equality uses Lemma \ref{lem:omega<m} and the equation ($*$) uses a substitution $a=m-l$.
This finishes the proof of the proposition.
\end{proof}

\subsection{$\haH(z)$ or $\haE(z)$ vs $\haP(z)$}

Recall $\haE(z)$,  $\haH(z)$ and  $\haP(z)$ from \eqref{def:haE}, \eqref{def:haH}, \eqref{def:haP}, respectively.

\begin{proposition}
\label{prop:EHP}
We have
\begin{align*}
\haE(z)
&=\exp \Big(\sum_{r\geq1} \frac{(-1)^{r-1}}{r} \haP_r z^{r} \Big) \exp_q \Big(\frac{ [K_S] z^2}{1-q} \Big)^{\frac12},
\\
\haH(z)
&= \exp \Big( \sum_{r\geq1} \frac{1}{r} \haP_r z^{r} \Big) \exp_q \Big(\frac{ [K_S] z^2}{1-q} \Big)^{\frac12}.
\end{align*}
\end{proposition}

\begin{proof}
These 2 formulas are equivalent by Proposition~\ref{prop:HE}, and so we only need to prove the first one for $\haE(z)$. Proposition~\ref{prop:QE} can be rewritten as
\begin{align*}
\frac{\haE(- z)}{\haE(-qz)} =& \frac{\haT(z)}{1-q[K_S] z^2}.
\end{align*}
Replacing $z$ by $q^{m}z$ for $m\ge 0$ above, we obtain
\begin{align}  \label{eq:EQK}
\haE(- z) & =
\prod_{m\ge 0} \frac{\haE(- q^{m} z)}{\haE(-q^{m+1}z)}
= \frac{\prod_{m\ge 0} \haT(q^{m} z)}{\prod_{m\ge 0} (1-q^{2m+1}[K_S] z^2)}.
\end{align}
It follows by \eqref{eq:QP} that
\begin{align*}
\haT(z) &= \frac{\exp \Big(\sum_{r\geq1} \frac{q^r}{r} \haP_r z^{r} \Big)}{\exp \Big(\sum_{r\geq1} \frac{1}{r} \haP_r z^{r} \Big) }
\exp \Big(\sum_{k \geq1} \frac{1-q^{k}}{2k} [K_S]^k z^{2k} \Big).
\end{align*}
Therefore, 
we have
\begin{align*}
\prod_{m\ge 0} \haT(q^{m} z)
&=  \exp \Big(-\sum_{r\geq1} \frac{1}{r} \haP_r z^{r}  \Big) \cdot
 \exp \Big(\sum_{m\ge 0} \sum_{k \geq1} \frac{1-q^{k}}{2k} [K_S]^k (q^{m}z)^{2k} \Big)
\\
&=  \exp \Big(-\sum_{r\geq1} \frac{1}{r} \haP_r z^{r}  \Big) \cdot
 \exp \Big(\frac12 \sum_{m\ge 0} \sum_{k \geq1} \frac{q^{2mk} -q^{(2m+1)k}}{k} [K_S]^k z^{2k} \Big)
\\
&=  \exp \Big(-\sum_{r\geq1} \frac{1}{r} \haP_r z^{r}  \Big) \cdot
 \exp \Big(\frac12  \sum_{m\ge 0} \ln \frac{1 -q^{2m+1} [K_S] z^2}{1 -q^{2m} [K_S] z^2} \Big)
\\
&=  \exp \Big(-\sum_{r\geq1} \frac{1}{r} \haP_r z^{r}  \Big) \cdot
\Big(\prod_{m\ge 0} \frac{1 -q^{2m+1} [K_S] z^2}{1 -q^{2m} [K_S] z^2} \Big)^{\frac12  }.
\end{align*}
Plugging this formula into \eqref{eq:EQK} and using Euler's identity \eqref{eq:Euler}, we obtain
\begin{align*}
\haE(- z) &= \exp \Big(- \sum_{r\geq1} \frac{1}{r} \haP_r z^{r}  \Big) \cdot
\Big(\prod_{n\ge 0}  (1 -q^{n} [K_S] z^2) \Big)^{-\frac12  }
\\
&
=  \exp \Big(- \sum_{r\geq1} \frac{1}{r} \haP_r z^{r}  \Big) \exp_q \Big(\frac{ [K_S] z^2}{1-q} \Big)^{\frac12}.
\end{align*}
This proves Proposition~\ref{prop:EHP}.
\end{proof}

\begin{corollary}
 \label{lem:PE}
We have
\begin{align*}
\haP(-z) \haE(z)
&=  \haE'(z)  - \haE(z) \sum_{k\ge 1} \frac {[K_S]^k z^{2k-1}}{ 1-q^k}.
\end{align*}
\end{corollary}

\begin{proof}
It follows by first replacing $\exp_q$ in the first formula in Proposition~\ref{prop:EHP} via the identity \eqref{exp=exp}, and then applying the differentiation $\frac{d}{dz}$.
\end{proof}






\end{document}